\documentclass[12pt]{article}

\usepackage{amsmath,amsthm,amscd,amsfonts,amssymb,longtable}
\usepackage{ifthen}
\usepackage{caption2}
\usepackage{graphicx}

\theoremstyle{remark}
\newtheorem{case}{Case}
\newtheorem{subcase}{Subcase}[case]
    \newtheorem*{subcase*}{Subcase}

\theoremstyle{plain}
\newtheorem{theorem}{Theorem}[section]
\newtheorem{corollary}[theorem]{Corollary}

\newtheorem{lemma}{Lemma}[section]
\newtheorem{proposition}[lemma]{Proposition}
\newtheorem{observation}[lemma]{Observation}

\theoremstyle{definition}
\newtheorem{example}{Example}[section]

\theoremstyle{definition}
\newtheorem*{remark}{Remark}

\theoremstyle{definition}
\newtheorem{definition}{Definition}[section]

\theoremstyle{definition}
\newtheorem{question}{Question}

\DeclareMathOperator{\mr}{mr}
\DeclareMathOperator{\hmr}{hmr}
\DeclareMathOperator{\Mm}{M}
\DeclareMathOperator{\hMm}{hM}
\DeclareMathOperator{\rank}{rank}
\DeclareMathOperator{\pin}{pin}
\DeclareMathOperator{\diag}{diag}
\DeclareMathOperator{\Col}{Col}
\DeclareMathOperator{\md}{MD}
\DeclareMathOperator{\MXi}{M\Xi}

\newcommand{\MD}[2]{{\md_{#1}(#2)}}

\newcommand{\newcaseset}{\setcounter{case}{0}}

\newcommand{\lm}{\lambda}
\newcommand{\eps}{\varepsilon}
\newcommand{\wh}{\widehat}

\newcommand{\cC}      {{\mathcal C}}
\newcommand{\cE}      {{\mathcal E}}
\newcommand{\cH}      {{\mathcal H}}
\newcommand{\cI}      {{\mathcal I}}
\newcommand{\cP}      {{\mathcal P}}
\newcommand{\cR}      {{\mathcal R}}
\newcommand{\cS}      {{\mathcal S}}

\newcommand{\RR}       { \mathbb{R}}
\newcommand{\CC}       { \mathbb{C}}
\newcommand{\NN}       { \mathbb{N}}

\newcommand{\Sym}{\cS}
\newcommand{\Herm}{\cH}
\newcommand{\hcI}{\mathrm{h}\cI}

\newcommand{\B}{X}
\newcommand{\R}{Y}
\newcommand{\BR}{bicolored }
\newcommand{\bl}{x}
\newcommand{\re}{y}
\newcommand{\nea}[1]{{#1}^\nearrow}
\newcommand{\uptrunc}[2]{\bigl[\nea{#2}\bigr]_{#1}}

\newcommand{\con}{+}
\newcommand{\del}{-}
\newcommand{\Econ}{\cE_v^\con}
\newcommand{\Edel}{\cE_v^\del}
\newcommand{\Ccon}{\cC_v^\con}
\newcommand{\Cdel}{\cC_v^\del}
\newcommand{\Qcon}{Q^\con}
\newcommand{\Qdel}{Q^\del}
\newcommand{\NNone}{\cE(K_1)}

\newcommand{\parti}{\pi}
\newcommand{\partboxes}[1]{\raisebox{-2pt}{\includegraphics{part#1}}}
\newcommand{\hpi}{\mathrm{h}\pi}

\newcommand{\vsum}{\underset{v}{\oplus}}
\newcommand{\bigvsum}{\bigoplus}
\newcommand{\FvsumGparen}{( F \vsum G )}
\newcommand{\caseskip}{\medskip\noindent}

\begin{document}

\baselineskip 18pt

\begin{titlepage}

\title{The Inverse Inertia Problem for Graphs}

\author
{Wayne Barrett\footnote {Wayne Barrett conducted this research as a
Fulbright grantee during Winter semester 2007 at the Technion. The
support of the US -- Israel Educational Foundation and of the
Department of Mathematics of the Technion is gratefully
acknowledged.} \footnote {Department of Mathematics, Brigham Young
University, Provo, Utah 84602, United States}, H.\,Tracy
Hall\footnote {Research conducted while a guest of the Department of
Mathematics at the Technion, whose hospitality is appreciated.}
\footnotemark[2], Raphael Loewy\footnote {Department of Mathematics,
Technion -- Israel Institute of Technology, Haifa 32000, Israel} }

\date{$17$ November 2007}

\maketitle

\begin{abstract}
Let $ G $ be an undirected graph on $ n $ vertices and let $ \Sym(G)
$ be the set of all real symmetric $ n\times n $ matrices whose
nonzero off-diagonal entries occur in exactly the positions
corresponding to the edges of $ G $. The inverse inertia problem for
$ G $ asks which inertias can be attained by a matrix in $ \Sym(G)
$. We give a complete answer to this question for trees in terms of
a new family of graph parameters, the maximal disconnection numbers
of a graph. We also give a formula for the inertia set of a graph
with a cut vertex in terms of inertia sets of proper subgraphs.
Finally, we give an example of a graph that is not inertia-balanced,
and investigate restrictions on the inertia set of any graph.
\end{abstract}



\thispagestyle{empty}
\end{titlepage}


\section{Introduction}

Given a simple undirected graph $ G=(V,E) $ with vertex set $
V=\{1,\ldots,n\} $, let $ \Sym(G) $ be the set of all real symmetric
$ n\times n $ matrices $ A = \bigl [ a_{ij} \bigr ]$ such that for $
i\ne j $, $ a_{ij}\ne 0 $ if and only if $ ij\in E $. There is no
condition on the diagonal entries of~$ A $.

The set $\Herm(G)$ is defined in the same way over Hermitian
$n \times n$ matrices,
and every problem we consider comes in
two flavors: the real version, involving $\Sym(G)$, and the
complex version, involving $\Herm(G)$.
There are known examples where a question of the sort we
examine here has a different answer when considered
over Hermitian matrices rather than over real symmetric
matrices~\cite{BvdHL1},~\cite{Hall}, but for each
question that is completely resolved in the
present paper, the answer over $\Herm(G)$ proves
to be the same as that obtained over $\Sym(G)$.

The inverse eigenvalue problem for graphs asks: Given a graph $ G $
on $ n $ vertices and prescribed real numbers $
\lm_1,\lm_2,\ldots,\lm_n $, is there some $ A \in \Sym(G) $ (or $A
\in \Herm(G)$, alternatively) such that the eigenvalues of $A$ are
exactly the numbers prescribed? In general, this is a very difficult
problem. Some contributions to its solution appear in \cite{DJ},
\cite{JD2}, \cite{JS2}, \cite{JDS}.

A more modest goal is to determine the maximum multiplicity $ \Mm(G) $ of an eigenvalue
of a matrix in $ \Sym(G) $.
This is easily seen to be equivalent to determining the minimum rank $ \mr(G) $
of a matrix in $ \Sym(G) $ since $ \mr(G)+\Mm(G)=n $.
This problem has been intensively studied.
Some of the major contributions appear in the papers
\cite{F}, \cite{N}, \cite{JD1}, \cite{Hs},  \cite{JS1}, \cite{vdH},  \cite{BFH1}, \cite{BvdHL1},
\cite{BvdHL2},  \cite{BFH2}, \cite{BFH3},  \cite{BF}, \cite{BGL},  \cite{JLS}, \cite{Hall}.
A variant of this problem is the study of $ \mr_+(G) $,
the minimum rank of all positive semidefinite $ A\in \Sym(G) $.
The Hermitian maximum multiplicity $\hMm(G)$,
Hermitian minimum rank $\hmr(G)$, and
Hermitian positive semidefinite rank
$\hmr_+(G)$ are defined analagously.

A problem whose level of difficulty lies between the inverse eigenvalue problem
and the minimum rank problem for graphs is the inverse inertia problem,
which we now explain.

\begin{definition}
\label{D1}
Given a Hermitian $ n\times n $ matrix $ A $, the \textit{inertia} of $ A $ is the triple
\[
\bigl(\pi(A),\nu(A),\delta(A)\bigr),
 \]
where
$ \pi(A) $ is the number of positive eigenvalues of $ A $,
$ \nu(A) $ is the number of negative eigenvalues of $ A $, and
$ \delta(A) $ is the multiplicity of the eigenvalue $ 0 $ of $ A $.
Then $ \pi(A)+\nu(A)+\delta(A)=n $ and $ \pi(A)+\nu(A)=\rank A $.
\end{definition}

If the order of $A$ is also known then the third number
of the triple is superfluous.
The following definition discards $\delta(A)$.

\begin{definition}
\label{D2}
Given a Hermitian matrix $ A $, the \textit{partial inertia} of $ A $ is the
ordered pair
\[
\bigl(\pi(A),\nu(A)\bigr).
 \]
We denote the partial inertia of $ A $ by $ \pin(A) $.
\end{definition}

We are interested in the following problem:
\begin{question}
[Inverse Inertia Problem]
Given a graph $G$ on $n$ vertices, for which ordered pairs $(r,s)$
of nonnegative integers with $r+s \le n$ is there a matrix
$A \in \Sym(G)$ such that $\pin(A) = (r,s)$?
\end{question}

The Hermitian Inverse Inertia Problem is the same question with
$\Herm(G)$ in the place of $\Sym(G)$. It is well known
\cite[p.\,8]{JD2} that in the case of a tree $T$ most questions over
$\Herm(T)$ are equivalent to questions over $\Sym(T)$, and in
particular if $F$ is a forest and $A \in \Herm(F)$, then there
exists a diagonal matrix $D$ with diagonal entries from the unit
circle such that $DAD^{-1} = DAD^* \in \Sym(F)$. In those sections
concerned with the inverse inertia problem for trees and forests we
thus assume without loss of generality that every matrix in
$\Herm(F)$ is in fact in $\Sym(F)$.

In this paper we give a complete solution to the inverse inertia problem
for trees and forests.
The statement of our solution is
a converse to an easier pair of lemmas
that apply not just to forests but to any graph.

\begin{lemma}
[Northeast Lemma] \label{L5} Let $G$ be a graph and suppose that $A
\in \Herm(G)$ with $\pin(A) = (\pi,\nu)$. Then for every pair of
integers $r \ge \pi$ and $s \ge \nu$ satisfying $r + s \le n$, there
exists a matrix $B \in \Herm(G)$ with $\pin(B) = (r,s)$. If in
addition $A$ is real, then $B$ can be taken to be real.
\end{lemma}
In other words, thinking of partial inertias or Hermitian partial
inertias as points in the Cartesian plane, the existence of a
partial inertia $(\pi,\nu)$ implies the existence of every partial
inertia $(r, s)$ anywhere ``northeast'' of $(\pi,\nu)$, as long as
$r + s$ does not exceed $n$. We prove this lemma in Section~\ref{S2}
by perturbing the diagonal entries of $A$.

To state the second lemma we need
to introduce an indexed family of graph parameters.
\begin{definition}
\label{DMD}
Let $G$ be a graph with $n$ vertices.
For any $k \in \{0, \ldots, n\}$
we define $\MD{k}{G}$,
the \textit{maximal disconnection} of $G$ by $k$ vertices,
as the maximum,
over all induced subgraphs $F$ of $G$
having $n-k$ vertices, of the number of
components of $F$.
\end{definition}
For example, $\MD{0}{G}$ is the number of components of $G$, and if
$T$ is a tree then $\MD{1}{T}$ is the maximum vertex degree of $T$.
Since an induced subgraph cannot have more components than vertices,
we always have $k + \MD{k}{G} \le n$.

\begin{remark}
As far as we can determine, $\MD{k}{G}$ is not a known family of
graph parameters. It is, however, related to the toughness $t$ of a
graph, which can be defined \cite{C} as
\[
  t(G) = \min \left \{ \frac{k}{\MD{k}{G}} :\ \MD{k}{G} \ge 2 \right \}.
\]
For a recent survey of results related to toughness of graphs,
see~\cite{BBS}. There is also some relation between $\MD{k}{G}$ and
vertex connectivity: a graph $G$ on $n$ vertices is $k$-connected,
$k<n$, if and only if $\MD{j}{G} = 1$ whenever $0 \le j < k$.
\end{remark}

\begin{lemma}
  [Stars and Stripes Lemma]
  \label{LMD}
  Let $G$ be a graph with $n$ vertices,
  let $k \in \{0, \ldots, n\}$ be such that
  $\MD{k}{G} \ge k$,
  and choose any pair of integers $r$ and $s$ such that
  $r \ge k$, $s \ge k$, and $r+s = n - \MD{k}{G} + k$.
  Then there exists a matrix $A \in \Sym(G)$
  such that $\pin(A) = (r,s)$.
\end{lemma}
This lemma is proved in Section~\ref{S2};
the idea of the proof is that
each partial inertia in the diagonal ``stripe'' from
$(n-\MD{k}{G},k)$ northwest to $(k,n-\MD{k}{G})$
can be obtained by combining
the adjacency matrices of ``stars'' at each
of the $k$ disconnection vertices
together with, for each of the remaining
components, a matrix
of co-rank $1$ and otherwise arbitrary inertia.

These two lemmas provide a partial solution to the
Inverse Inertia Problem for any graph.
Our main result for trees and forests is that for such graphs,
and exactly such graphs, the partial solution is complete.
\begin{definition}
  Let $G$ be a graph on $n$ vertices.
  Then $(r,s)$ is an \textit{elementary inertia} of $G$
  if for some integer $k$ in the range $0 \le k \le n$ we have
  $k \le r$, $k \le s$, and $n - \MD{k}{G} + k \le r+s \le n$.
\end{definition}
The elementary inertias of a graph $G$ are exactly those partial
inertias that can be obtained from $G$ by first applying the Stars
and Stripes Lemma and then applying the Northeast Lemma. The partial
solution given by these lemmas is the following: if $(r,s)$ is an
elementary inertia of a graph $G$, then there exists a matrix $A \in
\Sym(G)$ with $\pin(A) = (r,s)$.  This is proved as
Observation~\ref{Ob_elinI} in Section~\ref{S2}.

\begin{theorem}
\label{Th_treemain}
The Stars and Stripes Lemma and the Northeast Lemma
characterize the partial inertias of exactly forests, as follows:
\begin{enumerate}
\item
Let $F$ be a forest, and let $A \in \Sym(F)$ with $\pin(A) = (r,s)$.
Then $(r,s)$ is an elementary inertia of $F$.
\item
Conversely, let $G$ be a graph and suppose that for every $A \in
\Sym(G)$, $\pin(A)$ is an elementary inertia of $G$. Then $G$ is a
forest.
\end{enumerate}
\end{theorem}
Of course Claim 1 also applies for $A \in \Herm(F)$, since for $F$ a
forest any matrix in $\Herm(F)$ is diagonally congruent to a matrix
in $\Sym(F)$ having the same partial inertia.  Claim 2 of
Theorem~\ref{Th_treemain} is a corollary to known results, here
called Theorem~\ref{Th1}. We prove Claim 1 of
Theorem~\ref{Th_treemain} at the end of Section~\ref{S5}.

In Section~\ref{S4} we show that determining the set
of possible inertias of any graph with a cut vertex
can be reduced to the problem of determining the possible inertias
of graphs on a smaller number of vertices.
The formula we obtain is a generalization of the known formula for the
minimum rank of a graph with a cut vertex.
In Section~\ref{S5} we describe elementary inertias in terms of
certain edge-colorings of subgraphs, and we show that
the same cut-vertex formula proven in Section~\ref{S4} for inertias
also holds when applied to the (usually smaller) set of elementary inertias.
Applying these parallel formulas
inductively to trees and forests
then gives us a
proof of Claim 1 of Theorem~\ref{Th_treemain}.
In Section~\ref{S6} we outline an effective procedure for
calculating the set of partial inertias of any tree,
using the results of Section~\ref{S3} to justify
some simplications, and we calculate a few examples.
In Section~\ref{S7} we again consider more general graphs,
and demonstrate both an
infinite family of forbidden inertia patterns,
and the first example of
a graph that is not inertia-balanced.
The concept of an
inertia-balanced graph was introduced in~\cite{BF},
and determining whether a graph is inertia-balanced
is a special case of the inverse inertia problem.

\begin{definition}
\label{D3}
A Hermitian matrix $ A $ is \textit{inertia-balanced} if
\[
|\pi(A)-\nu(A)|\le 1.
 \]
A graph $ G  $ is \textit{inertia-balanced}
 if there is an inertia-balanced $ A\in \Sym(G) $ with $ \rank A=\mr(G) $.
A graph $ G  $ is \textit{Hermitian inertia-balanced}
 if there is an inertia-balanced $ A\in \Herm(G) $ with $ \rank A=\hmr(G) $.
\end{definition}

\begin{remark}
Our formulation, unlike the
definition in~\cite{BF}, is symmetric in allowing $ \nu(A)=\pi(A)+1 $.
This doubles the set of inertia-balanced matrices of odd rank,
but the two definitions are equivalent when applied to graphs
since $A \in \Sym(G)$ if and only if $-A \in \Sym(G)$.
\end{remark}

Barioli and Fallat \cite{BF} proved that every tree is
inertia-balanced.  Theorem~\ref{Th_treemain}, once proved, will
imply a slightly stronger result. The intuition for expecting a
graph to be inertia-balanced comes from many small examples in which
achieving an eigenvalue of high multiplicity appears to become
increasingly difficult as the imbalance increases between the number
of eigenvalues that are higher and the number that are lower than
the target multiple eigenvalue. The behavior observed in these small
examples can be stated formally in terms of the following
definitions.
\begin{definition}
  \label{Dstripe} A set $S$ of ordered pairs of integers is called
  \textit{symmetric} if whenever $(r,s) \in S$, then $(s,r) \in S$. A
  symmetric nonempty set $S$ of ordered pairs of nonnegative integers
  is called a \textit{stripe} if there is some integer $m$ such that
  $r+s = m$ for every $(r,s) \in S$, and we specify the particular
  constant sum by saying that $S$ is a \textit{stripe of rank $m$.}
  A stripe $S$ is \textit{convex} if the projection
  $\{ r :\ (r,s) \in S \}$ is a set of consecutive integers.
\end{definition}

\begin{example}
  The set $\{(2,2),(2,3),(3,2),(3,4),(4,3)\}$ is symmetric,
  the set $\{(6,0),(3,3),(0,6)\}$ is a stripe,
  and the stripe $\{(4,2),(3,3),(2,4)\}$
  is convex.
\end{example}

\begin{observation}
\label{Ob_symmetry}
Given a graph $G$ of order $n$
and an integer $m$ in the range
$\mr(G) \le m \le n$, the set
\[
\{ \pin(A) :\ A \in \Sym(G) \mbox{ \em{and} } \rank A = m \}
\]
is a stripe of rank $m$.
The same is true for $ A \in \Herm(G)$ with $m$
in the range $\hmr(G) \le m \le n$.
\end{observation}
\begin{proof}
Symmetry comes from the fact that $-A \in \Sym(G)$
if and only if $A \in \Sym(G)$,
and similarly for $\Herm(G)$.
The sets are
nonempty by
the definitions of $\mr(G)$ and $\hmr(G)$
and
the Northeast Lemma.
\end{proof}

\begin{definition}
A graph $G$ is
\textit{inertia-convex on stripes}
or
\textit{Hermitian inertia-convex on stripes}
if each of the stripes defined in
Observation~\ref{Ob_symmetry}
(with $A \in \Sym(G)$ or $A \in \Herm(G)$, respectively)
is convex.
\end{definition}

In other words, a graph is inertia-convex on stripes if each stripe
of possible partial inertias does not contain a gap.

\begin{corollary}
[Corollary to Theorem~\ref{Th_treemain}]
\label{Cor_stripes}
Every forest is inertia-convex on stripes.
\end{corollary}
\begin{proof}
Let $F$ be a forest. By Theorem~\ref{Th_treemain}, each of the
stripes defined in Observation~\ref{Ob_symmetry} is the set of
elementary inertias of some fixed rank $m$. For each fixed $k$ with
$\MD{k}{F} \ge k$ we obtain a set of elementary inertias which is a
union of convex stripes.  It follows that for any fixed $m$, the set
of elementary inertias of rank $m$ is the union of convex stripes of
rank $m$ as $k$ varies over all allowed integers. Since a union of
convex stripes of the same rank is a single convex stripe, each of
the stripes defined in Observation~\ref{Ob_symmetry} is convex.
\end{proof}

It has been an outstanding question if there is
any graph that is not inertia-balanced.
At the AIM Workshop in Palo Alto in October 2006,
the prevailing opinion
was that such a graph does not exist \cite{BHS}.

In Section~\ref{S7} we give an example of a graph
that is not inertia-balanced.
First we show that every graph satisfies
a condition that is
much weaker than inertia-balanced
(except in the case of minimum
semidefinite rank~$2$).
The counterexample graph and new condition together
allow us to completely determine
which sets can occur as the complement of the set
of possible partial inertias of a graph $G$
with $\mr_+(G) \le 3$.
The possible excluded partial inertia sets giving
minimum semidefinite
rank $4$ or greater remain unclassified.

For the most part our notation for graphs follows Diestel~\cite{D}.
We make use specifically of the following notation throughout:
\begin{itemize}
\item All graphs are simple,
and a graph is formally an ordered pair $G=(V,E)$
where $V$ is a finite set and $E$ consists of pairs from $V$.
When referring to an individual edge, we
abbreviate $\{u,v\}$ to $uv$ or $vu$.
The vertex set of a graph $G$ is also
referred to as $V(G)$, and the edge set as $E(G)$.
\item For $S \subseteq V$, $G[S]$ is the subgraph of $G$ induced by $S$ and
   $ G-S $ is the induced subgraph on $V(G) \setminus S$.  We write
   $G-F$ rather than $G-V(F)$ and $G-v$ rather than $G-\{v\}$.
\item The number of vertices of a graph $G$ is denoted $|G|$.
\item $ K_n $ is the complete graph on $ n $ vertices.
\item $ S_n = (\{1,2,3,\ldots,n\},\{12,13,\ldots,1n\})$ is called
the \textit{star graph} on $ n $ vertices.  This is the same as
the complete bipartite graph $K_{1,n-1}$.
\item $ P_n $ is the path on $ n $ vertices.  Paths are described explicitly
by concatenating the names of the vertices in order; for example, $uvw$ denotes
the graph $(\{u,v,w\},\{uv,vw\})$.
\item If $ v $ is a vertex of~$ G $, $ d(v) $ is the degree of~$ v $.
\item $ \Delta(G)=\max\left \{d(v) :\ v \in V(G)\right \} $.
\end{itemize}

We conclude the introduction with some elementary facts about inertia,
and include short proofs to keep the paper self-contained.

\begin{proposition}
\label{P1}
Let $ A $ be a Hermitian $ n\times n $ matrix and
let $ B $ be a principal submatrix of~$A$
of size $ (n-1)\times(n-1)$.
Then
\[
\pi(A)-1\le\pi(B)\le\pi(A) \ \ \textit{and} \ \ \nu(A)-1\le\nu(B)\le\nu(A).
 \]
\end{proposition}
\begin{proof}
By the interlacing inequalities \cite{B}
\[
\lm_1\ge\mu_1\ge\lm_2\ge\mu_2\ge\cdots\ge\lm_{n-1}\ge\mu_{n-1}\ge\lm_n,
 \]
where $ \lm_1,\ldots,\lm_n $ are the eigenvalues of $ A $
and $ \mu_1,\ldots,\mu_{n-1} $ are the eigenvalues of~$ B $,
arranged in decreasing order.  If $\mu_1 \le 0, \pi(B)=0 \le \pi(A)$.  Otherwise,
let $ m\in\{1,\ldots,n-1\} $ be the largest integer with $ \mu_m>0 $.
Then $ \lm_m>0 $ and $ \pi(B)=m\le\pi(A) $.

If $\lm_1 \le 0, \pi(B)=0 > \pi(A)-1$.  Otherwise, let $ \ell\in\{1,\ldots,n\} $ be the largest integer with $ \lm_\ell>0 $.
Then $ \mu_{\ell-1}\ge\lm_\ell>0 $ and $ \pi(B)\ge\ell-1=\pi(A)-1 $.

Similarly, $ \nu(A)-1\le\nu(B)\le\nu(A) $.
\end{proof}

\begin{proposition}
[Subadditivity]
\label{P2}
Let $ A $, $ B $ be Hermitian $ n\times n $ matrices and let $ C=A+B $.
Then
\[
 \pi(C)\le\pi(A)+\pi(B) \ \ \textit{and} \ \  \nu(C)\le\nu(A)+\nu(B).
\]
\end{proposition}
\begin{proof}
If $ \pi(A)+\pi(B)\ge n $, the first inequality is true, so assume
that $ \pi(A)+\pi(B)<n $. Let $ i=\pi(A) $ and $ j=\pi(B) $. Then $
\lm_{i+1}(A)\le 0 $ and $ \lm_{j+1}(B)\le 0 $. By the Weyl
inequalities \cite{B},
\[
\lm_{i+j+1}(C)= \lm_{i+1+j+1-1}(C)\le \lm_{i+1}(A)+\lm_{j+1}(B)\le 0.
\]
Therefore $ \pi(C)\le i+j=\pi(A)+\pi(B) $.

Similarly, $ \nu(C)\le\nu(A)+\nu(B) $.
\end{proof}

\begin{proposition}
\label{P3}
Let $ A $ be a Hermitian $ n\times n $ matrix and let $ cxx^* $ be a Hermitian rank~1 matrix
(so $c$ is real-valued).
Then
\[
\pi\left (A+cxx^*\right )\le
\left\lbrace
\begin{array}{ccc}
\pi(A)+1 & \textit{if} & c>0
\\
\pi(A)     & \textit{if} & c<0
\end{array}
\right.
 \]
and
\[
\nu\left (A+cxx^*\right )\le
\left\lbrace
\begin{array}{ccc}
\nu(A)+1 & \textit{if} & c<0
\\
\nu(A)     & \textit{if} & c>0
\end{array}
\right.
 \]
\end{proposition}
\begin{proof}
Let $ c>0 $.
Then $ \pi\left (cxx^*\right )=1 $, $ \nu\left (cxx^*\right )=0 $.
By Proposition~\ref{P2},
\begin{eqnarray*}
\pi\left (A+cxx^*\right )\le\pi(A)+ \pi\left (cxx^*\right )=\pi(A)+1,
\\
\nu\left (A+cxx^*\right )\le\nu(A)+ \nu\left (cxx^*\right )=\nu(A).
\end{eqnarray*}

The argument is similar if $ c<0 $.
\end{proof}


\section{The inertia set of a graph}
\label{S2}

\begin{definition}
\label{D4}
Let $ \NN $ be the set of nonnegative integers, and let $ \NN^2=\NN\times\NN $.
We define the following sets:
\begin{eqnarray*}
&& \NN^2_{\le k}=\left \{(r,s)\in\NN^2 :\ r+s\le k\right \},
\\
&& \NN^2_{\ge k}=\left \{(r,s)\in\NN^2 :\ r+s\ge k\right \},
\\
&& \NN^2_{[i,j]}=\NN^2_{\ge i}\cap\NN^2_{\le j},
\\
&& \NN^2_i=\NN^2_{[i,i]} \ \mbox{(the \textit{complete stripe} of rank $i$)}.
\end{eqnarray*}
\end{definition}
We note that a stripe of rank $i$ is
a nonempty symmetric subset of $\NN^2_i$.

\begin{definition}
\label{D5}
Given a graph $ G $, we define
\[
\cI(G)=\left \{(r,s) :\ \pin(A)=(r,s) \ \text{for some} \ A\in \Sym(G)\right \},
 \]
and
\[
\hcI(G)=\left \{(r,s) :\ \pin(A)=(r,s) \ \text{for some} \ A\in \Herm(G)\right \}.
 \]
We call $ \cI(G) $ the \textit{inertia set} of~$ G $
and
$ \hcI(G) $ the \textit{Hermitian inertia set} of~$ G $.
\end{definition}

Now suppose $ (r,s)\in\cI(G) $ and let $ A\in \Sym(G) $ with $ \pin(A)=(r,s) $.
Since $ r+s=\rank A $, we have $ \mr(G)\le r+s\le|G| $.
We record this as
\begin{observation}
\label{Ob4}
Given a graph $ G $ on $ n $ vertices,
$ \cI(G)\subseteq\NN^2_{[\mr(G),n]} $ and
$ \hcI(G)\subseteq\NN^2_{[\hmr(G),n]} $.
\end{observation}

The fact that every real symmetric matrix is also
Hermitian immediately gives us:
\begin{observation}
\label{Obvious}
For any graph $G$, $\cI(G)\subseteq\hcI(G)$ and $\hmr(G) \le \mr(G)$.
\end{observation}

The Northeast Lemma, as stated in the Introduction,
substantially shortens the calculation of the inertia set of a graph.
\begin{proof}
[Proof of Northeast Lemma] Let $ G $ be a graph and suppose that $
(\pi,\nu)\in\hcI(G) $, and let $ (r,s)\in\NN^2_{\le n} $ be given
with $ r\ge\pi $ and  $ s\ge\nu $. We wish to show that $
(r,s)\in\hcI(G) $. If in addition $ (\pi,\nu)\in\cI(G) $, we must
show that $ (r,s)\in\cI(G) $.

Let $ A\in \Herm(G) $ with $ \pin(A)=(\pi,\nu) $.
If $ \pi+\nu=n $ there is nothing to prove, so assume $ \pi+\nu<n $.
It suffices to prove that there exists a
$ B\in \Herm(G) $ with $ \pin(B)=(\pi+1,\nu) $,
because then an analogous argument
can be given to prove that there is a $ C\in \Herm(G) $
with $ \pin(C)=(\pi,\nu+1) $ and these two facts may be
applied successively to reach $ (r,s) $.
We also need to ensure that when $A$ is real
symmetric $B$ is also real symmetric.
Choose $ \eps>0 $ such that $ A+\eps I $ is invertible
and $ \nu(A+\eps I)=\nu $.
Then $ \pi(A+\eps I)=n-\nu $.
Let $A_0 = A$ and then perturb the diagonal
entries in order: for any $i \in \{1,\ldots,n\}$
let $A_i = A_{i-1}+\eps e_i^{\,}e_i^*$, so that
$A_n = A+\eps I$.
Then $ A_i\in \Herm(G) $ for $ i=0,\ldots,n $ and by Propositions~\ref{P2} and~\ref{P3},
\[ \pi(A_{i-1})\le\pi(A_i)\le\pi(A_{i-1})+1 \] for $ i=1,\ldots,n $.
It follows that
every integer in $\left\{\pi, \pi+1, \ldots, n-\nu \right\}$
is equal to $\pi(A_i)$ for some $i \in \left\{ 0, \ldots, n \right\}$.
Since
\[
 \nu=\nu(A+\eps I)\le\nu(A_{n-1})\le\cdots\le\nu(A_2)\le\nu(A_1)\le\nu(A)=\nu
\]
by Proposition~\ref{P3}, $ \nu(A_i)=\nu $ for $ i=0,1,\ldots,n $.
Then for some $ i $ we have $ \pin(A_i)=(\pi+1,\nu) $,
and we can take $B= A_i$.
As desired, $B$ is real symmetric if $A$ is real symmetric,
which completes the
$\Sym(G)$ version of the
Northeast Lemma
as well as the $\Herm(G)$
version:
Within either one of the two inertia sets
$\cI(G)$ or $\hcI(G)$,
the existence of a partial inertia $(\pi, \nu)$ implies
the existence of every partial inertia $(r,s)$
within the triangle
\[ r \ge \pi,\ s \ge \nu,\ r+s \le n, \]
or in other words every partial inertia
to the ``northeast'' of $(\pi, \nu)$.
\end{proof}

\begin{definition}
\label{D6}
If a graph $ G $ on $ n $ vertices satisfies $ \cI(G)=\NN^2_{[\mr(G),n]} $
we say that $ G $ is \textit{inertially arbitrary}.
If a graph $ G $ on $ n $ vertices satisfies $ \hcI(G)=\NN^2_{[\hmr(G),n]} $
we say that $ G $ is \textit{Hermitian inertially arbitrary}.
\end{definition}
\begin{example}
\label{Ex1}
The complete graph $K_n$, $n \ge 2$.
Since $ \pm J_n $ (the all ones matrix) $ \in \Sym(K_n) $, $ (1,0),(0,1)\in\cI(K_n) $.
By the Northeast Lemma $ \NN^2_{[1,n]}\subseteq\cI(K_n) $.
Since $ \cI(K_n)\subseteq\NN^2_{[\mr(K_n),n]}=\NN^2_{[1,n]} $ by Observation~\ref{Ob4},
$ K_n $ is inertially arbitrary.
\end{example}
\begin{example}

\label{Ex2}
The path $P_n$, $n \ge 2$.
A consequence of a well-known result of Fiedler~\cite{F}
is that for a graph $ G $ on $ n $ vertices, $ \mr(G)=n-1 $ if and only if $ G=P_n $.
It follows from a Theorem in~\cite{Hald} that there is an $ A\in \Sym(P_n) $
with eigenvalues $ 1,2,3,\ldots,n $.
Then for $ k=1,\ldots,n $, $ \pin(A-kI)=(n-k,k-1) $.
By the Northeast Lemma, $ \cI(P_n)=\NN^2_{[n-1,n]}=\NN^2_{[\mr(P_n),n]} $,
so $ P_n $ is also inertially arbitrary.
\end{example}

The partial inertia set for a graph on $ n $ vertices can never be smaller
than the partial inertia set for $ P_n $.
\begin{proposition}
\label{P6}
If $ G $ is any graph on $ n $ vertices, $\NN^2_{[n-1,n]} \subseteq \cI(G) $.
\end{proposition}
\begin{proof}
Let $ r,s\in\{0,1,\ldots,n\} $ with $ r+s=n $.
Let
\[
D=\diag(r,r-1,\ldots,2,1,-1,-2,\ldots,-s)
\]
and let $ A_G $ be the adjacency matrix of $ G $. By Gershgorin's
theorem, $ B=D+\tfrac{1}{2n}A_G\in \Sym(G) $ has eigenvalues $
\lm_1>\lm_2>\cdots>\lm_r>0>\lm_{r+1}>\cdots>\lm_n $, so $
\pin(B)=(r,s) $. Furthermore for $r<n$, $ B-\lm_{r+1}I_n\in \Sym(G)
$ has partial inertia $ (r,s-1) $. It follows that $\NN^2_{[n-1,n]}
\subseteq \cI(G) $.
\end{proof}
The fact that inertia sets are additive on disconnected unions of graphs
(Observation~\ref{Ob21})
gives us an immediate corollary.
\begin{corollary}
\label{Cor_ellcomp}
If $G$ is any graph on $n$ vertices and $G$ has $\ell$ components,
$\NN^2_{[n-\ell,n]} \subseteq \cI(G)$.
\end{corollary}

The existence of a complete stripe of
partial inertias of rank $n-\ell$
plays a role in the proof of
our second lemma from the Introduction.
\begin{proof}
  [Proof of the Stars and Stripes Lemma]
  Let $G$ be a graph with $n$ vertices,
  and let $S\subseteq V(G)$ be such that $|S| = k$
  and $G-S$ has $\MD{k}{G}$ components,
  with $\MD{k}{G} \ge k$.
  Also, let $(r,s)$ be any pair of integers such that
  $k \le r$, $k \le s$, and $r+s = n - \MD{k}{G} + k$.

Without loss of generality label the vertices of $G$ so that
$S = \{1, \ldots, k\}$,
and for each vertex $1 \le v \le k$
let $A_v$ be the $n \times n$
adjacency matrix of the subgraph of $G$
that retains all vertices of $G$,
but only those edges that include the vertex $v$.
If $v$ is isolated in $G$ then
$\pin(A_v) = (0,0)$;
otherwise the subgraph is
a star plus isolated vertices and
$\pin(A_v) = (1,1)$.

Now $G-S$ is a graph with  $n-k$ vertices and $\MD{k}{G}$
components, so by Corollary~\ref{Cor_ellcomp} there exists a matrix
$B \in \Sym(G-S)$ with $\pin(B) = (r-k, s-k)$. Let $C$ be the direct
sum of the $k \times k$ zero matrix with $B$, so that the rows and
columns of $C$ are indexed by the full set $V(G)$, as is the case
with the matrices $A_1, \ldots, A_k$. Let $M = A_1 + A_2 + \cdots +
A_k + C$. Then $M \in \Sym(G)$, and by subadditivity of partial
inertias (Proposition~\ref{P2}) we also have $\pi(M) \leq r-k+k = r$
and $\nu(M) \leq s-k+k = s$. Since $ r + s \le n$ we can apply the
Northeast Lemma to conclude that $(r,s) \in \cI(G)$.
\end{proof}

As mentioned in the introduction, the partial inertias which can be
deduced from Lemmas~\ref{L5} and~\ref{LMD} are precisely the
elementary inertias.
\begin{definition}
\label{D_elset} Let $G$ be a graph on $n$ vertices. Then the
\textit{set of elementary inertias of $G$}, $\cE(G)$, is given by
\[
  \cE(G) =
    \{ (r,s) \in \NN^2 :\
    (r, s)
    \ \mbox{ is an elementary inertia of $G$}
  \}.
\]
\end{definition}
  We may also think of $\cE(G)$ as follows:
  For each integer $k$, $0, \le k \le n$, let
  \[
    T_k = \{ (x,y) \in \RR^2 :\ k \le x, \ k \le y,\
      n-\MD{k}{G} + k \le x+y \le n \},
  \]
  and let $T = \bigcup\limits_{k=0}^n T_k$.
  Each nonempty $T_k$ is a possibly degenerate trapezoid, and
  \[
    \cE(G) = \NN^2_{\le n} \cap T.
  \]
\begin{observation}
\label{Ob_elinI}
  For any graph $G$, we have $\cE(G) \subseteq\cI(G)$.
\end{observation}
\begin{proof}
  Let $G$ be a graph on $n$ vertices, and suppose $(r,s) \in \cE(G)$.
  Then for some integer $k$ we have
  \[k \le r, \ k \le s, \mbox{ and } n - \MD{k}{G} + k \le r+s \le n.\]
  (Note that this implies $\MD{k}{G} \ge k$.)
  Recall that $k + \MD{k}{G} \le n$, so
  \[
    k+k \le n-\MD{k}{G} + k \le r+s.
  \]
  It follows that there is an
  ordered pair of integers $(x, y)$ satisfying
  \[
    k \le x \le r, \ k \le y \le s,
      \mbox{ and } x+y = n - \MD{k}{G} +k.
  \]
  The Stars and Stripes Lemma gives us
  $(x, y) \in \cI(G)$,
  after which the Northeast Lemma gives us
  $(r, s) \in \cI(G)$ since $r+s \le n$.
\end{proof}

\begin{remark}
  Given a graph $F$ on $m$ vertices there is a smallest integer $a$ such that
  $\NN^2_{a} \subseteq \cI(F)$.
  If $F$ is inertia-convex on stripes then $a$ is the same as $\mr_+(F)$,
  and if $F$ is inertially arbitrary then $a$ is the same as $\mr(F)$.
  Suppose that $F$ is $G-S$ as in the definition of $\MD{k}{G}$,
  with $|S| = k$ and $m = n-k$.
  Then some trapezoid of elementary inertias of $G$ comes from the
  easy estimate that the maximum
  co-rank of arbitrary inertia for $F$, i.e.~$m-a$, is at least the
  number of components $\ell$ of $F$ (Corollary~\ref{Cor_ellcomp}).
  Suppose we had an improved lower bound $\Xi(F)$ for this co-rank,
  a graph parameter that always
  satisfies $\ell \le \Xi(F) \le m-a$.
  (The improvement $\ell \le \Xi(F)$ will be guaranteed, for
  example, if $\Xi$ is additive on the components of $F$
  and is at least $1$ on each component.)
  We could then define a family of graph parameters analogous to
  $\MD{k}{G}$ by defining $\MXi_{k}(G)$ to be the maximum, over all
  subsets $S \subseteq V(G)$ of size $|S| = k$, of $\Xi(G-S)$.
  Replacing $\MD{k}{G}$ by
  $\MXi_{k}(G)$ would then give a stronger version of the Stars and
  Stripes Lemma, and an expanded set of not-as-elementary inertias.
\end{remark}
For any graph $G$, the Stars and Stripes Lemma gives
us a bound on the maximum eigenvalue multiplicity $\Mm(G)$.
\begin{corollary}
\label{Cor_Mbound}
Let $G$ be a graph on $n$ vertices.  Then for
any $0 \le k \le n$, $\Mm(G) \ge \MD{k}{G} - k$.
\end{corollary}
When this bound is attained,
it is attained in particular on a set that
includes the center of the stripe $\NN^2_{\mr(G)}$.
\begin{corollary}
\label{Cor_inertiabalanced}
  Let $G$ be a graph. If $\MD{k}{G} - k=\Mm(G)$ for some $k$,
  then $G$ is inertia-balanced.
\end{corollary}
\begin{example}
\label{ExHn} The \textit{$n$-sun} $H_n$ is defined as the graph on
$2n$ vertices obtained by attaching a pendant vertex to each vertex
of an $n$-cycle \cite{BFH1}. We have $\MD{0}{H_n} = 1$ and
$\MD{k}{H_n} = 2k$ for $1 \le k \le  \lfloor \frac n2 \rfloor$. It
follows that, in addition to $(2n-1,0)$ and $(0,2n-1)$, $\cI(H_n)$
contains every integer point $(r,s)$ within the trapezoid
\[
r+s \le 2n,\ 2n \le r + 2s,\ 2n \le 2r+s, 3n \le 2r + 2s.
\]
Since for $n > 3$ it is known that $\mr(H_n) = 2n - \lfloor \frac n2 \rfloor$
\cite{BFH1}, this shows that the $n$-sun is inertia-balanced for
$n > 3$.
\end{example}

It is useful to note the following connection
between the inverse inertia problem and the minimum semidefinite rank problem.

\begin{observation}
\label{Ob7}
The inertia set of a graph restricted to an axis gives
\begin{eqnarray*}
\cI(G)\cap\left (\NN\times\{0\}\right ) &=& \left \{(k,0) :\ k\in\NN, \mr_+(G)\le k\le n\right \},
\\
\cI(G)\cap\left (\{0\}\times\NN\right ) &=& \left \{(0,k) :\ k\in\NN, \mr_+(G)\le k\le n\right \},
\end{eqnarray*}
and similarly for $\hcI(G)$ and $\hmr_+(G)$.
\end{observation}
In other words, solving the inverse inertia problem for a graph $ G
$ on the $ x $-axis (or $ y $-axis) is equivalent to solving the
minimum semidefinite rank problem for $ G $. One well-known result
about minimum semidefinite rank is:
\begin{theorem}
[$\hmr_+$ \cite{vdH}, $\mr_+$ \cite{BFH3}]
\label{Th1}
Given a graph $ G $ on $ n $ vertices, $ \hmr_+(G)=n-1 $ if and only if $ G $ is a tree,
and  $ \mr_+(G)=n-1 $ if and only if $ G $ is a tree.
\end{theorem}

As noted in Example~\ref{Ex2}, if $G$ is not $P_n$ then $\mr(G) \ne n-1$,
and therefore $\mr(G) < n-1$.
It follows that $ \bigl\{P_n\bigr\}_{n=1}^\infty $ are the only inertially arbitrary trees.

If $G$ is not connected then any matrix in $\Sym(G)$ is a direct sum of smaller matrices,
which shows that $\mr_+$ is additive on the components of a graph.
\begin{observation}
\label{Obmainconverse}
Let $G$ be a graph on $n$ vertices and let $\ell$ be the number
of components of $G$.
Then $\mr_+(G) = n-\ell$ if and only if $G$ is a forest.
\end{observation}
This gives us a statement that implies the second claim of
Theorem~\ref{Th_treemain}.
\begin{corollary}
\label{Comainconverse} Let $G$ be a graph. Then $(\mr_+(G),0)$ is an
elementary inertia of $G$ if and only if $G$ is a forest.
\end{corollary}
\begin{proof}
Let $G$ be a graph on $n$ vertices and let $\ell$ be the number of
components of $G$. Since $\MD{0}{G} = \ell$, $(i,0)$ is an
elementary inertia of $G$ exactly for those integers $i$ in the
range $n-\ell \le i \le n$.  In particular, $(\mr_+(G),0)$ is an
elementary inertia if and only if $\mr_+(G) = n-\ell$.  By
Observation~\ref{Obmainconverse}, this is true if and only if $G$ is
a forest.
\end{proof}

Although the Stars and Stripes Lemma only gives the correct value
of $\mr_+(G)$ when $G$ is a forest, we have already seen that it
can give the correct values of $\mr(G)$ and $\Mm(G)$
for some graphs containing a cycle.
\begin{question}
What is the class of graphs for which
\[\Mm(G) = \max\limits_{0 \le k \le n} \left \{ \MD{k}{G} - k \right \}?\]
Theorem~\ref{Th_treemain} implies that this class includes all forests,
and Example~\ref{ExHn} shows that the class includes the $n$-sun
graphs $H_n$ for $n > 3$.
\end{question}

\begin{example}
\label{Ex3}
$ G=S_n $, $ n\ge 4 $.
Let $ A $ be the adjacency matrix of $ S_n $.
Then $ \pin(A)=(1,1) $.
Since $ \mr(S_n)=2 $ and $ \mr_+(S_n)=n-1 $, by the Northeast Lemma we have
\begin{eqnarray*}
\cI(S_n) &=& \bigl \{(n-1,0), (n,0), (0,n-1), (0,n)\bigr \}
\\
&& \cup \  \bigl \{(r,s) :\ r\ge1,\, s\ge 1,\, r+s\le n \bigr \}.
\end{eqnarray*}
It follows that $ S_n $ is not inertially arbitrary.
\end{example}

As has already been noted,
if $ A\in \Herm(G) $ with $ \pin(A)=(r,s) $,
then $ -A\in \Herm(G) $ with  $ \pin(A)=(s,r) $,
and if $A$ is real then $-A$ is real.
A consequence of Observation~\ref{Ob_symmetry} is
\begin{observation}
[Symmetry property] \label{Ob8} The sets $ \cI(G) $ and $\hcI(G)$
are symmetric about the line $ y=x $.
\end{observation}


\section{Principal parameters for the inertia set of a tree}
\label{S3}

The purpose of this section is to define some basic parameters
associated with a tree, establish their fundamental properties, and
relate them to the maximal disconnection numbers $\MD{k}{G}$. In
Section~\ref{S6} we will use these results to simplify the
application of Theorem~\ref{Th_treemain}.

In~\cite{JD1}, Johnson and Duarte computed the minimum rank of all matrices in $ \Sym(T) $,
where $ T $ is an arbitrary tree.
One of the graph parameters used by them, the \textit{path cover number} of~$ T $,
is also needed in our work.
It is defined as follows.
\begin{definition}
\label{D7}
Let $ T $ be a tree.
\begin{enumerate}
\item [(a)]
A \textit{path cover} of $ T $ is a collection of vertex disjoint paths,
occurring as induced subgraphs of~$ T $, that covers all the vertices of~$ T $.
\item [(b)]
The \textit{path cover number} of~$ T $, $ P(T) $,
is the minimum number of paths occurring in a path cover of~$ T $.
\item [(c)]
A \textit{path tree} $ \cP $ is a path cover of~$ T $ consisting of $ P(T) $ disjoint paths,
say $ Q_1, Q_2,\ldots,Q_{P(T)} $.
An \textit{extra edge} is an edge of~$ T $ that is incident to vertices on two distinct $ Q $'s.
Clearly there are exactly $ P(T)-1 $ extra edges.
\end{enumerate}
\end{definition}

The Theorem of Duarte and Johnson is
\begin{theorem}
\label{Th2}
For any tree $ T $ on $ n $ vertices,
\[
\mr(T)+P(T)=n.
 \]
\end{theorem}

As indicated, $ P(T) $ will also be used in our work.
This is not surprising, because inertia is a refinement of rank.
Our use of $ P(T) $ will be made precise now.
First, we need another definition.
\begin{definition}
\label{D8}
Let $ G=(V,E) $ be a graph, and let $ S\subseteq V $.
Let
\[
E_G(S)=\left \{xy\in E :\ x\in S \ \ \text{or} \ \ y\in S\right \},
 \]
that is, $ E_G(S) $ consists of all edges of $ G $ that are incident to at least one vertex in $ S $.
\end{definition}

We define now an integer-valued mapping $ f $ on the set of all subsets of $ V $ by:
\begin{equation*}
f_G(S)=|E_G(S)|-2|S|+1.
\end{equation*}
\begin{observation}
\label{Ob9}
For any graph $G$, $ f_G(\emptyset)=1 $.
\end{observation}

\begin{observation}
  \label{Ob9b}
  Let $T$ be a tree on $n$ vertices, and choose an
  integer $k$ in the range $0 \le k \le n$. Then
\begin{itemize}
\item
For every $S \subseteq V(T)$ with $|S| = k$,
$
f_T(S) \le \MD{k}{T} - k,
$
and
\item
For some $S \subseteq V(T)$ with $|S| = k$,
$
f_T(S) = \MD{k}{T} - k.
$
\end{itemize}
\end{observation}
\begin{proof}
In any forest, the number of components plus the number of edges
equals the number of vertices.  Let $S \subseteq V$ with $\left| S
\right| = k$. The forest $T-S$ has $n-1-E_T(S)$ edges and $n-k$
vertices, so it has $E_T(S)-k+1$ components. By definition of
$\MD{k}{T}$, $E_T(S) - k+1 \le \MD{k}{T}$, or equivalently, $f_T(S)
\le \MD{k}{T}-k$. Since $E_T(S) - k + 1 = \MD{k}{T}$ for some $S$
with $\left| S \right| = k$, the second statement follows.
\end{proof}

Our first theorem in this section is the following:
\begin{theorem}
\label{Th3}
Let $ T $ be a tree with $|T| = n$ and let $ P(T) $ denote its path cover number.
Then
\[
P(T)=\max\limits_{S\subseteq V} \left \{ f_T(S) \right \}
= \max\limits_{0 \le k \le n} \left \{\MD{k}{T} - k \right \}.
 \]
\end{theorem}

The second equality is a direct consequence of Observation~\ref{Ob9b}.
The first equality will be proved by induction on $ |T| $,
but first we prove it directly for several special cases.
These special cases will also be used in the proof of Theorem~\ref{Th3}.
\begin{observation}
\label{Ob10}
Theorem~\ref{Th3} holds for any path.
\end{observation}
\begin{proof}
Let $ P_n $ denote the path on $ n $ vertices, $n \ge 1$.
The degree of any vertex in $ P_n $ is at most two, so for any $ \emptyset\ne S\subseteq V(P_n) $,
$ f_{P_n}(S)\le 2|S|-2|S|+1=1 $.
The result  follows by Observation~\ref{Ob9}.
\end{proof}
\begin{corollary}
\label{Cor4}
Theorem~\ref{Th3} holds for any tree $ T $ with $ |T|\le 3 $.
\end{corollary}
\begin{observation}
\label{Ob11}
Theorem~\ref{Th3} holds for $ S_n $, for any $ n\ge 3 $.
\end{observation}
\begin{proof}
Label the pendant vertices of $ S_n $ by $ 2,3,\ldots,n $ and the vertex of degree $ n-1 $ by $ 1 $.
It is known that $ P(S_n)=n-2 $.
For $ S=\{1\} $, we have $ f_{S_n}(S)=n-1-2+1=n-2 $,
while for any $ S\subseteq V(S_n) $ it is straightforward to see that $ f_{S_n}(S)\le n-2 $.
\end{proof}
\begin{lemma}
\label{L12}
Let $ T $ be a tree, and let $ \emptyset\ne S\subseteq V(T) $.
Then $ f_T(S)\le P(T) $.
\end{lemma}
\begin{proof}
The proof is by induction on $ |T| $.
Corollary~\ref{Cor4} covers the base of the induction, so we proceed to the general induction step.

Let $ \cP $ be any path tree of $ T $, consisting of paths
$ Q_1,Q_2,\ldots,Q_{P(T)} $.
There are $ P(T)-1 $ extra edges.
We can assume without loss of generality that $ Q_1 $ is a pendant path in $ \cP $
(so exactly one extra edge emanates from it),
and we denote by $ v $ the vertex of $ Q_1 $ that is incident to an extra edge.

We can also assume without loss of generality that no vertex of $ S $
has degree $ 1 $ or $ 2 $, since deleting such a vertex cannot
increase the value of the function $f_T(S)$ that we are trying to
bound from above.
\begin{case}
$ v $ is an end vertex of $ Q_1 $.
In this case $ S\subseteq\bigcup\limits_{i=2}^{P(T)}V(Q_i) $.
Also, $ P(T-Q_1)=P(T)-1 $.
Applying the induction hypothesis, we get
\[
f_T(S)\le f_{T-Q_1}(S)+1\le P(T-Q_1)+1=P(T).
 \]
\end{case}
\begin{case}
$ v $ is an internal vertex of $ Q_1 $.
Suppose first that one of the two end vertices of $ Q_1 $ (call it $ z $)
is at distance (in $ Q_1 $)
of at least two from $ v $.
Then $ P(T-z)=P(T) $ and $ S\subseteq V(T-z) $.
Moreover, by the induction hypothesis,
\[
f_T(S)= f_{T-z}(S)\le P(T-z)= P(T).
 \]
Hence we may assume that $ Q_1 $ has the form $yvz$:
\[
\setlength{\unitlength}{1mm}
\begin{picture}(34,15)
 \put(6,7.7){\circle{5}}
 \put(5,6.6){$y$}
 \put(8.5,8){\line(1,0){6}}
 \put(17,7.8){\circle{5}}
 \put(16,6.5){$v$}
 \put(19.5,8){\line(1,0){6}}
 \put(28,7.8){\circle{5}}
 \put(27,6.5){$z$}
\end{picture}
\]
We have $ P(T-Q_1)=P(T)-1 $.
If $ v\notin S $ then, by induction,
\[
f_T(S)\le f_{T-Q_1}(S)+1\le P(T-Q_1)+1=P(T).
 \]
If $ v\in S $, then
\[
f_T(S)=f_{T-Q_1}\left (S\backslash\{v\}\right )+3-2\le P(T-Q_1)+1=P(T).
 \]
\end{case}
\end{proof}
\pagebreak[4]
\begin{observation}
\label{Ob13}
Let $ T $ be a tree that is not a star and for which $ \Delta(T)\ge 3 $.
Then there exists $ v\in V $ that has a unique non-pendant neighbor
and at least one pendant neighbor.
\end{observation}
\begin{proof}
Let $ r $ be a vertex of degree $ \Delta(T) $.
Let $ Q $ be a path starting at $ r $, and of maximum length.
Denote by $ u_1 $ the terminal vertex of $ Q $,
by $ v $ the predecessor of $ u_1 $ in $ Q $ (note that $ v\ne r $),
and by $ w $ the predecessor of $ v $ in $ Q $ (it is possible that $ w=r $).
Then $ u_1 $ is a pendant neighbor of $ v $,
and $ w $ is the unique non-pendant neighbor of $ v $.
\end{proof}
\begin{remark}
A similar result appears as Lemma~13 in~\cite{S}.
\end{remark}
\begin{proof}
[Proof of Theorem~\ref{Th3}]
As previously mentioned, the second equality comes from
Observation~\ref{Ob9b}.
The proof of the first equality is by induction on $ |T| $.
The base of the induction is ensured by Corollary~\ref{Cor4}.
The theorem holds for any path and any star,
by Observations~\ref{Ob10} and~\ref{Ob11}.
Hence we may assume that $ T $ is not a star and $ \Delta(T)\ge 3 $.
Let $ v $ be as in Observation~\ref{Ob13},
and let $ u_1,u_2,\ldots,u_m $ $ (m\ge 1) $ be its pendant neighbors.
\newcaseset
\begin{case}
$ m=1 $.
\\
In this case, $ P(T-u_1)=P(T) $.
By the induction hypothesis, there exists $ S\subseteq V(T-u_1) $
such that $ f_{T-u_1}(S)=P(T-u_1)=P(T) $.
Hence,
\[
f_T(S)\ge f_{T-u_1}(S)=P(T),
 \]
so $ f_T(S)=P(T) $ by Lemma~\ref{L12}.
\end{case}
\begin{case}
$ m=2 $.
\\
In this case it is straightforward to see that $ P\left (T-\{u_1,u_2,v\}\right )=P(T)-1 $.
By the induction hypothesis, there exists $ S\subseteq V(T-\{u_1,u_2,v\}) $
such that $ f_{T-\{u_1,u_2,v\}}(S)=P(T)-1 $.
Hence, for $ S_1=S\cup\{v\} $
we have
\[
f_T(S_1)=f_{T-\{u_1,u_2,v\}}(S)+3-2=P(T).
 \]
\end{case}
\begin{case}
$ m\ge 3 $.
\\
In this case it is straightforward to see that $ P(T-u_m)=P(T)-1 $.
By the induction hypothesis, there exists $ S\subseteq V(T-u_m) $
such that $ f_{T-u_m}(S)=P(T-u_m)=P(T)-1 $.
If $ v\in S $ then $ f_T(S)=f_{T-u_m}(S)+1=P(T) $,
so we may assume that $ v\notin S $.
We claim that $ u_1,u_2,\ldots,u_{m-1}\notin S $.
Suppose otherwise that $ u_1\in S $.
Then
\[
f_{T-u_m}\left (S\backslash\{u_1\}\right )=f_{T-u_m}(S)+2-1>P(T-u_m),
 \]
contradicting Lemma~\ref{L12}.
Let $ S_1=S\cup\{v\} $. Then
\[
f_T(S_1)\ge f_T(S)+m-2\ge f_{T-u_m}(S)+1=P(T),
 \]
so $ f_T(S_1)=P(T) $ by Lemma~\ref{L12}.
\end{case}

This completes the proof that
\[
P(T)=\max\limits_{S\subseteq V} \left \{ f_T(S) \right \}
= \max\limits_{0 \le k \le n} \left \{\MD{k}{T} - k \right \}.
 \]
\end{proof}
We pause to note a similar result to Theorem~\ref{Th3}.
Given a tree $T$, Johnson and Duarte~\cite{JD1} ascertained
that $P(T)$ is the maximum of $p-q$ such that there exist
$q$ vertices whose deletion leaves $p$ components
each of which is a path (possibly including singleton paths).
It is obvious that this maximum is at most
$ \max\limits_{0 \le k \le n} \left \{ \MD{k}{T} - k \right \} $
since any components are allowed in determining $\MD{k}{T}$,
and the converse is also true:
if any component of the remaining forest is not a path, then
deleting a vertex of degree greater than $2$ increases
the value of $\MD{k}{T} - k$.
The observation of Johnson and Duarte can thus be seen
as a corollary of Theorem~\ref{Th3}.
We will see the usefulness of allowing non-path components
in Section~\ref{S6},
where we show that the lower values of
$\MD{k}{T}$ provide an exact description of part of the
boundary of $\cI(T)$.
\begin{definition}
\label{D9} Let $ T $ be a tree.
\begin{enumerate}
\item [(a)]
A set $ S\subseteq V(T) $ is said to be \textit{optimal} if $ f_T(S)=P(T) $.
\item [(b)]
Let $ c(T)=\min\left \{|S| :\ S~\text{is optimal}\right \} $.
\item [(c)]
We say $ S $ is \textit{minimal optimal} if $ S $ is optimal and $ |S|=c(T) $.
\end{enumerate}
\end{definition}
\begin{observation}
\label{Ob_cMD}
For a tree $T$,
\[
  c(T) = \min\limits_{0 \le k \le n} \left \{ k :\
      n - \MD{k}{T} + k = \mr(T) \right \}.
\]
\end{observation}
\begin{proof}
Theorem~\ref{Th2}, Theorem~\ref{Th3}, Definition~\ref{D9}, and Observation~\ref{Ob9b}.
\end{proof}
\begin{observation}
\label{Ob_c_mr}
For a tree $T$,
\[
  c(T) \le \left \lfloor \frac{\mr (T)}{2} \right \rfloor .
\]
\end{observation}
\begin{proof}
  Let $h = c(T)$ so $n-\MD{h}{T}+h = \mr (T)$.
  Recall that $k+\MD{k}{T} \le n$ for any integer
  $k$, $0 \le k \le n$, so in particular $h \le n-\MD{h}{T}$
  and therefore $2h \le \mr(T)$.
\end{proof}
\begin{observation}
\label{Ob14}
Let $ T $ be a tree and let $ S\subseteq V(T) $ be minimal optimal.
Then $ d(v)\ge 3 $ for every $ v\in S $.
\end{observation}
\begin{example}
\label{Ex4}
We calculate $c(T)$ for paths and stars.
\begin{itemize}
\item
$T=P_n$:  Then $f_{P_n}(\emptyset)=1=P(P_n)$ so $c(P_n)=0$.
\item
$T=S_n, n \ge 4$: Let $v$ be the degree $n-1$ vertex of $S_n$.  Then
$f_{S_n}(\{v\})=n-1-2+1=P(S_n)>1=f_{S_n}(\emptyset)$.  So $c(S_n)=1$.
\end{itemize}
\end{example}
\begin{proposition}
\label{P15}
Let $ T $ be a tree and let $ v\in V(T) $ be adjacent to $ m\ge 2 $ pendant vertices $ u_1,u_2,\ldots,u_m $,
and at most one non-pendant vertex $ w $.
Then there is a path tree $ \cP $ of $ T $ in which
$ u_1vu_2 \in\cP  $.
\end{proposition}
\begin{proof}
The claim is obvious if $ T $ is a star, $ |T|\ge 3 $, so assume this is not the case.
Then $ v $ is adjacent to exactly one non-pendant vertex.

Let $ \cP $ be a path tree of $ T $.
Then at least $ m-2 $ of the vertices $ u_1,u_2,\ldots,u_m $ give single-vertex paths in $ \cP $.
Let $ Q $ be a path in $ \cP $ containing $ v $.
Then $ Q $ contains at least one pendant neighbor of $ v $, say $ u_1 $.
Then $ Q= u_1vv_1v_2\ldots v_k  $.
Note that $ k\ge 1 $, as $ \cP $ is a path tree.
If $ v_1=u_2 $, then $ Q=u_1vu_2 $.
Otherwise, $ u_2 $ is a single-vertex path in $ \cP $.
We can form a new path tree $ \cP_1 $
by replacing the path $ u_1vv_1v_2\ldots v_k $ and the singleton path $ u_2 $ of $ \cP $
by the pair of paths $ u_1vu_2 $ and $ v_1v_2\ldots v_k $.
\end{proof}
\begin{proposition}
\label{P16}
Let $ T $ be a tree and let $ v\in V(T) $ be adjacent to $ m\ge 2 $ pendant vertices $ u_1,u_2,\ldots,u_m $,
and at most one non-pendant vertex $ w $.
Let $ T_1=T-\left \{u_1,u_2,\ldots,u_m,v\right \} $.
Then
\[
P(T)=P(T_1)+m-1,
 \]
and
\[
c(T)\le c(T_1)+1.
 \]
If $ m\ge 3 $,
\[
c(T)=c(T_1)+1.
 \]
\end{proposition}
\begin{proof}
The proposition clearly holds if $ T $ is a star,
so we may assume that this is not the case.
Let $ \cP_1 $ be a path tree for $ T_1 $.
Then $ \cP_1\cup u_1vu_2 \cup u_3 \cup\cdots\cup u_m  $ is a path cover for $ T $,
so $ P(T)\le P(T_1)+m-1 $.

Now let $ \cP $ be a path tree for $ T $ containing the path $ u_1vu_2 $ (see Proposition~\ref{P15}).
Then $ \cR=\left \{u_1vu_2,u_3,\ldots,u_m\right \}\subseteq\cP $, and $ \cP\backslash\cR $
is a path cover for $ T_1 $.
Therefore,
\[
P(T_1)\le\left |\cP\backslash\cR\right |=|\cP|-(m-1)=P(T)-(m-1).
 \]
Hence $ P(T)=P(T_1)+m-1 $.

Now let $ S $ be a minimal optimal set for $ T_1 $, so $ |S|=c(T_1) $.
This implies that $ \left |E_{T_1}(S)\right |-2|S|+1=P(T_1) $.
Let $ S_v=S\cup\{v\} $.
Since $ T $ is not a star $ v $ has a unique non-pendant neighbor $ w $.
The vertices $ w,u_1,u_2,\ldots,u_m $ are adjacent
to $ v $, so $ \left |E_T(S_v)\right |=\left |E_{T_1}(S)\right |+m+1 $.
Then
\begin{eqnarray*}
f_T(S_v) &=& \left |E_T(S_v)\right |-2|S_v|+1=\left |E_{T_1}(S)\right |+m+1-2|S|-2+1
\\
&=& P(T_1)+m-1=P(T),
\end{eqnarray*}
so $ S_v $ is an optimal set for $ T $. It follows that
\[
c(T)\le|S_v|=|S|+1=c(T_1)+1.
 \]
Now assume that $ m\ge 3 $ and that $ S $ is a minimal optimal set for $ T $.
By Observation~\ref{Ob14}, none of the vertices $ u_1,u_2,\ldots,u_m $ is in $S$.
If $ v\notin S $, Lemma~\ref{L12} implies
\begin{eqnarray*}
P(T) &\ge & \left |E_T\left (S\cup\{v\}\right )\right |-2\left |S\cup\{v\}\right |+1
\\
&\ge & \left |E_T\left (S\right )\right |+3-2|S|-2+1=P(T)+1,
\end{eqnarray*}
a contradiction. Therefore $ v\in S $.

Let $ S^\prime=S\backslash\{v\} $.
Then $ \left |E_{T_1}\left (S^\prime\right )\right |\ge\left |E_T\left (S\right )\right |-(m+1) $,
so
\begin{eqnarray*}
f_{T_1}(S^\prime) &\ge &\left |E_T\left (S\right )\right |-(m+1)-2\left |S^\prime\right |+1
\\
&=& \left |E_T\left (S\right )\right |-(m+1)-2|S|+2+1
\\
&=& P(T)-(m-1)=P(T_1).
\end{eqnarray*}
It follows from Lemma~\ref{L12} that $ S^\prime $ is optimal for $
T_1 $, implying
\[
c(T_1)\le\left |S^\prime\right |=|S|-1=c(T)-1.
 \]
Hence $ c(T)=c(T_1)+1 $.
\end{proof}

Proposition~\ref{P16} gives us a simple algorithm to calculate
$P(T)$ and thus the minimum rank of a tree.
We will use the fact that if $u$ is a pendant vertex whose neighbor
$v$ has degree $2$, then any path in a minimal path cover that
includes the vertex $v$ will also include the vertex $u$,
and $P(T)$ = $P(T-u)$.

\begin{observation}
\label{Ob_Palgorithm}
Let $T$ be a tree.  Then $P(T)$ may be calculated as follows:
\begin{enumerate}
\item
Set $G$ to $T$ and set $p$ to $0$.
\item If $G$ has a pendant vertex $u$ whose neighbor $v$ has degree $2$,
then replace $G$ by $G - u$ and repeat step 2.
\item If $G$ consists of a single edge or single vertex, then $P(T) = p+1$.
If $G$ is a star on $m+1$ vertices, then $P(T) = p+m-1$.
\item
In all other cases (by Observation \ref{Ob13}) there will be some $v
\in V(G)$ that is adjacent to $m \ge 2$ pendant vertices $u_1, u_2,
\ldots, u_m$ and exactly one non-pendant vertex $w$. Replace $G$ by
$G - \{u_1, u_2, \ldots, u_m, v \}$, replace $p$ by $p + m - 1$, and
return to step 2.
\end{enumerate}
\end{observation}

The calculation of $c(T)$ is not quite as straightforward as that of
$P(T)$, although we can show one special case in which it is
additive on subgraphs. For this we need the following definition.

\begin{definition}
\label{D14} Let $ F $ and $ G $ be graphs on at least two vertices,
each with a vertex labeled $ v $. Then $ F\vsum G $ is the graph on
$ |F|+|G|-1 $ vertices obtained by identifying the vertex $ v $ in $
F $ with the vertex $ v $ in $ G $.
\end{definition}

The vertex $ v $ in Definition~\ref{D14} is commonly referred to as a
\textit{cut vertex} of the graph $ F\vsum G $.
The next result determines $ c(T) $ when $ d(v)=2 $.

\begin{theorem}
\label{Th5}
Let $ T_1 $ and $ T_2 $ be trees each with a pendant vertex labeled $ v $.
Let $ T=T_1\vsum T_2 $.
Then $ c(T)=c(T_1)+c(T_2) $.
\end{theorem}
\begin{proof}
Let $ R_1 $, $ R_2 $ be minimal optimal sets for $ T_1 $, $ T_2 $, respectively.
Then
\[
f_{T_i}(R_i)=
\left |E_{T_i}(R_i)\right |-2|R_i|+1=P(T_i), \ \ \ i=1,2.
 \]
Since $ P(T)\le P(T_1)+P(T_2)-1 $, and $ v\notin R_1,R_2 $
by Observation~\ref{Ob14}, by Lemma~\ref{L12}
\begin{eqnarray*}
P(T)\ge f_T(R_1\cup R_2)
&=&\left |E_T\left (R_1\cup R_2\right )\right |-2\left |R_1\cup R_2\right |+1
\\
&=& \sum\limits_{i=1}^2\bigl (\left |E_{T_i}(R_i)\right |-2|R_i|+1\bigr )-1
\\
&=& P(T_1)+P(T_2)-1\ge P(T).
\end{eqnarray*}
Therefore, $ R_1\cup R_2 $ is an optimal set for $ T $ by
Lemma~\ref{L12} and $ c(T)\le\left |R_1\cup R_2\right |=
|R_1|+|R_2|=c(T_1)+c(T_2) $. We also see that $ P(T)=P(T_1)+P(T_2)-1
$.

Suppose now $ S $ is a minimal optimal set for $ T $.
By Observation~\ref{Ob14}, $ v\notin S $.
Let $ S_i=S\cap V(T_i) $, $ i=1,2 $.
Since $ v\notin S $, $ S_1\cap S_2=\emptyset $.
Now
\begin{eqnarray*}
P(T) &=& f_T(S)=\left |E_T(S)\right |-2|S|+1,
\\
P(T_1) &\ge & f_{T_1}(S_1)=\left |E_{T_1}(S_1)\right |-2|S_1|+1,
 \\
\text{and} \qquad\qquad &&
\\
P(T_2) &\ge & f_{T_2}(S_2)=\left |E_{T_2}(S_2)\right |-2|S_2|+1.
\end{eqnarray*}
Then
\begin{eqnarray*}
1 &=& P(T_1)+P(T_2)-P(T)
\\
&\ge & \sum\limits_{i=1}^2\left (\left |E_{T_i}(S_i)\right |-2|S_i|+1\right )
-\left (\left |E_T(S)\right |-2|S|+1\right )=1,
\end{eqnarray*}
so we must have
\[
P(T_i)=f_{T_i}(S_i), \ \ \ i=1,2,
 \]
and
\[
c(T_1)+c(T_2)\le|S_1|+|S_2|=|S|=c(T).
 \]
\end{proof}
\begin{corollary}
\label{Cor6}
Let $ p $ be a pendant vertex in a tree $ T $ and suppose the neighbor of~$ p $ has degree $ 2 $.
Then $ c(T)=c(T-p) $.
\end{corollary}
\begin{corollary}
\label{Cor7}
If a tree $ T $ has exactly one vertex of degree $ d>2 $, then \mbox{$ c(T)=1 $}.
\end{corollary}
\begin{proof}
It is straightforward to see, by repeated application of Corollary~\ref{Cor6},
that $ c(T)=c(S_{d+1})=1 $.
\end{proof}
\begin{definition}
\label{D10}
Let $T$ be a tree and let $ k $ be an integer such that $ 0\le k\le c(T) $.
Then
\[
r_k(T)=\max\left \{|E_T(S)| :\ S\subseteq V(T),\,  |S|=k\right \}.
 \]
\end{definition}
\begin{observation}
\label{Obr_k}
For a tree $T$, $r_k(T) = \MD{k}{T} + k - 1$.
\end{observation}

The next theorem will play an important role in
simplifying the computation of
$\cI(T)$.
\begin{theorem}
\label{Th8}
Let $ T $ be a tree with $ c(T)\ge 1 $.
Then
\[
r_k(T)-r_{k-1}(T)\ge\left \{
\begin{array}{ccc}
3 & \ \ \textit{if} \ \ & k=1 \ \textit{or} \ k=c(T),
\\ \\
2 & \ \ \textit{if} \ \ & 1<k<c(T).
\end{array}
\right .
 \]
\end{theorem}
\begin{proof}
Since $ c(T)\ge 1 $, $ T $ is not a path.
Therefore $ \Delta(T)\ge 3 $, implying $ r_1(T)-r_0(T)=\Delta(T)-0\ge 3 $.
If $ k=c(T) $,
\[
r_k(T)-2k+1=P(T),
 \]
while
\[
r_{k-1}(T)-2(k-1)+1<P(T).
 \]
Then
\[
r_k(T)-r_{k-1}(T)-2\ge 1.
 \]
Thus, the stronger conclusion in the special cases
$ k=1 $ and $ k=c(T) $ has been established.
We proceed by induction on $ |T| $.
Since $ T $ cannot be a path, the base of the induction is $ |T|=4 $,
and the only relevant tree $ T $ with $ |T|=4 $ is $ S_4 $.
Since $c(S_4)=1 $ the theorem holds in this case.

Consider now the general induction step.
Let $ T $ be a tree on $ n $ vertices,
and let $ k\in\{2,\ldots,c(T)-1\} $.
Note that if $ c(T)\le 2 $ we are done.
In particular, we can assume $ T $ is not a star.
We have to show that $ r_k(T)-r_{k-1}(T)\ge 2 $.

By Observation~\ref{Ob13} there exists $ v\in V $
that is adjacent to a unique non-pendant vertex $ w $,
and to pendant vertices $ u_1,u_2,\ldots,u_m $, where $ m\ge 1 $.
\newcaseset
\begin{case}
$ m\ge 2 $.
\\
Let $ T_1=T-\{u_1,u_2,\ldots,u_m,v\} $.
Then $ c(T_1)\ge c(T)-1 $ by Proposition~\ref{P16}.
This tells us both that we are allowed to assume
the induction hypothesis on the tree $T_1$
(which requires $ c(T_1) \ge 1 $) and that
$ k\le c(T_1) $.

Now choose $ Q\subseteq V(T) $ with $ |Q| \le k-1 $ and $ |E_T(Q)| \ge r_{k-1}(T) $.
This choice is possible (with equality) by
the definition of $r_{k-1}$.  We can assume without loss of
generality that $Q$ contains none of the vertices $\{u_1, \ldots, u_m\}$
as follows:
If $ v\in Q $ we delete all $ u_i $'s that belong to $ Q $, possibly
decreasing $|Q|$ without changing $|E_T(Q)|$.
If at least one of $ u_1,u_2,\ldots,u_m $, say $ u_1 $, belongs to $ Q $ but $ v\notin Q $,
we replace $ u_1 $ by $ v $ in $ Q $ and delete from $ Q $ all remaining $ u_i $,
possibly decreasing $|Q|$ and possibly increasing $|E_T(Q)|$.
We give the name $\ell$ to $|Q|$, so $\ell \le k-1$.
\begin{subcase}
Suppose that $ v\notin Q $.
Let $ R=Q\cup\{v\} $.
Then $ |R|=\ell+1\le k $,
and $ E_T(R)\supseteq E_T(Q)\cup\bigl \{ vu_1, vu_2,\ldots,vu_m\bigr \} $.
Hence
\[
r_k(T)\ge r_{\ell+1}(T)\ge|E_T(R)|\ge|E_T(Q)|+m\ge r_{k-1}(T)+m\ge r_{k-1}(T)+2.
 \]
\end{subcase}
\begin{subcase}
Suppose that $ v\in Q $.
Let $ Q^\prime=Q\backslash\{v\} $.
Then $ |Q^\prime|=\ell-1\le k-2 $, and $ r_{\ell-1}(T_1)\ge|E_{T_1}(Q^\prime)| $.
By the induction hypothesis,
\[
r_{k-1}(T_1)-r_{\ell-1}(T_1)\ge r_{k-1}(T_1)-r_{k-2}(T_1)\ge 2.
 \]
Choose $ R\subseteq V(T_1) $ with $ |R|=k-1 $
such that $ r_{k-1}(T_1)=|E_{T_1}(R)| $.
Let $ R_v=R\cup\{v\} $.
Then $ |R_v|=k $, and
\[
E_T(R_v)=E_{T_1}(R)\cup\bigl\{vu_1, vu_2,\ldots,vu_m,vw\bigr \}.
 \]
Also,
\[
E_T(Q)=E_{T_1}(Q^\prime)\cup\bigl\{vu_1, vu_2,\ldots,vu_m,vw\bigr \},
 \]
so
\[
|E_T(R_v)|=|E_{T_1}(R)|+m+1, \ \ \ |E_T(Q)|=|E_{T_1}(Q^\prime)|+m+1.
 \]
Then
\begin{eqnarray*}
r_k(T) &\ge & |E_T(R_v)|=|E_{T_1}(R)|+m+1=r_{k-1}(T_1)+m+1
\\
&\ge & r_{\ell-1}(T_1)+2+m+1
\\
&\ge & |E_{T_1}(Q^\prime)|+m+1+2=|E_T(Q)|+2\ge r_{k-1}(T)+2.
\end{eqnarray*}
\end{subcase}
\end{case}
\begin{case}
$ m=1 $.
\\
Let $ T_1=T-u_1 $.
By Corollary~\ref{Cor6}, we have $ c(T)=c(T_1) $, and clearly $ P(T)=P(T_1) $.
Since $ k\le c(T)-1 $, $ k\le c(T_1)-1 $.

As in Case~1, we choose $Q\subseteq V(T)$ with $|Q| \le k-1$
and $|E_T(Q)| \ge r_{k-1}(T)$, and can assume without
loss of generality that $u_1 \notin Q$.
\begin{subcase}
Suppose that $ v\notin Q $.
Then $ |E_{T_1}(Q)|=|E_T(Q)| $,
and applying the induction hypothesis, we have
\[
r_k(T)\ge r_k(T_1)\ge r_{k-1}(T_1)+2\ge|E_{T_1}(Q)|+2=|E_T(Q)|+2 \ge r_{k-1}(T)+2.
 \]
\end{subcase}
\begin{subcase}
Suppose that $ v\in Q $. Let $ Q^\prime=Q\backslash\{v\} $. Then $
|Q^\prime| \le k-2 $ and $ E_{T_1}(Q^\prime)=E_T(Q^\prime) $. Hence
\[
|E_{T_1}(Q^\prime)|=|E_T(Q^\prime)|\ge|E_T(Q)|-2 \ge r_{k-1}(T)-2.
 \]
Applying the induction hypothesis to $ T_1 $,
\[
r_{k-1}(T_1)\ge r_{k-2}(T_1)+2\ge|E_{T_1}(Q^\prime)|+2\ge r_{k-1}(T).
 \]
Applying the induction hypothesis again to $ T_1 $,
\[
r_{k}(T)\ge r_{k}(T_1)\ge r_{k-1}(T_1)+2 \ge r_{k-1}(T)+2.
 \]
\end{subcase}
\end{case}
\end{proof}
\begin{corollary}
\label{Cor_MDmonotone}
Let $T$ be a tree.  Then for $0 \leq j \leq k \leq c(T)$,
\[
\MD{k}{T} \ge \MD{j}{T} + (k-j).
\]
\end{corollary}
\begin{proof}
It suffices to prove the case $k-j = 1$.
Here we have $1 \leq k \leq c(T)$, and
Theorem~\ref{Th8} gives us
$
r_k(T) \ge r_{k-1}(T)+2
$.
Making the substitution
$r_{k}(T) = \MD{k}{T} + k - 1$
from Observation~\ref{Obr_k}
gives us the desired result.
\end{proof}
\begin{proposition}
\label{P17}
Let $ T $ be a tree on $ n\ge 3 $ vertices.
Then $ c(T)\le\dfrac{n-1}{3} $.
\end{proposition}
\begin{proof}
We prove the proposition by induction on $ n $.
The cases $ n=3, 4 $ are obvious, so we consider the general induction step.
The proposition holds if $ T=S_n $, as $ c(S_n)=1 $, so assume $ T $ is not a star.
By Observation~\ref{Ob13}, there exists a vertex $ v $
that is adjacent to exactly one non-pendant vertex $ w $
and to pendant vertices $ u_1,u_2,\ldots,u_m $, where $ m\ge 1 $.
\newcaseset
\begin{case}
$ m=1 $.
\\
It follows from Corollary~\ref{Cor6} and the induction hypothesis that
\[
c(T)=c(T-u_1)\le\dfrac{n-1-1}{3}<\dfrac{n-1}{3}.
 \]
\end{case}
\begin{case} $ m\ge 2 $.
\\
Let $ T_1=T-\{u_1,u_2,\ldots,u_m,v\} $.
Then $ |T_1|\le n-3 $, and by Proposition~\ref{P16} $ c(T)\le c(T_1)+1 $.
Then, by induction hypothesis,
\[
c(T)-1\le c(T_1)\le\dfrac{|T_1|-1}{3}\le\dfrac{n-4}{3}=\dfrac{n-1}{3}-1.
 \]
\end{case}
\end{proof}

We conclude this section with a partial result toward the
first claim of Theorem~\ref{Th_treemain}.

\begin{definition}
\label{DL_T}
For a tree $T$ we define $L_T$,
the \textit{minimum-rank stripe} of $T$, as the set
\[
L_T = \{(r,s) \in \NN^2_{\mr(T)} :\  r \ge c(T), \, s \ge c(T) \}.
\]
\end{definition}
For the moment the name ``minimum-rank stripe''
is not entirely justified, since it suggests that
$L_T = \cI(T) \cap \NN^2_{\mr(T)}$.  In Section~\ref{S6}
we will show that this is the case,
but we can already show one direction of containment.
\begin{theorem}
\label{Th9}
For any tree $T$, $L_T \subseteq \cI(T)$.
\end{theorem}
\begin{proof}
  Let $k = c(T)$.
  Given any $(r,s) \in L_T$, we have
  $r \ge k$, $s \ge k$,
  and $r+s = \mr(T) = n - \MD{k}{T} + k$ by
  Observation~\ref{Ob_cMD}.
  Then by the Stars and Stripes Lemma we have $(r,s) \in \cI(T)$.
\end{proof}
\begin{corollary}
\label{Cor_rightrank}
Theorem~\ref{Th_treemain} gives the correct value of
$\mr(T)$ for $T$ a tree.
\end{corollary}


\section{Inertia formulae for a graph with a cut vertex}
\label{S4}

In this section we interrupt our discussion of inertia sets of trees
in order to derive basic formulae
about the inertia set of any graph with a cut vertex.
We obtain formulae for inertia sets that are the analogue of Theorem~16 in~\cite{Hs}
and Theorem~2.3 in~\cite{BFH1} for minimum rank.
\begin{definition}
\label{D13}
If $ Q $, $ R $ are subsets of $ \NN^2 $, then
\[
Q+R=\left \{(a+c,b+d) :\ (a,b)\in Q \ \text{and} \ (c,d)\in R\right \}.
\]
Addition of $ 3 $ or more sets is defined similarly.
\end{definition}

\begin{definition}
\label{D13b}
If $ Q $ is a subset of $ \NN^2 $ and $ n $ is a positive integer,
we let
\[
\bigl[Q\bigr]_n=Q\cap\NN^2_{\le n}.
\]
\end{definition}

We first consider the case of disconnected graphs.
Since the inertia of a direct sum of matrices is the sum of the inertias of the summands,
we have:
\begin{observation}
\label{Ob21}
Let $ G=\bigcup\limits_{i=1}^k G_i $.
Then
\[
\cI(G)=\cI(G_1)+\cI(G_2)+\cdots+\cI(G_k),
\]
and similarly for $\hcI(G)$.
\end{observation}
We now determine the inertia set of a graph with a cut vertex---see
Definition~\ref{D14}. We first recall the following useful result
\cite{Hs}, \cite{BFH1}, which reduces the minimum rank problem for
graphs to the case of $ 2 $-connected graphs.

\begin{theorem}
[Hsieh; Barioli, Fallat, Hogben]
\label{Th11}
With $ F $, $ G $ and $ F\vsum G $ as in Definition~\ref{D14},
we have
\[
\mr \FvsumGparen =
\min\bigl\{\mr(F)+\mr(G), \mr(F-v)+\mr(G-v)+2\bigr\}.
\]
\end{theorem}

Our next result generalizes this to inertia sets.
\begin{theorem}
\label{Th12} Let $ F $ and $ G $ be graphs on at least two vertices
with a common vertex $ v $ and let $ n=|F|+|G|-1 $. Then
\[
\cI \FvsumGparen =\bigl [\cI(F)+\cI(G)\bigr ]_n\cup \bigl
[\cI(F-v)+\cI(G-v)+\{(1,1)\}\bigr ]_n
\]
and
similarly for $\hcI \FvsumGparen$.
\end{theorem}
\begin{proof}
We prove the complex Hermitian version of the theorem;
the proof of the real symmetric version is the same but
with the assumption that all matrices and vectors are real.

Let $ v $ be the last vertex of $F $ and the first vertex of $ G $.

\textit{Reverse containment}:
\\
\textbf{I.}
~~Let $ (r,s)\in\bigl [\hcI(F)+\hcI(G)\bigr ]_n $.
Then $ r+s\le n $ and there exist $ (i,j)\in\hcI(F) $ and $ (k,\ell)\in\hcI(G) $
such that $ i+k=r $, $ j+\ell=s $.
Let
\[
M=
\begin{bmatrix}
A & b \\ b^* & c_1
\end{bmatrix}
\in \Herm(F)
\ \ \text{and} \ \
N=
\begin{bmatrix}
c_2 & d^* \\ d & E
\end{bmatrix}
\in \Herm(G)
 \]
with $ \pin(M)=(i,j) $ and $ \pin(N)=(k,\ell) $, and let
\[
\wh M=
\begin{bmatrix}
A & b & 0 \\ b^* & c_1 & 0 \\ 0 & 0 & 0
\end{bmatrix}
\ \ \text{and} \ \
\wh N=
\begin{bmatrix}
0 & 0 & 0 \\ 0 & c_2 & d^* \\ 0 & d & E
\end{bmatrix}
 \]
be matrices of order $ n $.
Then
\[
\pin(\wh M)=\pin(M)=(i,j), \ \ \ \pin(\wh N)=\pin(N)=(k,\ell),
 \]
and $ \wh M+\wh N\in \Herm \FvsumGparen $.
By the subadditivity of partial inertias (Proposition~\ref{P2}),
\[
\pi (\wh M+\wh N )\le\pi (\wh M )+\pi (\wh N )=i+k=r,
 \]
and
\[
\nu (\wh M+\wh N )\le\nu (\wh M )+\nu (\wh N )=j+\ell=s.
 \]
Since $  (\pi (\wh M+\wh N ),\nu (\wh M+\wh N ) ) \in\hcI
\FvsumGparen$ by definition, and $ r+s\le n $, $ (r,s)\in\hcI
\FvsumGparen $ by the Northeast Lemma (Lemma~\ref{L5}). Thus, we
have $ \bigl[\hcI(F)+\hcI(G)\bigr]_n\subseteq\hcI \FvsumGparen $.

\caseskip
\textbf{II.}
~~Now let $ (r,s)\in\bigl[\hcI(F-v)+\hcI(G-v)+\left \{(1,1)\right \}\bigr]_n $.
Then $ r+s\le n $ and there exist
$ (i,j)\in\hcI(F-v) $ and $ (k,\ell)\in\hcI(G-v) $
with $ (i,j)+(k,\ell)+(1,1)=(r,s) $.
Let $ A\in \Herm(F-v) $ with $ \pin(A)=(i,j) $ and
let $ E\in \Herm(G-v) $ with $ \pin(E)=(k,\ell) $.
Choose $ b $, $ c $, $ d $ such that
\[
M=
\begin{bmatrix}
A & b & 0 \\ b^* & c & d^* \\ 0 & d & E
\end{bmatrix}
\in \Herm \FvsumGparen.
 \]
By Proposition~\ref{P1},
\[
\pi(M)\le\pi\left (
\begin{bmatrix}
A & 0 \\ 0 & E
\end{bmatrix}\right )
+1=\pi(A)+\pi(E)+1=i+k+1=r,
 \]
and, similarly,
\[
  \nu(M)\le j+\ell+1=s.
\]
Since $ \bigl(\pi(M),\nu(M)\bigr)\in\hcI \FvsumGparen $,
and $ r+s\le n $,
by the Northeast Lemma, $ (r,s)\in\hcI \FvsumGparen $.

So we have
\[
\bigl[\hcI(F-v)+\hcI(G-v)+\bigl\{(1,1)\bigr\}\bigr]_n
\subseteq\hcI \FvsumGparen.
 \]

\textit{Forward containment}:
\\
Now let $ (i,j)\in\hcI \FvsumGparen $.
By Observation~\ref{Ob4}, $ i+j\le n $.
Let
\[
M=
\begin{bmatrix}
A & b & 0 \\ b^* & c & d^* \\ 0 & d & E
\end{bmatrix}
\in \Herm \FvsumGparen.
 \]
with $ \pin(M)=(i,j) $.
Then
\[
\rank A+\rank E\le\rank
\begin{bmatrix}
A & b & 0 \\ 0 & d & E
\end{bmatrix}
\le\rank M\le\rank A+\rank E+2.
 \]
If the first and third inequalities are strict,
then
\[
\rank A+\rank E+1=\rank
\begin{bmatrix}
A & b & 0 \\ 0 & d & E
\end{bmatrix}
=\rank M.
 \]
The first equality implies that either $ b\notin\Col(A) $ or else $ d\notin\Col(E) $,
while the second equality implies that $ b\in\Col(A) $ and $ d\in\Col(E) $.
So this case does not occur and either
\[
\rank A+\rank E=\rank\begin{bmatrix}
A & b & 0 \\ 0 & d & E
\end{bmatrix}
 \]
or else
\[
  \rank M=\rank A+\rank E+2.
\]

\caseskip
\textbf{I.}
~~$ \rank A+\rank E=\rank\begin{bmatrix}
A & b & 0 \\ 0 & d & E
\end{bmatrix} $.
\\
Then $ b\in\Col(A) $ and $ d\in\Col(E) $.
So $ b=Au $, $ d=Ev $ for some $ u\in\CC^{|F|-1} $, $ v\in \CC^{|G|-1} $.

Define
\[
\wh A=
\begin{bmatrix}
A & Au \\ u^*A & u^*Au
\end{bmatrix}
\in \Herm(F);
\ \ \
\wh E=
\begin{bmatrix}
v^*Ev & v^*E \\ Ev & E
\end{bmatrix}
\in \Herm(G).
 \]
Then $ \wh A $ is congruent to $ A\oplus\underset{1\times 1}{[~0~]} $;
$ \wh E $ is congruent to $ E\oplus\underset{1\times 1}{[~0~]} $.
Hence
\[
\begin{array}{lp{1cm}l}
\pi(\wh A\, )=\pi(A); && \nu(\wh A\, )=\nu(A);
\\
\pi(\wh E\, )=\pi(E); && \nu(\wh E\, )=\nu(E).
\end{array}
 \]
Also,
\begin{eqnarray}
\rank\wh A=\rank A\le|F|-1,
\label{eq1}
\\
\rank\wh E=\rank E\le|G|-1.
\label{eq2}
\end{eqnarray}
By Proposition~\ref{P1}, $ \exists\,a,b\in\{0,1\} $ such that
\begin{eqnarray}
i=\pi(M)=\pi(A)+\pi(E)+a=\pi(\wh A\,)+\pi(\wh E\,)+a,
\label{eq3}
\\
j=\nu(M)=\nu(A)+\nu(E)+b=\nu(\wh A\,)+\nu(\wh E\,)+b.
\label{eq4}
\end{eqnarray}
It follows from~\eqref{eq1} and~\eqref{eq2} that
\begin{eqnarray*}
\pi(\wh A\,)+a+\nu(\wh A\,)=\rank\wh A+a\le|F|-1+a\le|F|,
\\
\pi(\wh E\,)+\nu(\wh E\,)+b=\rank\wh E+b\le|G|-1+b\le|G|.
\end{eqnarray*}
Hence, by the Northeast Lemma,
\begin{eqnarray*}
\bigl(\pi(\wh A\,)+a,\nu(\wh A\,)\bigr)=\bigl(\pi(A)+a,\nu(A)\bigr)\in\hcI(F),
\\
\bigl(\pi(\wh E\,),\nu(\wh E\,)+b\bigr)=\bigl(\pi(E),\nu(E)+b\bigr)\in\hcI(G),
\end{eqnarray*}
and since these two vectors add up to $ (i,j) $, by~\eqref{eq3}
and~\eqref{eq4}, we conclude that $ (i,j)\in\hcI(F)+\hcI(G) $.

Since $ i+j\le n $ we get $ (i,j)\in\bigl[\hcI(F)+\hcI(G)\bigr]_n $.
So in this case, $ \hcI \FvsumGparen \subseteq\bigl[\hcI(F)+\hcI(G)\bigr]_n $.

\caseskip
\textbf{II.}
~~$ \rank M=\rank A+\rank E+2 $.
\\
By Proposition~\ref{P1}, we have
\begin{eqnarray*}
i+j &\le & \pi\left (\begin{bmatrix} A & 0 \\ 0 & E \end{bmatrix} \right )+1
                +\nu\left (\begin{bmatrix} A & 0 \\ 0 & E \end{bmatrix} \right )+1
\\
        &=& \pi(A)+\pi(E)+1+\nu(A)+\nu(E)+1
\\
        &=& \pi(A)+\nu(A)+\pi(E)+\nu(E)+2
\\
        &=& \rank(A)+\rank(E)+2=\rank M=i+j.
\end{eqnarray*}
It follows that $ i=\pi(A)+\pi(E)+1 $ and $ j=\nu(A)+\nu(E)+1 $ and
$ (i,j)=\bigl(\pi(A),\nu(A)\bigr)+\bigl(\pi(E),\nu(E)\bigr)+(1,1) $.
By definition, $ (i,j)\in\hcI(F-v)+\hcI(G-v)+\{(1,1)\} $, and since
$ i+j\le n $, $ (i,j)\in\bigl[\hcI(F-v)+\hcI(G-v)+\{(1,1)\}\bigr]_n
$. So in this case,
\[
\hcI \FvsumGparen \subseteq\bigl[\hcI(F-v)+\hcI(G-v)+\{(1,1)\}\bigr]_n.
 \]
This completes the proof of the forward containment.
\end{proof}
It is straightforward to show that Theorem~\ref{Th11} is a corollary
of Theorem~\ref{Th12}. The proof is not illuminating, so we do not
include it.
\begin{example}
\label{Ex7}
Let $ F=S_4 $ and $ G=P_3 $ with $ v $ a pendant vertex in $ S_4 $
and the degree $ 2 $ vertex in $ P_3 $.
Then $ T=F\vsum G $ is the graph below.
\setlength{\unitlength}{1pt}
\[
\begin{picture}(77,50)
 \put(0,0){\includegraphics{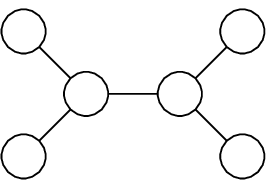}}
 \put(49,22){$v$}
\end{picture}
\]
From Examples~\ref{Ex3} and~\ref{Ex2} we have

\setlength{\unitlength}{1mm}
\begin{picture}(50,38)
 \put(5,20){$\mathcal I(S_4)$:}
 \put(22,4){\vector(0,1){30}}
 \put(20,6){\vector(1,0){31}}
 \put(20.75,28){\line(1,0){2.5}}
 \put(22,28){\circle*{1.5}}
 \put(20.75,22.5){\line(1,0){2.5}}
 \put(22,22.5){\circle*{1.5}}
 \put(27.5,22.5){\circle*{1.5}}
 \put(20.75,17){\line(1,0){2.5}}
  \put(27.5,17){\circle*{1.5}}
 \put(33,17){\circle*{1.5}}
 \put(20.75,11.5){\line(1,0){2.5}}
 \put(27.5,11.5){\circle*{1.5}}
 \put(27.5,4.75){\line(0,1){2.5}}
 \put(33,11.5){\circle*{1.5}}
 \put(33,4.75){\line(0,1){2.5}}
 \put(38.5,11.5){\circle*{1.5}}
 \put(38.5,4.75){\line(0,1){2.5}}
 \put(43.5,4.75){\line(0,1){2.5}}
 \put(38.5,6){\circle*{1.5}}
 \put(44,6){\circle*{1.5}}
 \put(65,20){$\mathcal I(P_3)$:}
 \put(82,4){\vector(0,1){30}}
 \put(80,6){\vector(1,0){27}}
 \put(80.75,22.5){\line(1,0){2.5}}
 \put(82,22.5){\circle*{1.5}}
 \put(80.75,17){\line(1,0){2.5}}
 \put(82,17){\circle*{1.5}}
 \put(87.5,17){\circle*{1.5}}
 \put(80.75,11.5){\line(1,0){2.5}}
 \put(87.5,11.5){\circle*{1.5}}
 \put(87.5,4.75){\line(0,1){2.5}}
 \put(93,11.5){\circle*{1.5}}
 \put(93,4.75){\line(0,1){2.5}}
 \put(93,6){\circle*{1.5}}
 \put(98.5,4.75){\line(0,1){2.5}}
 \put(98.5,6){\circle*{1.5}}
\end{picture}

\setlength{\unitlength}{1mm}
\begin{picture}(50,38)
 \put(0,10){$\mathcal I(S_4 - v) = \mathcal I(P_3)$:}
 \put(40,4){\vector(0,1){27}}
 \put(38,6){\vector(1,0){27}}
 \put(38.75,22.5){\line(1,0){2.5}}
 \put(40,22.5){\circle*{1.5}}
 \put(38.75,17){\line(1,0){2.5}}
 \put(40,17){\circle*{1.5}}
 \put(45.5,17){\circle*{1.5}}
 \put(38.75,11.5){\line(1,0){2.5}}
 \put(45.5,11.5){\circle*{1.5}}
 \put(45.5,4.75){\line(0,1){2.5}}
 \put(51,11.5){\circle*{1.5}}
 \put(51,4.75){\line(0,1){2.5}}
 \put(51,6){\circle*{1.5}}
 \put(56.5,4.75){\line(0,1){2.5}}
 \put(56.5,6){\circle*{1.5}}
 \put(70,10){$\mathcal I(P_3 - v) = \mathcal I(2K_1)$:}
 \put(112,4){\vector(0,1){21}}
 \put(110,6){\vector(1,0){20}}
 \put(110.75,17){\line(1,0){2.5}}
 \put(112,17){\circle*{1.5}}
 \put(110.75,11.5){\line(1,0){2.5}}
 \put(112,11.5){\circle*{1.5}}
 \put(117.5,11.5){\circle*{1.5}}
 \put(117.5,4.75){\line(0,1){2.5}}
 \put(112,6){\circle*{1.5}}
 \put(117.5,6){\circle*{1.5}}
 \put(123,4.75){\line(0,1){2.5}}
 \put(123,6){\circle*{1.5}}
\end{picture}

\noindent It follows that

\setlength{\unitlength}{1mm}
\begin{picture}(50,50)
 \put(5,20){$[\mathcal I(S_4) + \mathcal I(P_3)]_6$ is}
 \put(42,4){\vector(0,1){40}}
 \put(40,6){\vector(1,0){41}}
 \put(40.75,39){\line(1,0){2.5}}
 \put(42,39){\circle*{1.5}}
 \put(40.75,33.5){\line(1,0){2.5}}
 \put(42,33.5){\circle*{1.5}}
 \put(47.5,33.5){\circle*{1.5}}
 \put(40.75,28){\line(1,0){2.5}}
 \put(47.5,28){\circle*{1.5}}
 \put(53,28){\circle*{1.5}}
 \put(40.75,22.5){\line(1,0){2.5}}
 \put(47.5,22.5){\circle*{1.5}}
 \put(53,22.5){\circle*{1.5}}
 \put(58.5,22.5){\circle*{1.5}}
 \put(40.75,17){\line(1,0){2.5}}
 \put(53,17){\circle*{1.5}}
 \put(58.5,17){\circle*{1.5}}
 \put(64,17){\circle*{1.5}}
 \put(40.75,11.5){\line(1,0){2.5}}
 \put(58.5,11.5){\circle*{1.5}}
 \put(64,11.5){\circle*{1.5}}
 \put(69.5,11.5){\circle*{1.5}}
\put(47.5,4.75){\line(0,1){2.5}}
 \put(53,4.75){\line(0,1){2.5}}
 \put(58.5,4.75){\line(0,1){2.5}}
 \put(64,4.75){\line(0,1){2.5}}
 \put(69.5,4.75){\line(0,1){2.5}}
 \put(69.5,6){\circle*{1.5}}
 \put(75,4.75){\line(0,1){2.5}}
 \put(75,6){\circle*{1.5}}
 \put(100,20){and}
\end{picture}

\setlength{\unitlength}{1mm}
\begin{picture}(50,50)
 \put(0,20){$[\mathcal I(P_3) + \mathcal I(2K_1)+\{(1,1)\}]_6$ is}
 \put(62,4){\vector(0,1){40}}
 \put(60,6){\vector(1,0){41}}
 \put(60.75,39){\line(1,0){2.5}}
 \put(60.75,33.5){\line(1,0){2.5}}
 \put(67.5,33.5){\circle*{1.5}}
 \put(60.75,28){\line(1,0){2.5}}
 \put(67.5,28){\circle*{1.5}}
 \put(73,28){\circle*{1.5}}
 \put(60.75,22.5){\line(1,0){2.5}}
 \put(67.5,22.5){\circle*{1.5}}
 \put(73,22.5){\circle*{1.5}}
 \put(78.5,22.5){\circle*{1.5}}
 \put(60.75,17){\line(1,0){2.5}}
 \put(73,17){\circle*{1.5}}
 \put(78.5,17){\circle*{1.5}}
 \put(84,17){\circle*{1.5}}
 \put(60.75,11.5){\line(1,0){2.5}}
 \put(78.5,11.5){\circle*{1.5}}
 \put(84,11.5){\circle*{1.5}}
 \put(89.5,11.5){\circle*{1.5}}
 \put(67.5,4.75){\line(0,1){2.5}}
 \put(73,4.75){\line(0,1){2.5}}
 \put(78.5,4.75){\line(0,1){2.5}}
 \put(84,4.75){\line(0,1){2.5}}
 \put(89.5,4.75){\line(0,1){2.5}}
 \put(95,4.75){\line(0,1){2.5}}
\end{picture}

\noindent Then $\cI(T) = \cI \FvsumGparen = [\cI(S_4) + \cI(P_3)]_6
\cup [\cI(P_3) + \cI(2K_1) + \{(1, 1)\}]_6$ is:

\setlength{\unitlength}{1mm}
\begin{picture}(50,50)
 \put(42,4){\vector(0,1){40}}
 \put(40,6){\vector(1,0){41}}
 \put(40.75,39){\line(1,0){2.5}}
 \put(42,39){\circle*{1.5}}
 \put(40.75,33.5){\line(1,0){2.5}}
 \put(42,33.5){\circle*{1.5}}
 \put(47.5,33.5){\circle*{1.5}}
 \put(40.75,28){\line(1,0){2.5}}
 \put(47.5,28){\circle*{1.5}}
 \put(53,28){\circle*{1.5}}
 \put(40.75,22.5){\line(1,0){2.5}}
 \put(47.5,22.5){\circle*{1.5}}
 \put(53,22.5){\circle*{1.5}}
 \put(58.5,22.5){\circle*{1.5}}
 \put(40.75,17){\line(1,0){2.5}}
 \put(53,17){\circle*{1.5}}
 \put(58.5,17){\circle*{1.5}}
 \put(64,17){\circle*{1.5}}
 \put(40.75,11.5){\line(1,0){2.5}}
 \put(58.5,11.5){\circle*{1.5}}
 \put(64,11.5){\circle*{1.5}}
 \put(69.5,11.5){\circle*{1.5}}
\put(47.5,4.75){\line(0,1){2.5}}
 \put(53,4.75){\line(0,1){2.5}}
 \put(58.5,4.75){\line(0,1){2.5}}
 \put(64,4.75){\line(0,1){2.5}}
 \put(69.5,4.75){\line(0,1){2.5}}
 \put(69.5,6){\circle*{1.5}}
 \put(75,4.75){\line(0,1){2.5}}
 \put(75,6){\circle*{1.5}}
\end{picture}

\noindent
Since $ |T|=6 $ and $ P(T)=2 $, $ \mr(T)=4 $ by Theorem~\ref{Th2}.
Since $ \left |E_T\left (\{v\}\right )\right |=3 $,
$ \left |E_T\left (\{v\}\right )\right |-2\left |\{v\}\right |+1=2=P(T) $,
so $ \{v\} $ is a minimal optimal set for $ T $.
We observe that in this case $L_T = \cI(T) \cap \NN^2_{\mr(T)}$.
\end{example}
We pause to develop some additional fundamental properties of
inertia sets before generalizing Theorem~\ref{Th12}. The next result
generalizes the fact \cite{N} that $\mr(G-v) \le \mr(G) \le
\mr(G-v)+2$.
\begin{proposition}
\label{P22}
 Let $ G $ be any graph on $ n $ vertices
and let $ v $ be any vertex of $ G $.
Then we have:
\begin{enumerate}
\item [(a)]
$ \bigl[\cI(G)\bigr]_{n-1}\subseteq\cI(G-v) $.
\item [(b)]
$ \cI(G)\supseteq\bigl[\cI(G-v)\bigr]_{n-2}+\{(1,1)\} $.
\end{enumerate}
The same inclusions hold in the Hermitian case.
\end{proposition}
\begin{proof}
Let $ (r,s)\in\bigl[\cI(G)\bigr]_{n-1} $.
Then $ r+s\le n-1 $.
Let $ A\in \Sym(G) $ with $ \pin(A)=(r,s) $,
and let $B$ be the principal submatrix of $A$
obtained by deleting the row and column $v$.
Then $ B \in \Sym(G-v) $ and by the interlacing inequalities $ \pin\bigl(B\bigr) $
is one of $ (r,s) $, $ (r-1,s) $, $ (r,s-1) $, or $ (r-1,s-1) $.
Then one of these is in $ \cI(G-v) $
so by the Northeast Lemma, $ (r,s)\in\cI(G-v) $.
This proves (a).

Now let $ (r,s)\in\bigl[\cI(G-v)\bigr]_{n-2} $
so $ r+s\le n-2 $.
Choose $ A\in \Sym(G) $ in such a way that the principal submatrix $B$
obtained by deleting row and column $v$ satisfies
$ \pin\bigl(B\bigr)=(r,s) $.
Then by the interlacing inequalities, $ \pin(A) $
is one of $ (r,s) $, $ (r+1,s) $, $ (r,s+1) $, or $ (r+1,s+1) $.
Since $ r+1+s+1\le n $, $ (r+1,s+1)\in\cI(G) $
by the Northeast Lemma applied to $ G $.
This completes the proof of (b).

The proof of the Hermitian case is the same, but with Hermitian notation.
\end{proof}
\begin{proposition}
\label{P23}
If $ v $ is a pendant vertex of the graph $ G $ and $ (i,j)\in\cI(G-v) $,
then $ (i+1,j)\in\cI(G) $ and $ (i,j+1)\in\cI(G) $, and similarly
for $\hcI(G-v)$ and $\hcI(G)$.
\end{proposition}
\begin{proof}
As usual, the proofs of the real symmetric and Hermitian versions
do not differ materially.
Let $ v $ be the first vertex of $ G $ and let its neighbor $ u $ be the second.
Let $ A\in \Sym(G-v) $ with $ \pin(A)=(i,j) $.
Then
\[
M=\begin{bmatrix}
J_2 & 0 \\ 0 & 0
\end{bmatrix}
+
\begin{bmatrix}
0 & 0 \\ 0 & A
\end{bmatrix}
\in \Sym(G)
 \]
and $ \rank M=1+\rank A $.
By Proposition~\ref{P3}
\[
 \pi(M) \le\pi(A)+1 \ \ \text{and} \ \ \nu(M)\le\nu(A).
\]
Then $ \rank M=\pi(M)+\nu(M)\le\pi(A)+1+\nu(A)=\rank A+1=\rank M $
and $ (i+1,j)=\bigl(\pi(A)+1,\nu(A)\bigr)=\bigl(\pi(M),\nu(M)\bigr)\in\cI(G) $.
Similarly, $ (i,j+1)\in\cI(G) $.
\end{proof}
The following corollary of Theorem~\ref{Th11} is very useful in
simplifying the calculation of the minimum rank of a graph.
\begin{proposition}
[{\cite[Lemma~38]{S}}]
\label{P24}
If the degree of $ v $ is $ 2 $ in $ F\vsum G $,
then
\[
\mr \FvsumGparen =\mr(F)+\mr(G).
 \]
\end{proposition}
The following result generalizes this fact to inertia sets.
\begin{proposition}
\label{P25}
If the degree of $ v $ is $ 2 $ in $ F\vsum G $,
and $ n=|F|+|G|-1 $, then
\[
\cI \FvsumGparen =\bigl[\cI(F)+\cI(G)\bigr]_n,
 \]
and similarly for $\hcI \FvsumGparen$.
\end{proposition}
\begin{proof}
By Theorem~\ref{Th12} it suffices to show that
\[
\bigl[\cI(F-v)+\cI(G-v)+\{(1,1)\}\bigr]_n\subseteq\bigl[\cI(F)+\cI(G)\bigr]_n.
 \]
Let $ (r,s)\in\bigl[\cI(F-v)+\cI(G-v)+\{(1,1)\}\bigr]_n $.
Then $ r+s\le n $ and $ (r,s)=(i,j)+(k,\ell)+(1,1) $
with $ (i,j)\in\cI(F-v) $ and $ (k,\ell)\in\cI(G-v) $.
Since $ v $ is pendant in both $ F $ and $ G $,
by Proposition~\ref{P23}, $ (i+1,j)\in\cI(F) $ and $ (k,\ell+1)\in\cI(G) $,
so $ (r,s)=(i+1+k,j+\ell+1)\in\cI(F)+\cI(G) $.
Since $ r+s\le n $, $ (r,s)\in\bigl[\cI(F)+\cI(G)\bigr]_n $.

Replacing $\cI$ by $\hcI$ uniformly proves the Hermitian case.
\end{proof}
\begin{example}
\label{Ex8}
Let $ F=G=S_4 $ and let $ v $ be a pendant vertex in each of $F$ and $G$
so that $ T=F\vsum G $
is the graph below.
\setlength{\unitlength}{1pt}
\[
\begin{picture}(104,50)
 \put(0,0){\includegraphics{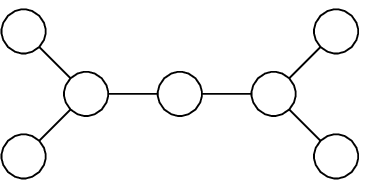}}
 \put(49,22){$v$}
\end{picture}
\]
By Proposition~\ref{P25},
$ \cI(T)=\bigl[\cI(F)+\cI(G)\bigr]_7=\bigl[\cI(S_4)+\cI(S_4)\bigr]_7 $.
Knowing that the inertia set $\cI(S_4)$ is
\[
\setlength{\unitlength}{1mm}
\begin{picture}(30,38)
 \put(2,4){\vector(0,1){30}}
 \put(0,6){\vector(1,0){30}}
 \put(0.75,28){\line(1,0){2.5}}
 \put(2,28){\circle*{1.5}}
 \put(0.75,22.5){\line(1,0){2.5}}
 \put(2,22.5){\circle*{1.5}}
 \put(7.5,22.5){\circle*{1.5}}
 \put(0.75,17){\line(1,0){2.5}}
  \put(7.5,17){\circle*{1.5}}
 \put(13,17){\circle*{1.5}}
 \put(0.75,11.5){\line(1,0){2.5}}
 \put(7.5,11.5){\circle*{1.5}}
 \put(7.5,4.75){\line(0,1){2.5}}
 \put(13,11.5){\circle*{1.5}}
 \put(13,4.75){\line(0,1){2.5}}
 \put(18.5,11.5){\circle*{1.5}}
 \put(18.5,4.75){\line(0,1){2.5}}
 \put(23.5,4.75){\line(0,1){2.5}}
 \put(18.5,6){\circle*{1.5}}
 \put(23.5,6){\circle*{1.5}}
\end{picture}
\]
allows us to calculate $ \cI(T) $, as depicted below.
\[
\setlength{\unitlength}{1mm}
\begin{picture}(52,54)
 \put(2,4){\vector(0,1){48}}
 \put(0,6){\vector(1,0){48}}
 \put(0.75,44.5){\line(1,0){2.5}}
 \put(2,44.5){\circle*{1.5}}
 \put(0.75,39){\line(1,0){2.5}}
 \put(2,39){\circle*{1.5}}
 \put(7.5,39){\circle*{1.5}}
 \put(0.75,33.5){\line(1,0){2.5}}
 \put(7.5,33.5){\circle*{1.5}}
 \put(13,33.5){\circle*{1.5}}
 \put(0.75,28){\line(1,0){2.5}}
 \put(7.5,28){\circle*{1.5}}
 \put(13,28){\circle*{1.5}}
 \put(18.5,28){\circle*{1.5}}
 \put(0.75,22.5){\line(1,0){2.5}}
 \put(13,22.5){\circle*{1.5}}
 \put(18.5,22.5){\circle*{1.5}}
 \put(24,22.5){\circle*{1.5}}
 \put(0.75,17){\line(1,0){2.5}}
 \put(13,17){\circle*{1.5}}
 \put(18.5,17){\circle*{1.5}}
 \put(24,17){\circle*{1.5}}
 \put(29.5,17){\circle*{1.5}}
 \put(0.75,11.5){\line(1,0){2.5}}
 \put(24,11.5){\circle*{1.5}}
 \put(29.5,11.5){\circle*{1.5}}
 \put(35,11.5){\circle*{1.5}}
\put(7.5,4.75){\line(0,1){2.5}}
 \put(13,4.75){\line(0,1){2.5}}
 \put(18.5,4.75){\line(0,1){2.5}}
 \put(24,4.75){\line(0,1){2.5}}
 \put(29.5,4.75){\line(0,1){2.5}}
 \put(35,4.75){\line(0,1){2.5}}
 \put(35,6){\circle*{1.5}}
 \put(40.5,4.75){\line(0,1){2.5}}
 \put(40.5,6){\circle*{1.5}}
\end{picture}
\]
Here, $ \mr(T)=4 $ is attained only at the partial inertia $ (2,2)
$, and one can easily check that $ c(T)=2 $, so that $ L_T=\{(2,2)\}
$.
\end{example}

We close this section with the generalization of Theorem~\ref{Th12}.
We first extend Definition~\ref{D14}.
\begin{definition}
\label{D15} Let $ G_1,G_2,\ldots,G_k $, $ k\ge 2 $, be graphs on at
least two vertices with a common vertex $ v $ and let $
G=\bigvsum\limits_{i=1}^k G_i $ be the graph on $
n=\sum\limits_{i=1}^k\left |G_i\right |-(k-1) $ vertices obtained by
identifying the vertex $ v $ in each of the $ G_i $. We call $ G $
the \textit{vertex sum} of the graphs $ G_1,G_2,\ldots,G_k $ at $ v
$.
\end{definition}
\begin{theorem}
\label{Th13}
Let $ G $  be a graph on $ n\ge 3 $ vertices and let $ v $ be a cut vertex of~$ G $.
Write $ G=\bigvsum\limits_{i=1}^k G_i $, $k \ge 2$,
the vertex sum of \ $ G_1,G_2,\ldots,G_k $ at $ v $.
Then
\begin{eqnarray}
\cI(G) &=&\bigl[\cI(G_1)+\cI(G_2)+\cdots+\cI(G_k)\bigr]_n
\label{eq5}
\\
&&\cup\ \bigl[\cI(G_1-v)+\cI(G_2-v)+\cdots+\cI(G_k-v)+\{(1,1)\}\bigr]_n,
\nonumber
\end{eqnarray}
and similarly for $\hcI(G)$.
\end{theorem}
\begin{proof}
The idea of the proof is the same as in the proof of Theorem~\ref{Th12},
which is showing that each side of equation~\eqref{eq5}
is contained in the other.
Since each of the theorems cited applies equally well to $\hcI$ as to $\cI$,
the same proof demonstrates both cases.

\textit{Forward containment}:
\\
We prove that
\begin{eqnarray}
\cI(G) &\subseteq &\bigl[\cI(G_1)+\cI(G_2)+\cdots+\cI(G_k)\bigr]_n
\label{eq6}
\\
&&\cup\
\bigl[\cI(G_1-v)+\cI(G_2-v)+\cdots+\cI(G_k-v)+\{(1,1)\}\bigr]_n
\nonumber
\end{eqnarray}
by induction on $ k $. For $ k=2 $ this follows from
Theorem~\ref{Th12}. Assume \eqref{eq6} holds for all integers $ j $
with $ 2\le j<k $. Let $ G^\prime=\bigvsum\limits_{i=1}^{k-1} G_i $,
the vertex sum of \ $ G_1,\ldots,G_{k-1} $ at $ v $ and let $
n^\prime=|G^\prime| $. Then by Theorem~\ref{Th12},
\begin{eqnarray*}
\cI(G) &=& \cI (G^\prime\vsum G_k )\subseteq
\bigl[\cI(G^\prime)+\cI(G_k)\bigr]_n
\\
&&\cup\ \bigl[\cI(G^\prime-v)+\cI(G_k-v)+\{(1,1)\}\bigr]_n.
\end{eqnarray*}
But
\[
\cI(G^\prime-v)=\cI\left (\bigcup\limits_{i=1}^{k-1}(G_i-v)\right )=\cI(G_1-v)+\cdots+\cI(G_{k-1}-v)
 \]
by Observation~\ref{Ob21}.
Applying the induction hypothesis to $ \cI(G^\prime) $ we have
\begin{eqnarray*}
\cI(G) &\subseteq &
\Bigl[\Bigl\{\bigl[\cI(G_1)+\cdots+\cI(G_{k-1})\bigr]_{n^\prime}
\\
&& \ \  \cup \ \bigl[ \cI(G_1-v)+\cdots+\cI(G_{k-1}-v)+\{(1,1)\}\bigr]_{n^\prime}\Bigr\}+\cI(G_k)\Bigr]_n
\\
&& \cup \ \bigl[ \cI(G_1-v)+\cdots+\cI(G_{k}-v)+\{(1,1)\}\bigr]_{n}
\\
&=& \Bigl[\Bigl(\bigl[\cI(G_1)+\cdots+\cI(G_{k-1})\bigr]_{n^\prime}+\cI(G_k)\Bigr)
\\
&& \ \  \cup \ \Bigl(\bigl[ \cI(G_1-v)+\cdots+\cI(G_{k-1}-v)+\{(1,1)\}\bigr]_{n^\prime}
+\cI(G_k)\Bigr)\Bigr]_n
\\
&& \cup \ \bigl[ \cI(G_1-v)+\cdots+\cI(G_{k}-v)+\{(1,1)\}\bigr]_{n}
\\
&=& \Bigl[\bigl[\cI(G_1)+\cdots+\cI(G_{k-1})\bigr]_{n^\prime}+\cI(G_k)\Bigr]_n
\\
&& \cup \ \Bigr[\bigl[ \cI(G_1-v)+\cdots+\cI(G_{k-1}-v)+\{(1,1)\}\bigr]_{n^\prime}
+\cI(G_k)\Bigr]_n
\\
&& \cup \ \bigl[ \cI(G_1-v)+\cdots+\cI(G_{k}-v)+\{(1,1)\}\bigr]_{n}
\\
&\subseteq & \bigl[\cI(G_1)+\cdots+\cI(G_k)\bigr]_n
\\
&& \cup \ \Bigl[\Bigl(\bigl[\cI(G_1-v)+\cdots+\cI(G_{k-1}-v)\bigr]_{n^\prime-2}+\{(1,1)\}\Bigr)
+\cI(G_k)\Bigr]_n
\\
&& \cup \ \bigl[ \cI(G_1-v)+\cdots+\cI(G_{k}-v)+\{(1,1)\}\bigr]_{n}.
\end{eqnarray*}
Let
\begin{eqnarray*}
Q_1 &=& \bigl[\cI(G_1)+\cdots+\cI(G_k)\bigr]_n,
\\
Q_2 &=& \bigl[ \cI(G_1-v)+\cdots+\cI(G_k-v)+\{(1,1)\}\bigr]_{n},
\\
Q_0 &=& \Bigl[\bigl[\cI(G_1-v)+\cdots+\cI(G_{k-1}-v)\bigr]_{n^\prime-2}+\{(1,1)\}
+\cI(G_k)\Bigr]_n.
\end{eqnarray*}
We show that $ Q_0\subseteq Q_2 $. Suppose that $ (r,s)\in Q_0 $.
Then
\[
(r,s)=(i_1,j_1)+(i_2,j_2)+\cdots+(i_{k-1},j_{k-1})+(1,1)+(i,j)
 \]
with $ (i_t,j_t)\in\cI(G_t-v) $, $ t=1,\ldots,k-1 $, $ (i,j)\in\cI(G_k) $,
\[
\sum\limits_{t=1}^{k-1}\left (i_t+j_t\right )\le n^\prime-2, \ \
\text{and} \ \ r+s\le n.
 \]
If $ i+j<|G_k| $, by Proposition~\ref{P22}(a),  $ (i,j)\in\cI(G_k-v)
$ and then $ (r,s)\in Q_2 $. So suppose that $ i+j=|G_k| $. At least
one of $ i $,  $ j $ is greater than $ 0 $. Without loss of
generality, assume $ i>0 $. By Proposition~\ref{P6}, $
(i-1,j)\in\cI(G_k) $, and by Proposition~\ref{P22}(a), $
(i-1,j)\in\cI(G_k-v) $. Since $
n^\prime=|G^\prime|=\sum\limits_{t=1}^{k-1}|G_t|-(k-2) $, we have
$$
n^\prime-2=\sum\limits_{t=1}^{k-1}\left (|G_t|-1\right )-1=\left
(\sum\limits_{t=1}^{k-1}|G_t-v|\right )-1.
$$
Therefore, $
\sum\limits_{t=1}^{k-1}(i_t+j_t)<\sum\limits_{t=1}^{k-1}|G_t-v| $.
Without loss of generality, assume $ i_1+j_1<|G_1-v| $. By the
Northeast Lemma $ (i_1+1,j_1)\in\cI(G_1-v) $. Since $ r+s\le n $,
and $
(r,s)=(i_1+1,j_1)+(i_2,j_2)+\cdots+(i_{k-1},j_{k-1})+(i-1,j)+(1,1)
$, we again have $ (r,s)\in Q_2 $. This completes the proof that $
Q_0\subseteq Q_2 $. Therefore $ \cI(G)\subseteq Q_1\cup Q_2 $, which
is~\eqref{eq6}.

\textit{Reverse containment}:
\\
A proof by induction is not straightforward.
However, one can show the two containments
$ \bigl[\cI(G_1)+\cdots+\cI(G_k)\bigr]_n\subseteq\cI(G) $, and
$ \bigl[\cI(G_1-v)+\cdots+\cI(G_k-v)+\{(1,1)\}\bigr]_n\subseteq\cI(G) $,
by simply imitating each step in the proof of Theorem~\ref{Th12}.
As there are no new ideas in the proof, we omit it.
\end{proof}


\section{The cut-vertex formula for elementary inertias}
\label{S5}

The results of the previous section give us a way to inductively
calculate the inertia set of any graph once we know the inertia sets
of $2$-connected graphs. In this section we prove that the same
inductive formula holds when calculating the set of elementary
inertias. Claim 1 of Theorem~\ref{Th_treemain} will then follow
because a forest is a graph with no $2$-connected subgraph on $3$ or
more vertices.

It is convenient to describe the elementary inertias of a graph $G$
in terms of \BR edge-colorings of certain subgraphs of $G$.

\begin{definition}
  Let $G$ be a graph on $n$ vertices, let $S$ be a subset of $V(G)$,
  and let $\B$ and $\R$ be disjoint subsets of $E(G-S)$. The ordered
  triple $(S,\B,\R)$ is called a \textit{\BR span of $G$} if
  $(V \setminus S,\B \cup \R)$ is a spanning forest of $G-S$.
  (A \textit{spanning forest} of a graph consists of a spanning
  tree for each connected component.)
  If $(S,\B,\R)$ is a \BR span of
  $G$, we say that the ordered pair $(|S|+|\B|,|S|+|\R|)$ is a
  \textit{color vector} of $G$.  The set of color vectors of $G$ is
  denoted $\cC(G)$.
\end{definition}

The color vector counts how many edges of the spanning forest have
been marked with either the first color or the second color, but it
also counts the set $S$ of excluded vertices twice, as though each
such vertex were marked simultaneously with both colors. Because
every spanning forest has the same number of edges, the quantity
$|\B|+|\R|$ depends only on $S$, and for a given size $|S| = k$,
$|\B| + |\R|$ is minimized if $G-S$ has $\MD{k}{G}$ components. If
$G$ is a graph on $n$ vertices, $\ell = \MD{0}{G}$ is the number of
components of $G$, and $(\emptyset,\B,\R)$ is a \BR span of $G$,
then $|\B|+|\R|+\ell=n$.

\begin{observation}
\label{Ob_cvec_stripe}
  If $G$ is a graph on $n$ vertices and $\ell = \MD{0}{G}$, then
  $\NN^2_{n-\ell} \subseteq \cC(G)$.
\end{observation}

\begin{definition}
  If $Q$ is a subset of $\NN^2$, we define the
  \textit{northeast expansion of $Q$} as
  \[
    \nea{Q} = Q + \NN^2.
  \]
\end{definition}

For example, the Northeast Lemma is equivalent to the statement
that, for $G$ a graph on $n$ vertices, $\uptrunc{n}{\cI(G)}
\subseteq \cI(G)$.  The prevalence of northeast expansions in this
section leads us to define the following equivalence relation:
\begin{definition}

  Given two sets $P, Q \subseteq \NN^2$, we say that
  $P$ is \textit{northeast equivalent to $Q$},
  written as $P \sim Q$, if  $\nea{P} = \nea{Q}$.
\end{definition}

\begin{definition}
  Let $G$ be a graph on $n$ vertices, and let $(\bl,\re)$ be an ordered
  pair of integers.  We say that $(\bl,\re)$ is a
  \textit{northeast color vector of $G$} if $\bl+\re \le n$ and if
  $\bl \ge \bl_0$ and $\re \ge \re_0$ for some color vector $(\bl_0,\re_0)$ of
  $G$.
\end{definition}

Note that the set of all northeast color vectors of $G$ is
$\uptrunc{n}{\cC(G)}$.  The term \textit{northeast color vector} is
actually a synonym for \textit{elementary inertia,} as we now
demonstrate.

\begin{proposition}
\label{P_elementary}
  Let $G$ be a graph on $n$ vertices.
  Then $\cE(G) = \uptrunc{n}{\cC(G)}$.
\end{proposition}

\begin{proof}
  We show both inclusions.

  \noindent \textit{Forward inclusion.}
  Let $(r,s)$ be an elementary inertia of $G$.  Then there exist
  a nonnegative integer $k$ and an ordered pair of integers $(r_0, s_0)$ such that
  \[
    k \le r_0 \le r, \ k \le s_0 \le s, \mbox{ and } r_0 + s_0 = n - \MD{k}{G} + k.
  \]
  Let $S$ be chosen such that $|S| = k$ and $G-S$ has $\MD{k}{G}$
  components, and let $F$ be a spanning forest of $G-S$,
  so that $F$ has $n-k$ vertices and $r_0 + s_0 - 2k$ edges.
  We partition the edges of $F$ into two sets $\B$ and $\R$ with $r_0-k$
  and $s_0-k$ edges respectively.
  It follows that $(k+|\B|,k+|\R|) = (r_0, s_0)$ is a color vector of
  $G$.  Since $r+s \le n$, $(r,s)$
  belongs to the set $\uptrunc{n}{\cC(G)}$ of northeast color
  vectors of $G$.

  \noindent \textit{Reverse inclusion.}
  Let $(\bl, \re)$ be a northeast color vector of $G$,
  and let $(S,\B,\R)$ be a \BR span of $G$ such that
  $\bl_0 = |S|+|\B| \le \bl$ and $\re_0 = |S|+|\R| \le \re$.
  Letting $k = |S|$, we can assume without loss of generality
  that $S$ is chosen among all sets of size $k$ in such a way as to
  minimize $|\B|+|\R|$, or in other words that $G-S$ has
  $\MD{k}{G}$ components.  Under this assumption
  we have $|\B| + |\R| + \MD{k}{G} = n-k$,
  so $n - \MD{k}{G} + k = \bl_0 + \re_0 \le \bl+\re \le n$.
  We further have $k \le \bl_0 \le \bl$ and $k \le \re_0 \le \re$, so
  $(\bl,\re)$ is an elementary inertia of $G$.
\end{proof}

We now state some set-theoretic results that allow us to simplify
certain expressions involving $\nea{Q}$ and $\bigl[ Q \bigr]_n$.
\begin{observation}
\label{Ob_settheoretic}
  For $Q \subseteq \NN^2$ and nonnegative integers $m \le n$, we
  have
  \begin{enumerate}
  \item
  \label{line1}
    $\Bigl[\bigl[Q\bigr]_n\Bigr]_m
    = \Bigl[\bigl[Q\bigr]_m\Bigr]_n
    = \bigl[Q\bigr]_m$.
  \item
  \label{line1b}
    $\bigl[ \nea{ \uptrunc{n}{Q} \!\!} \bigr]_m = \uptrunc{m}{Q}$.
  \item
  \label{line2}
    $\uptrunc{m}{Q} \sim \bigl[Q\bigr]_m$.
  \item
  \label{line3}
    If $P$ is a stripe of rank $m$, then
    $\bigl[ Q + P \bigr]_n = \bigl[ Q \bigr]_{n-m} + P$.
  \item
  \label{line4}
  $\NN^2_m \subseteq Q$ implies
    $Q \sim \bigl[Q\bigr]_m$.
  \item
  \label{line5}
    $Q \sim \bigl[Q\bigr]_m$ implies
  $Q \sim \bigl[Q\bigr]_n$.
  \end{enumerate}
\end{observation}
\begin{proof}
  These are all straightforward consequences of the definitions.
\end{proof}
\begin{proposition}
  \label{P_settheoretic}
  Let $\ell$, $m$, and $n$ be nonnegative integers with $0 \le \ell \le n$
  and $0 \le m \le n$,
  suppose that $Q \subseteq \NN^2$ satisfies
  $Q \sim \bigl[Q\bigr]_{n-\ell}$,
  and let $P = \uptrunc{n}{Q}$.
  Then
  \begin{enumerate}
  \item
  \label{setline1}
  $P \sim \bigl[P\bigr]_{n-\ell}$,
  \item
  \label{setline2}
  $\uptrunc{m}{P} = \bigl[P\bigr]_{m}$, \ and
  \item
  \label{setline3}
  $\uptrunc{n}{P} = P$.
  \end{enumerate}
\end{proposition}
\begin{proof}
  We have
  \[
    P
    = \uptrunc{n}{Q}
    \sim \bigl[Q\bigr]_n
    \sim Q
    \sim \bigl[Q\bigr]_{n-\ell}
    \sim \uptrunc{n-\ell}{Q}
    = \Bigl[ \uptrunc{n}{Q} \Bigr]_{n-\ell}
    = \bigl[P\bigr]_{n-\ell},
  \]
  \[
    \uptrunc{m}{P} = \bigl[\nea{\uptrunc{n}{Q}\!\!}\bigr]_{m} =
    \uptrunc{m}{Q} = \bigl[\uptrunc{n}{Q}\bigr]_{m} =
    \bigl[P\bigr]_{m},
  \]
  and
  \[
    \uptrunc{n}{P} = \bigl[\nea{\uptrunc{n}{Q}\!\!}\bigr]_{n} =
    \uptrunc{n}{Q} = P.
  \]
\end{proof}
We can apply this proposition immediately.  First note that
Observations~\ref{Ob_cvec_stripe}
and~\ref{Ob_settheoretic}~(\ref{line4}) give us
\begin{observation}
  \label{Ob_c_equiv}
  Let $G$ be a graph on $n$ vertices with $\ell$ components.
  Then $\cC(G) \sim \bigl[\cC(G)\bigr]_{n-\ell}$.
\end{observation}
Observation~\ref{Ob_c_equiv} and Proposition~\ref{P_elementary}
allow us to apply Proposition~\ref{P_settheoretic}, by substituting
$\cC(G)$ for $Q$.
\begin{observation}[Northeast equivalence I]
  \label{Ob_e_equiv}
  Let $G$ be a graph on $n$ vertices with $\ell$ components,
  and let $m$ be an integer in the range $0 \le m \le n$.
  Then
  \begin{enumerate}
  \item
  \label{eqline1}
  $\cE(G) \sim \bigl[\cE(G)\bigr]_{n-\ell}$,
  \item
  \label{eqline2}
  $\uptrunc{m}{\cE(G)} = \bigl[\cE(G)\bigr]_{m}$, \ and
  \item
  \label{eqline3}
  $\uptrunc{n}{\cE(G)} = \cE(G)$.
  \end{enumerate}
\end{observation}
Observation~\ref{Ob_e_equiv} (\ref{eqline3}) can be viewed as a
Northeast Lemma for elementary inertias.

\begin{lemma}
  \label{L_settheoretic}
  Let $Q_1$, \ldots, $Q_k$ be subsets of $\NN^2$,
  and suppose that for some collection $n_1, \ldots, n_k$
  of nonnegative integers we have
  $Q_i \sim \bigl [ Q_i \bigr ]_{n_i}$
  for $i = 1, \ldots, k$.
  Let $Q = Q_1 + \cdots + Q_k$ and
  let $n = n_1 + \cdots + n_k$.
  Then
  \[
    \uptrunc{n}{Q}
    =
    \bigl[ \nea{Q_1} + \cdots + \nea{Q_k} \bigr]_n
    =
    \uptrunc{n_1}{Q_1}
     + \cdots +
    \uptrunc{n_k}{Q_k}.
  \]
\end{lemma}
\begin{proof}
  The first equality comes from
  the observation that $\NN^2 + \NN^2 = \NN^2$.

  For the second equality, the reverse inclusion is easy to check.
  Suppose then that we are given
  \[
    (x,y) \in \bigl[ \nea{Q_1} + \nea{Q_2} + \cdots + \nea{Q_k}
    \bigr]_n,
  \]
  so there exist $k$ ordered pairs of integers
  $(x_i,y_i) \in \nea{Q_i}$
  with $x = \sum\limits_{i=1}^k x_i$,
  $y = \sum\limits_{i=1}^k y_i$, and $x+y \le n$.
  For any such collection $\{(x_i,y_i)\}$, we can define
  two quantities, a \textit{surplus}
  \[
    s = \sum\limits_{i=1}^m \max (x_i + y_i - n_i,\ 0 )
  \]
  and a \textit{deficit}
  \[
    d = \sum\limits_{i=1}^m \max (n_i - x_i - y_i,\ 0 ),
  \]
  so that $x + y - s + d = n$ and hence
  $s \le d$.
  If $s = 0$, then
  in every case
  we have
  $(x_i, y_i) \in \uptrunc{n_i}{Q_i}$, so
  \[
    (x,y) \in \uptrunc{n_1}{Q_1} + \uptrunc{n_2}{Q_2} + \cdots +
    \uptrunc{n_k}{Q_k}
  \]
  and we are done.
  But we can assume $s=0$
  without loss of generality
  for the following reason:
  If $s > 0$, then $d > 0$ also and for some integers
  $i$ and $j$ in the range $1 \le i,j \le k$
  we have $x_i + y_i > n_i$ and $x_j + y_j < n_j$.
  Since $\nea{Q_i} = \nea{\bigl [ Q_i \bigr ]_{n_i}\!\!}$,
  we can replace $(x_i,y_i)$ by
  either $(x_i-1,y_i)$ or $(x_i,y_i-1)$, one
  of which must belong to $\nea{Q_i}$, and
  simultaneously replace $(x_j,y_j)$
  with respectively either $(x_j+1,y_j) \in \nea{Q_j}$ or
  $(x_j,y_j+1) \in \nea{Q_j}$.
  This reduces both the value of $s$ and the value of $d$,
  so we can assume without loss of generality that $s=0$,
  giving the desired result.
\end{proof}
The following proposition is an immediate corollary.
\begin{proposition}
  \label{Cor_settheoretic}
  Given $Q \subseteq \NN^2$ and nonnegative integers $m \le n$,
  suppose that $Q \sim \bigl[ Q \bigr]_m$.  Then
    \[ \uptrunc{n}{Q} = \uptrunc{m}{Q} + \sum\limits_{i=1}^{n-m}
    \NNone. \]
\end{proposition}
\begin{proof}
  Apply Lemma~\ref{L_settheoretic} with $k = n-m+1$, $Q_1 = Q$, $n_1 = m$,
  and for $i > 1$, $Q_i = \{(0,0)\}$ and $n_i = 1$.
  We have abbreviated $\uptrunc{1}{\{(0,0)\}}$ by the equivalent
  expression $\NNone$.
\end{proof}

With the necessary set-theoretic tools in place, we can proceed to
demonstrate some properties of $\cE(G)$, starting with the fact that
it is additive on the connected components of $G$.

\begin{proposition}[Additivity on components]
  \label{P_21e}
  Let $ G=\bigcup\limits_{i=1}^k G_i $.
  Then
  \[
  \cE(G)=\cE(G_1)+\cE(G_2)+\cdots+\cE(G_k).
  \]
\end{proposition}
\begin{proof}
  We first observe that for any \BR span $(S,\B,\R)$ of $G$,
  each entry of the triple is a disjoint union of corresponding
  entries from \BR spans of the components $G_i$, so
  \[
  \cC(G)=\cC(G_1)+\cC(G_2)+\cdots+\cC(G_k).
  \]
  Now let $n = |G|$ and for each integer $i$ in the range $1 \le i \le k$,
  let $n_i = |G_i|$.
  From Observations~\ref{Ob_c_equiv} and~\ref{Ob_settheoretic} (\ref{line5})
  we can conclude that
  $\cC(G_i) \sim \bigl[\cC(G_i)\bigr]_{n_i}$.
  Since $n$ = $n_1 + n_2 + \cdots + n_k$, we can
  apply Lemma~\ref{L_settheoretic} to obtain
  \[
    \uptrunc{n}{\cC(G)} = \uptrunc{n_1}{\cC(G_1)} +
    \uptrunc{n_2}{\cC(G_2)} + \cdots + \uptrunc{n_k}{\cC(G_k)} ,
  \]
  which by Proposition~\ref{P_elementary} is equivalent to
  the desired conclusion.
\end{proof}

Before stating and proving the cut vertex formula for elementary
inertia sets, it will be useful to split the set $\cE(G)$ into two
specialized sets depending on a choice of vertex $v$, and establish
some of the properties of these sets.

\begin{definition}
  Let $G$ be a graph and let $v$ be a vertex of $G$.
  \begin{itemize}
  \item
  If $(S,\B,\R)$ is a \BR span of $G$ and $v \in S$, then
  we say that the ordered pair $(|S|+|\B|,|S|+|\R|)$
  is a \textit{$v$-deleting color vector} of $G$.
  The set of $v$-deleting color vectors of $G$ is
  denoted $\Cdel(G)$.
  \item
If $(S,\B,\R)$ is a \BR span of $G$ and $v \not \in S$, then
  we say that the ordered pair $(|S|+|\B|,|S|+|\R|)$
  is a \textit{$v$-keeping color vector} of $G$.
  The set of $v$-keeping color vectors of $G$ is
  denoted $\Ccon(G)$.
  \end{itemize}
\end{definition}

\begin{definition}
\label{D_eplusminus}
  Let $G$ be a graph on $n$ vertices including $v$.
  We define the set of
  \textit{$v$-deleting elementary inertias of $G$} as
  \[
    \Edel(G) = \uptrunc{n}{\Cdel(G)}
  \]
  and the set of
  \textit{$v$-keeping elementary inertias of $G$} as
  \[
    \Econ(G) = \uptrunc{n}{\Ccon(G)}.
  \]
\end{definition}

The first result we need is an immediate consequence of these
definitions.
\begin{proposition}[Splitting at $v$]
  \label{P_L1}
  Let $G$ be a graph with $v \in V(G)$.
  Then
  \[
    \cE(G) = \Edel(G) \cup \Econ(G).
  \]
\end{proposition}

There are equivalent, simpler expressions for the set of
$v$-deleting color vectors and $v$-deleting elementary inertias of
$G$.
\begin{proposition}[The $v$-deleting formula]
  \label{P_del_equiv}
  Let $G$ be a graph on $n \ge 2$ vertices with $v \in V(G)$.
  Then
  \[
    \Cdel(G) = \cC(G-v) + \{(1,1)\}
  \]
  and
  \[
    \Edel(G) = \bigl[\cE(G-v)\bigr]_{n-2} + \{(1,1)\} = \bigl[\cE(G-v) + \{(1,1)\}\bigr]_n.
  \]
\end{proposition}
\begin{proof}
  The triple $(S,\B,\R)$ is a \BR span of $G$
  with $v \in S$ if and only if the triple
  $(S-\{v\},\B,\R)$ is a \BR span of $G-v$.
  It follows that the $v$-deleting color vectors
  $(r,s)$ in $\Cdel(G)$ are exactly the vectors
  $(1+\bl,1+\re)$ where $(\bl,\re)$ is a color vector
  of $G-v$.  This gives us our first conclusion
  \[
    \Cdel(G) = \cC(G-v) + \{(1,1)\}.
  \]

  With the first conclusion as our starting point, we now have
  \[
    \Edel(G) = \uptrunc{n}{\Cdel(G)}
    = \bigl[ \nea{\cC(G-v)} + \{(1,1)\} \bigr]_{n}.
  \]
  Since $\{(1,1)\}$ is a stripe of rank $2$, by Observation~\ref{Ob_settheoretic}
  this simplifies to
  \begin{eqnarray*}
    \Edel(G) &=&
    \uptrunc{n-2}{\cC(G-v)} + \{(1,1)\} \\
    &=& \Bigl[\uptrunc{n-1}{\cC(G-v)}\Bigr]_{n-2} + \{(1,1)\} \\
    &=& \bigl[\cE(G-v)\bigr]_{n-2} + \{(1,1)\} \\
    &=& \bigl[\cE(G-v) + \{(1,1)\}\bigr]_n
  \end{eqnarray*}
    which completes the proof.
\end{proof}
\begin{observation}
\label{Ob_c_minus_equiv}
  Let $G$ be a graph whose $n$ vertices include $v$, and let $\ell$ be the
  number of components of $G-v$.  Then
  $\Cdel(G) \sim \bigl[\Cdel(G)\bigr]_{n+1-\ell}$.
\end{observation}
\begin{proof}
  By Observation~\ref{Ob_c_equiv},
  $\cC(G-v) \sim \bigl[ \cC(G-v) \bigr]_{n-1-\ell}$.
  Proposition~\ref{P_del_equiv}
  and Observation~\ref{Ob_settheoretic} (\ref{line3})
  then give us
  $\Cdel(G) \sim \bigl[\Cdel(G)\bigr]_{n+1-\ell}$.
\end{proof}
Substituting $Q = \Cdel(G)$ into Proposition~\ref{P_settheoretic}
now gives us a result about $v$-deleting elementary inertias.
\begin{observation}[Northeast equivalence II]
  \label{Ob_e_minus_equiv}
  Let $G$ be a graph whose $n$ vertices include $v$, let $\ell$ be the
  number of components of $G-v$, and let $m$ be an integer
  in the range $0 \le m \le n$.
  Then
  \begin{enumerate}
  \item
  \label{eqminusline1}
    $\Edel(G) \sim \bigl[\Edel(G)\bigr]_{n+1-\ell}$,
  \item
  \label{eqminusline2}
    $\uptrunc{m}{\Edel(G)} = \bigl[\Edel(G)\bigl]_m$, \ and
  \item
  \label{eqminusline3}
    $\uptrunc{n}{\Edel(G)} = \Edel(G)$.
  \end{enumerate}
\end{observation}
Similar results hold for the $v$-keeping color vectors and
$v$-keeping elementary inertias:
\begin{observation}
\label{Ob_c_plus_equiv}
  Let $G$ be a graph whose $n$ vertices include $v$, and let $\ell =
  \MD{0}{G}$.  Then
  $\Ccon(G) \sim \bigl[\Ccon(G)\bigr]_{n-\ell}$.
\end{observation}
\begin{proof}
  It suffices to consider \BR spans of the form $(\emptyset,\B,\R)$,
  which of course satisfy $v \not \in \emptyset$.  The set of $v$-keeping color
  vectors arising from such \BR spans is exactly $\NN^2_{n-\ell}$,
  from which the desired result follows by
  Observation~\ref{Ob_settheoretic} (\ref{line4}).
\end{proof}
  Proposition~\ref{P_settheoretic} now gives us:
\begin{observation}[Northeast equivalence III]
\label{Ob_e_plus_equiv}
  Let $G$ be a graph whose $n$ vertices include $v$, let $\ell =
  \MD{0}{G}$, and let $m$ be an integer in the range $0 \le m \le n$.
  Then
  \begin{enumerate}
  \item
  \label{eqplusline1}
    $\Econ(G) \sim \bigl[\Econ(G)\bigr]_{n-\ell}$,
  \item
  \label{eqplusline2}
    $\uptrunc{m}{\Econ(G)} = \bigl[\Econ(G)\bigr]_m$, \ and
  \item
  \label{eqplusline3}
    $\uptrunc{n}{\Econ(G)} = \Econ(G)$.
  \end{enumerate}
\end{observation}

It is possible to restrict the set of allowable \BR spans that
define $\Ccon(G)$ and still obtain the full set of $v$-keeping color
vectors of $G$.
\begin{proposition}
  \label{P_con_equiv}
  Let $G$ be a graph with vertex $v$, and let
  $E^\prime = E(G-v)$.
  Suppose that $(\bl,\re)$ belongs to $\Ccon(G)$.
  Then there exists a \BR span
  $(S,\B,\R)$ of $G$ with $v \not \in S$
  such that $(\bl,\re) = (|S|+|\B|,|S|+|\R|)$ and
  such that $(S,\B\cap E^\prime,\R\cap
  E^\prime)$ is a \BR span of $G-v$.
\end{proposition}
\begin{proof}
  By the definition of $\Ccon(G)$,
  there exists a \BR span
  $(S,\B,\R)$ of $G$ with $v \not \in S$
  such that $(\bl,\re) = (|S|+|\B|,|S|+|\R|)$.
  The vertex $v$ thus belongs to some component $G_i$ of
  $G-S$, and those edges in $\B$ and $\R$ which are part of $G_i$
  give a spanning tree $T_i$ of $G_i$.
  There is no loss of generality if we assume that $T_i$
  is constructed as follows:
  First, a spanning tree is obtained for
  each component of $G_i - v$.
  Each subtree is then connected to $v$
  by way of a single edge, so that the degree of $v$
  in $T_i$ is equal to $\MD{0}{G_i - v}$.
  With this assumption, $(S, \B\cap E^\prime, \R\cap
  E^\prime)$ is a \BR span of $G-v$.
\end{proof}

The next key ingredient is a consequence of
Propositions~\ref{P_del_equiv} and~\ref{P_con_equiv}.
\begin{proposition}[Domination by $G-v$]
  \label{P_L2}
  Let $G$ be a graph on $n$ vertices, one of which is $v$.
  Then for $\epsilon \in \{\del,\con\}$ we have
  \[
    \bigl[\cE_v^\epsilon(G)\bigr]_{n-1} \subseteq \cE(G-v).
  \]
\end{proposition}

Given Proposition~\ref{P_L1}, Proposition~\ref{P_L2} is equivalent
to an inclusion on elementary inertia sets which has already been
proven for inertia sets as Proposition~\ref{P22} (a):
\begin{proposition}
  \label{P22e}
  For any graph $G$ and any vertex $v \in V(G)$,
  \[
    \bigl[\cE(G)\bigr]_{n-1} \subseteq \cE(G-v).
  \]
\end{proposition}

We need one more result before stating and proving the cut vertex
formula for elementary inertias.

\begin{proposition}[The $v$-keeping cut vertex formula]
  \label{P_L3}
  Let $G=\bigvsum\limits_{i=1}^k G_i$ be a graph on $n$ vertices
  which is a vertex sum of graphs $ G_1,G_2,\ldots,G_k $ at $ v $,
  for $k \ge 2$.
  Then
  \[
    \Econ(G)
    =
    \bigl[\Econ(G_1) + \Econ(G_2) + \cdots + \Econ(G_k)\bigr]_n.
  \]
\end{proposition}
\begin{proof}
  Let $G$, $v$, $n$, and $G_1, \ldots, G_k$ be as in the statement of
  the proposition.
  We first establish a related identity,
  \[
    \Ccon(G)
    =
    \Ccon(G_1) + \Ccon(G_2) + \cdots + \Ccon(G_k).
  \]
  This holds because
  \begin{enumerate}
  \item
    The sets $V(G_i) - \{v\}$ are disjoint,
    and their union is $V(G) - \{v\}$,
    so subsets $S \subseteq V(G)$ with $v \not \in S$
    are in bijective correspondence with collections
    of subsets $S_i \subseteq V(G_i)$ none of which
    contain $v$.
  \item
    For any such set $S$ partitioned as a union of $S_i$,
    $G-S$ is a vertex sum at $v$ of the graphs
    $G_i-S_i$, and so
    the set $E(G-S)$ is a disjoint union
    of $E(G_i - S_i)$.
  \item
    A subgraph $F$ of the vertex sum $G-S$ is
    a spanning forest of $G-S$ if and only if $F$
    is a vertex sum of graphs $F_i$
    each of which is a spanning forest of $G_i-S_i$.
  \end{enumerate}
  For each graph $G_i$, let $n_i = |G_i|$, so that
  $(n-1) = \sum\limits_{i=1}^k(n_i-1)$.
  Since each graph $G_i$ contains the vertex $v$,
  $\MD{0}{G_i} \ge 1$.  Observations~\ref{Ob_c_plus_equiv}
  and~\ref{Ob_settheoretic} (\ref{line5}) then give us
  $\Ccon(G_i) \sim \bigl[\Ccon(G_i)\bigr]_{n_i-1}$.
  Thus by Lemma~\ref{L_settheoretic} we have
  \[
    \uptrunc{n-1}{\Ccon(G)}
    =
    \uptrunc{n_1-1}{\Ccon(G_1)} +
    \cdots
    + \uptrunc{n_k-1}{\Ccon(G_k)}.
  \]
  We also have $\Ccon(G) \sim \bigl[\Ccon(G)\bigr]_{n-1}$,
  so by Proposition~\ref{Cor_settheoretic} we can add
  $k$ copies of $\NNone$ to both sides to obtain
  \[
    \uptrunc{n-1+k}{\Ccon(G)}
    =
    \uptrunc{n_1}{\Ccon(G_1)} +
    \cdots
    + \uptrunc{n_k}{\Ccon(G_k)}
  \]
  which gives the desired formula by
  Observation~\ref{Ob_settheoretic} (\ref{line1})
  and Definition~\ref{D_eplusminus}.
\end{proof}

The proof of the cut vertex formula depends on the following
properties of $\cE(G)$, $\Econ(G)$, and $\Edel(G)$:
\begin{itemize}
\item
  Northeast equivalence I and III (Observations~\ref{Ob_e_equiv} and~\ref{Ob_e_plus_equiv}),
\item
  Additivity on components (Proposition~\ref{P_21e}),
\item
  Splitting at $v$ (Proposition~\ref{P_L1}),
\item
  The $v$-deleting formula (Proposition~\ref{P_del_equiv}),
\item
  Domination by $G-v$ (Proposition~\ref{P_L2}), and
\item
  The $v$-keeping cut vertex formula (Proposition~\ref{P_L3}).
\end{itemize}

\begin{theorem}[The cut vertex formula for elementary inertias]
  \label{Th14e} Let $G$  be a graph on $ n\ge 3 $ vertices and let
  $v$ be a cut vertex of~$G$.
  Write $ G=\bigvsum\limits_{i=1}^k G_i$, $k \ge 2$,
  the vertex sum of \ $ G_1,G_2,\ldots,G_k $ at $ v $.  Then
\begin{eqnarray*}
\cE(G) &=&\bigl[\cE(G_1)+\cE(G_2)+\cdots+\cE(G_k)\bigr]_n
\label{eq5e}
\\
&&\cup\
\bigl[\cE(G_1-v)+\cE(G_2-v)+\cdots+\cE(G_k-v)+\{(1,1)\}\bigr]_n.
\nonumber
\end{eqnarray*}
\end{theorem}

\begin{proof}
  We manipulate both sides to obtain the same set.

  Define two sets
  \[
    \Qdel = \bigl[ \cE(G_1-v) + \cdots + \cE(G_k-v) + \{(1,1)\} \bigr]_n
  \]
  and
  \[
    \Qcon = \bigl[ \Econ(G_1) + \cdots + \Econ(G_k) \bigr]_n.
  \]

  By Propositions~\ref{P_del_equiv} and~\ref{P_21e},
  $\Edel(G) = \Qdel$
  and by Proposition~\ref{P_L3},
  $\Econ(G) = \Qcon$,
  so by Proposition~\ref{P_L1}, $\cE(G) = \Qdel \cup \Qcon$.

  The right hand side is
  \[
    \mathrm{RHS} = \bigl[ \cE(G_1) + \cdots + \cE(G_k) \bigr]_n \cup \Qdel.
  \]
  For each $i = 1, \ldots, k$, let $n_i = |G_i|$, so that
  $\cE(G_i) \sim \bigl[ \cE(G_i) \bigr]_{n_i-1}$
  (Observations~\ref{Ob_e_equiv} (\ref{eqline1}) and~\ref{Ob_settheoretic} (\ref{line5}), since
  in each case $\ell \ge 1$).
  Starting with Observation~\ref{Ob_e_equiv} (\ref{eqline3}) and
  then applying Lemma~\ref{L_settheoretic} both backwards and forwards, we have
  \begin{eqnarray*}
    \bigl[ \cE(G_1) + \cdots + \cE(G_k) \bigr]_n
    &=&
    \Bigl[ \uptrunc{n_1}{\cE(G_1)} + \cdots + \uptrunc{n_k}{\cE(G_k)}
     \Bigr]_n \\
    &=&
    \Bigl[ \bigl[\nea{\cE(G_1)} + \cdots + \nea{\cE(G_k)}
    \bigr]_{n-1+k} \Bigr]_n \\
    &=&
    \bigl[\nea{\{(0,0)\}} + \nea{\cE(G_1)} + \cdots + \nea{\cE(G_k)}
    \bigr]_{n} \\
    &=&
    \NNone + \uptrunc{n_1-1}{\cE(G_1)} + \cdots
      + \uptrunc{n_k-1}{\cE(G_k)}
  \end{eqnarray*}
  which by Observation~\ref{Ob_e_equiv} (\ref{eqline2}) gives us
  \[
    \mathrm{RHS} = \Qdel \cup \Bigl(
      \NNone +
      \bigl[\cE(G_1)\bigr]_{n_1-1} +
      \cdots +
      \bigl[\cE(G_k)\bigr]_{n_k-1}
    \Bigr).
  \]
  By applying Proposition~\ref{P_L1} to each term
  $\bigl[\cE(G_i)\bigr]_{n_i-1}$ we obtain
  \[
    \mathrm{RHS} =
    \Qdel \cup
    \Bigl (
      \NNone +
      \sum\limits_{i=1}^k
      \bigl (
        \bigl[\Edel(G_i)\bigr]_{n_i-1} \cup
        \bigl[\Econ(G_i)\bigr]_{n_i-1}
      \bigr )
    \Bigr ) .
  \]
  For any $\alpha = (\epsilon_1, \ldots, \epsilon_m) \in \{\del,\con\}^k$
  we will define
  \[
    \cE_v^\alpha = \sum\limits_{i=1}^k
    \bigl[\cE_v^{\epsilon_i}(G_i)\bigr]_{n_i-1}.
  \]
  This gives us
  \[
    \mathrm{RHS} = \!\!
    \bigcup\limits_{\alpha \in \{\del,\con\}^k} \!\!
      \Qdel \cup \bigl( \NNone + \cE_v^\alpha\bigr).
  \]
  We divide the $2^k$ choices for $\alpha$ into two cases:
  either $\epsilon_j$ is  ``$\del$''
  for some $j \in \{1, \ldots, k\}$, or
  $\epsilon_i$ is ``$\con$'' for all $i$.
  In the first case, by Proposition~\ref{P_del_equiv} we have
  \begin{eqnarray*}
    \NNone + \cE_v^\alpha
    &=&
    \NNone +
    \bigl[\cE_v^{\epsilon_1}(G_1)\bigr]_{n_1-1} + \cdots \\
    &&\cdots +
      \bigl[\cE(G_j - v) + \{(1,1)\} \big]_{n_j-1}
    + \cdots \\
    &&\cdots + \bigl[\cE_v^{\epsilon_k}(G_k)\bigr]_{n_k-1}.
  \end{eqnarray*}
  We wish to show that this is a subset of $\Qdel$.
  For every $i$ besides $j$, we have
  $\bigl[\cE_v^{\epsilon_i}(G_i)\bigr]_{n_i-1} \subseteq \cE(G_i-v)$
  by Proposition~\ref{P_L2}.
  The remaining terms we regroup as
  \begin{eqnarray*}
    \NNone +
    \bigl[\cE(G_j - v) + \{(1,1)\} \big]_{n_j-1} \!
    &=& \!
    \NNone + \bigl[\cE(G_j - v)\bigr]_{n_j-3} + \{(1,1)\} \\
    &  \subseteq & \! \NNone + \bigl[\cE(G_j - v)\bigr]_{n_j-2} +
    \{(1,1)\}.
  \end{eqnarray*}
  Observation~\ref{Ob_e_equiv} and
  Proposition~\ref{Cor_settheoretic} give us
  \[
    \NNone + \bigl[\cE(G_j - v)\bigr]_{n_j-2} = \cE(G_j - v).
  \]
  We have thus shown that
  \[
    \NNone + \cE_v^\alpha \, \subseteq \, \cE(G_1-v) + \cdots + \cE(G_k-v) +
    \{(1,1)\},
  \]
  and since
  \[
    \NNone + \cE_v^\alpha  = \bigl[ \NNone + \cE_v^\alpha \bigr]_n,
  \]
  this gives us $\Qdel \cup \bigl( \NNone + \cE_v^\alpha \bigr) =
  \Qdel$ in the case where $\alpha$ has at least one sign $\epsilon_j =
  $ ``$\del$''.

  This leaves the case where $\alpha$ has all signs $\epsilon_j = $ ``$\con$''.
  By Observations~\ref{Ob_e_plus_equiv}~(\ref{eqplusline1})
  and~\ref{Ob_settheoretic}~(\ref{line5}), $\Econ(G_i) \sim \bigl[ \Econ(G_i) \bigr]_{n_i-1}$.
  Starting with Observation~\ref{Ob_e_plus_equiv}~(\ref{eqplusline2}),
  applying Lemma~\ref{L_settheoretic}
  both backwards and forwards, and finally using
  Observation~\ref{Ob_e_plus_equiv}~(\ref{eqplusline3}), we have
  \begin{eqnarray*}
    \NNone + \cE_v^\alpha
    &=&
    \NNone + \uptrunc{n_1-1}{\Econ(G_1)} + \cdots
      + \uptrunc{n_k-1}{\Econ(G_k)} \\
    &=&
    \bigl[\nea{\{(0,0)\}} + \nea{\Econ(G_1)} + \cdots + \nea{\Econ(G_k)}
    \bigr]_{n} \\
    &=&
    \Bigl[ \bigl[\nea{\Econ(G_1)} + \cdots + \nea{\Econ(G_k)}
    \bigr]_{n-1+k} \Bigr]_n \\
    &=&
    \Bigl[ \uptrunc{n_1}{\Econ(G_1)} + \cdots + \uptrunc{n_k}{\Econ(G_k)} \Bigr]_n \\
    &=& \Qcon.
  \end{eqnarray*}
  The entire union thus collapses to
  $
    \mathrm{RHS} =
    \Qdel \cup \Qcon,
  $
  and the left and right hand expressions are equal.
\end{proof}

\begin{remark}
  We can generalize the splitting of $\cE(G)$ into
  $\Econ(G)$ and $\Edel(G)$ for non-elementary inertias:
  Given a graph $G$ with vertex $v$ and $A \in \Herm(G)$,
  order the vertices of $G$ such that $v=1$ and decompose
  $A$ as
  \[
  A = \left [
    \begin{array}{cc}
       a_{11} & b^* \\
       b & B
    \end{array}
    \right ].
  \]
  If $b$ is in the column space of $B$, then
  say that $\pin(A) \in \hcI_v^\con(G)$, and define
  $\hcI_v^\del(G)$ as $\bigl[\hcI(G-v)+\{(1,1)\}\bigr]_n$.
  Define $\cI_v^\con(G)$ and $\cI_v^\del(G)$ analogously.
  Under these definitions we can uniformly replace
  $\cE$ with $\hcI$ or $\cI$
  in Observations~\ref{Ob_e_equiv},
  and~\ref{Ob_e_plus_equiv}
  and in
  each of Propositions~\ref{P_21e}, \ref{P_L1},
  \ref{P_del_equiv}, \ref{P_L2}, and~\ref{P_L3},
  and we claim that in every case the result still holds.
  We will not prove these statements,
  as we already have a proof of Theorem~\ref{Th13},
  but given those observations and propositions,
  the proof of Theorem~\ref{Th14e}
  demonstrates the same cut vertex formula
  for inertia sets.
\end{remark}

We now state and prove the main result of the section.
\begin{theorem}
\label{Th15}
For any tree $ T $, $ \cI(T)=\cE(T) $.
\end{theorem}
\begin{proof}
  Let $ n=|T| $.
\\
  If~$ n=1 $, \ $ T=K_1 $ and $ \cI(T)=\NN^2_{[0,1]} $.
  Since $(\emptyset,\emptyset,\emptyset) $ is a \BR span of $ K_1 $,
  the origin $(0,0)$ is a color vector of $T$ and
  $ \cE(T)=\NN^2_{[0,1]} $ also.
  If $n=2$, then $T=K_2$, and
  $\cI(K_2)=\NN^2_{[1,2]}=\cE(K_2)$.

  Proceeding by induction, assume that $\cI(T)=\cE(T)$ for all trees
  $ T $ on fewer than $ n $ vertices and let $ T $ be a tree on $ n $
  vertices, $n \ge 3$.
  Let $ v $ be a cut vertex of~$ T $ of degree $ k \ge 2 $.
  Write $ T=\bigvsum\limits_{i=1}^k T_i $, the vertex sum of $
  T_1,\ldots,T_k $ at $ v $. By Theorem~\ref{Th13},
\begin{eqnarray*}
\cI(T) &=& \bigl[\cI(T_1)+\cdots+\cI(T_k)\bigr]_n
\\
&&\cup \ \bigl[\cI(T_1-v)+\cdots+\cI(T_k-v)+\{(1,1)\}\bigr]_n
\end{eqnarray*}
and by Theorem~\ref{Th14e},
\begin{eqnarray*}
\cE(T) &=& \bigl[\cE(T_1)+\cdots+\cE(T_k)\bigr]_n
\\
&&\cup \ \bigl[\cE(T_1-v)+\cdots+\cE(T_k-v)+\{(1,1)\}\bigr]_n.
\end{eqnarray*}
Corresponding terms on the right hand side of these last two
equations are equal by the induction hypothesis, so $ \cI(T)=\cE(T)
$.
\end{proof}
\begin{corollary}
\label{Cor16}
For any forest $ F $, \ $ \cI(F)=\cE(F) $.
\end{corollary}
\begin{proof}
  By Theorem~\ref{Th15}, $\cI(T) = \cE(T)$ for every component
  $T$ of $F$, and by additivity on components for both
  $\cI(G)$ (Observation~\ref{Ob21}) and
  $\cE(G)$ (Proposition~\ref{P_21e}),
  $\cI(F) = \cE(F)$.
\end{proof}
Claim 1 of Theorem~\ref{Th_treemain}, which says $\cI(F) \subseteq
\cE(F)$ for any forest $F$, has now been verified.

We restate Theorem~\ref{Th_treemain} compactly as
\begin{theorem}
  Let $G$ be a graph.  Then $\cI(G) = \cE(G)$ if and only if
  $G$ is a forest.
\end{theorem}


\section{Graphical determination of the inertia set of a tree}
\label{S6}

Tabulating the full inertia set of a tree $T$ on $n$ vertices by
means of Theorem~\ref{Th_treemain} appears, potentially, to require
a lot of calculation:  Every integer $k$ in the range $0 \le k \le
n$ with $\MD{k}{T} \ge k$ gives a trapezoid (possibly degenerate) of
elementary inertias, and the full elementary inertia set is the
union of those trapezoids.  (One could also construct every possible
\BR span of the tree, which is even more cumbersome.)  In fact the
calculation is quite straightforward once we have the first few
values of $\MD{k}{T}$. In this section we present the necessary
simplifications and perform the calculation for a few examples.

\begin{definition}
  \label{D19}
  For any graph $ G $ on $ n $ vertices and $k\in\{0,\ldots,n\} $,
  let
  \[
    \pi_k(G)=\min\left \{r :\ (r,k)\in\cI(G)\right \},
  \]
  \[
    \nu_k(G)=\min\left \{s :\ (k,s)\in\cI(G)\right \}.
  \]
Since $ \pi_k(G)=\nu_k(G) $ for each $k$, we will deal exclusively
with $ \pi_k(G) $.
\end{definition}

The main simplification toward calculating the inertia set of a tree is the following:
\begin{theorem}
  \label{Th17}
  Let $ T $ be a tree on $ n $ vertices and let
  $k\in\{0,1,\ldots,c(T)\}$. Then
\begin{equation}
  \label{eq7}
  \pi_k(T)=n-\MD{k}{T}.
\end{equation}
\end{theorem}
\begin{proof}
  For $k$ in the given range, we can apply Corollary~\ref{Cor_MDmonotone}
  with $j=0$ to obtain $\MD{k}{T} \ge \MD{0}{T} + k$
  and in particular $\MD{k}{T} \ge k$.
  We can thus apply the Stars and Stripes Lemma to obtain
  $(n-\MD{k}{T},k) \in \cI(T)$ and hence $\pi_k(T) \le n-\MD{k}{T}$.
  It remains to prove, for $k \le c(T)$,
  \[ n - \MD{k}{T} \le \pi_k(T). \]
  Suppose by way of contradiction that $(r,s) \in \cI(T)$
  with $s = k$ and $r < n - \MD{k}{T}$.
  By Theorem~\ref{Th_treemain}, every element of $\cI(T)$ is an
  elementary inertia, and thus there is some integer $j$ for which
  $j \le r$, $j \le s$, and $n - \MD{j}{T} + j \le r+s$.
  This implies, for $0 \le j \le k \le c(T)$, that
  \[
    \MD{k}{T}
    <
    \MD{j}{T} + (k-j),
  \]
  which contradicts Corollary~\ref{Cor_MDmonotone}.
\end{proof}
\begin{corollary}
\label{Cor18}
$ \Bigl\{\pi_k(T)\Bigr\}_{k=0}^{\mr(T)-c(T)} $ is a strictly decreasing sequence.
\end{corollary}
\begin{proof}
This follows from Theorem~\ref{Th17},
Corollary~\ref{Cor_MDmonotone} (with $k-j=1$),
and Theorem~\ref{Th9}.
\end{proof}
\begin{theorem}
\label{Th19}
Let $T$ be a tree.
Then $L_T = \cI(T) \cap \NN^2_{\mr(T)}$.
\end{theorem}
In other words, every partial inertia of minimum rank
is in the minimum-rank stripe already defined.
\begin{proof}
  We already have $L_T \subseteq \cI(T)$ by Theorem~\ref{Th9},
  and $L_T \subseteq \NN^2_{\mr(T)}$ by definition.
  To show equality, it suffices by symmetry (Observation~\ref{Ob8})
  to show that for $ k<c(T) $, \
  $ k + \pi_k(T)>\mr(T) $.
  Let $c = c(T)$.
  If $c=0$, we are done.
  By Observation~\ref{Ob_cMD}, $n - \MD{c}{T} + c = \mr(T)$
  but $n - \MD{k}{T} + k > \mr(T)$ for $k < c$.
  It follows by Theorem~\ref{Th17} that $k + \pi_k(T) > \mr(T)$ for $k < c$,
  which completes the proof.
\end{proof}
Theorem~\ref{Th12} already gives a method for determining the
inertia set $\cI(T)$ for any tree $T$, but with
Theorem~\ref{Th_treemain} and the simplifications above there is a
much easier method, which we summarize in the following steps:
\begin{enumerate}
\item   Use the algorithm of Observation~\ref{Ob_Palgorithm} to find $P(T)$.
\item   Since $T$ is connected, $\MD{0}{T} = 1$.
If $T$ is a path then $c(T) = 0$; otherwise $c(T) \ge 1$ and
$\MD{1}{T} = \Delta(T)$.  Continue to calculate higher values
of $\MD{k}{T}$ until $\MD{k}{T} - k = P(T)$, at which point
$k = c(T)$.
\item The defining southwest corners of $\cI(T)$ are $(n-\MD{k}{T}, k)$
and its reflection $(k,n-\MD{k}{T})$, for $0 \le k < c(T)$,
together with the stripe $L_T$ of partial inertias
from $(n-P(T)-c(T),c(T))$ to $(c(T),n-P(T)-c(T))$.
\item Every other point of $\cI(T)$ is a result of the Northeast Lemma
applied to the defining southwest corners.
\end{enumerate}

We give three examples.
\begin{example}
\label{Ex13}
Let $T$ be the tree in Example~\ref{Ex7}, whose inertia set we have already calculated.
\[
\includegraphics{tree2}
\]
Here $P(T) = 2$ and $\mr(T) = 4$.
We have
\[ \MD{1}{T} = 3, \ \MD{1}{T} - 1 = 2,\]
so $c(T) = 1$, and from
$\pi_0(T) = 5$ we go immediately to $L_T$, starting at height $1$,
which is the
convex stripe of three partial inertias from
$(3,1)$ to $(1,3)$.

\setlength{\unitlength}{1mm}
\begin{picture}(50,50)
 \put(42,4){\vector(0,1){40}}
 \put(40,6){\vector(1,0){41}}
 \put(40.75,39){\line(1,0){2.5}}
 \put(42,39){\circle*{1.5}}
 \put(40.75,33.5){\line(1,0){2.5}}
 \put(42,33.5){\circle*{1.5}}
 \put(47.5,33.5){\circle*{1.5}}
 \put(40.75,28){\line(1,0){2.5}}
 \put(47.5,28){\circle*{1.5}}
 \put(53,28){\circle*{1.5}}
 \put(40.75,22.5){\line(1,0){2.5}}
 \put(47.5,22.5){\circle*{1.5}}
 \put(53,22.5){\circle*{1.5}}
 \put(58.5,22.5){\circle*{1.5}}
 \put(40.75,17){\line(1,0){2.5}}
 \put(53,17){\circle*{1.5}}
 \put(58.5,17){\circle*{1.5}}
 \put(64,17){\circle*{1.5}}
 \put(40.75,11.5){\line(1,0){2.5}}
 \put(58.5,11.5){\circle*{1.5}}
 \put(64,11.5){\circle*{1.5}}
 \put(69.5,11.5){\circle*{1.5}}
\put(47.5,4.75){\line(0,1){2.5}}
 \put(53,4.75){\line(0,1){2.5}}
 \put(58.5,4.75){\line(0,1){2.5}}
 \put(64,4.75){\line(0,1){2.5}}
 \put(69.5,4.75){\line(0,1){2.5}}
 \put(69.5,6){\circle*{1.5}}
 \put(75,4.75){\line(0,1){2.5}}
 \put(75,6){\circle*{1.5}}
\end{picture}
\end{example}
\begin{example}
\label{Ex14}
Let $T$ be the tree
\[ \includegraphics{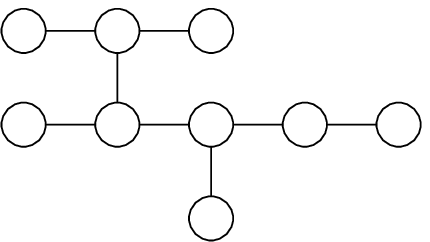} \]
whose horizontal paths realize the path cover
number $P(T) = 3$, so $\mr(T) = 6$.
Taking any vertex of degree $3$ we have
\[ \MD{1}{T} = 3, \  \MD{1}{T} - 1 = 2, \]
and taking the non-adjacent pair of degree-$3$ vertices we have
\[ \MD{2}{T} = 5, \ \MD{2}{T}-2 = 3, \]
so $c(T) = 2$.
Starting as always from $\pi_0(T) = n - 1 = 8$, we
need only one more value $\pi_1(T) = 9 - 3 = 6$ before
reaching the minimum-rank stripe $L_T$ from $(4, 2)$ to $(2,4)$.
The complete set $\cI(T)$ is:

\setlength{\unitlength}{1mm}
\begin{picture}(50,70)
 \put(42,4){\vector(0,1){58}}
 \put(40,6){\vector(1,0){58}}
 \put(40.75,55.5){\line(1,0){2.5}}
 \put(42,55.5){\circle*{1.5}}
 \put(40.75,50){\line(1,0){2.5}}
 \put(42,50){\circle*{1.5}}
 \put(47.5,50){\circle*{1.5}}
 \put(40.75,44.5){\line(1,0){2.5}}
 \put(47.5,44.5){\circle*{1.5}}
 \put(53,44.5){\circle*{1.5}}
 \put(40.75,39){\line(1,0){2.5}}
 \put(47.5,39){\circle*{1.5}}
 \put(53,39){\circle*{1.5}}
 \put(58.5,39){\circle*{1.5}}
 \put(40.75,33.5){\line(1,0){2.5}}
 \put(53,33.5){\circle*{1.5}}
 \put(58.5,33.5){\circle*{1.5}}
 \put(64,33.5){\circle*{1.5}}
 \put(40.75,28){\line(1,0){2.5}}
 \put(53,28){\circle*{1.5}}
 \put(58.5,28){\circle*{1.5}}
 \put(64,28){\circle*{1.5}}
 \put(69.5,28){\circle*{1.5}}
 \put(40.75,22.5){\line(1,0){2.5}}
 \put(58.5,22.5){\circle*{1.5}}
 \put(64,22.5){\circle*{1.5}}
 \put(69.5,22.5){\circle*{1.5}}
 \put(75,22.5){\circle*{1.5}}
 \put(40.75,17){\line(1,0){2.5}}
 \put(64,17){\circle*{1.5}}
 \put(69.5,17){\circle*{1.5}}
 \put(75,17){\circle*{1.5}}
 \put(80.5,17){\circle*{1.5}}
 \put(40.75,11.5){\line(1,0){2.5}}
 \put(75,11.5){\circle*{1.5}}
 \put(80.5,11.5){\circle*{1.5}}
 \put(47.5,4.75){\line(0,1){2.5}}
 \put(53,4.75){\line(0,1){2.5}}
 \put(58.5,4.75){\line(0,1){2.5}}
 \put(64,4.75){\line(0,1){2.5}}
 \put(69.5,4.75){\line(0,1){2.5}}
 \put(75,4.75){\line(0,1){2.5}}
 \put(80.5,4.75){\line(0,1){2.5}}
 \put(86,4.75){\line(0,1){2.5}}
 \put(86,11.5){\circle*{1.5}}
 \put(86,6){\circle*{1.5}}
 \put(91.5,4.75){\line(0,1){2.5}}
 \put(91.5,6){\circle*{1.5}}
\end{picture}
\end{example}

The examples we have shown so far appear to exhibit some sort of convexity.
For $F$ a forest we do at least have convexity of $\cI(F)$ on stripes
of fixed rank, as stated in Corollary~\ref{Cor_stripes}.
Based on small examples one may be led to believe that, in addition,
$\pi_k(T)$ is a convex function
in the range $0 \le k \le c(T)$, or in other words that
$$ \pi_k(T)-\pi_{k+1}(T) \le \pi_{k-1}(T) - \pi_k(T) \quad \text{for}\ 0<k<c(T).$$
However, this is not always the case, as seen in the following
example:
\begin{example}
\label{Ex15}
Given $S_4$ with $v$ a pendant vertex,
let $T$ be the tree constructed as a
vertex sum of four copies of the marked $S_4$:
\[
T= \bigvsum_{i=1}^4 S_4
= \
\raisebox{-.5\totalheight}{
\includegraphics{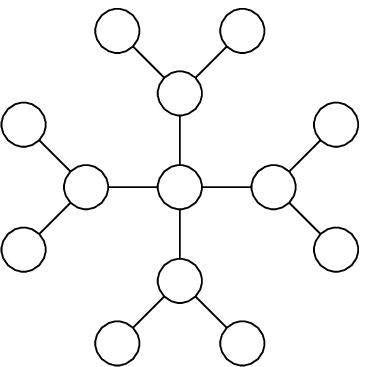}}.
\]
Here $P(T) = 5$ and $\mr(T) = 8$.
To find $\MD{1}{T}$ we always take a vertex of maximum degree;
here
\[ \MD{1}{T} = 4, \ \MD{1}{T} - 1 = 3. \]
For $\MD{2}{T}$ we can either add the center of a branch
or leave out the degree-$4$ vertex and take two centers of branches;
either choice gives us
\[ \MD{2}{T} = 5, \ \MD{2}{T} - 2 = 3. \]
At $k=3$ something odd happens: to remove $3$ vertices and maximize
the number of remaining components, we must not include the
single vertex of maximum degree.  Taking the centers of three branches,
we obtain
\[ \MD{3}{T} = 7, \ \MD{3}{T} - 3 = 4, \]
and finally
taking all four vertices of degree $3$ gives us
\[ \MD{4}{T} = 9, \ \MD{4}{T} - 4 = 5 = P(T),\]
so $c(T) = 4$.
The sequence $\pi_k(T)$ thus starts $(12, 9, 8, 6, 4)$.
As is the case with the stars $S_n$ and Example~\ref{Ex8},
we here have a tree where the minimum-rank stripe $L_T$
is a singleton,
in this case the point $(4, 4)$.
The full plot of $\cI(T)$ is:

\setlength{\unitlength}{1mm}
\begin{picture}(50,90)
 \put(22,4){\vector(0,1){80}}
 \put(20,6){\vector(1,0){80}}
 \put(20.75,77.5){\line(1,0){2.5}}
   \put(22,77.5){\circle*{1.5}}
 \put(20.75,72){\line(1,0){2.5}}
   \put(22,72){\circle*{1.5}}
   \put(27.5,72){\circle*{1.5}}
 \put(20.75,66.5){\line(1,0){2.5}}
   \put(27.5,66.5){\circle*{1.5}}
   \put(33,66.5){\circle*{1.5}}
 \put(20.75,61){\line(1,0){2.5}}
   \put(27.5,61){\circle*{1.5}}
   \put(33,61){\circle*{1.5}}
   \put(38.5,61){\circle*{1.5}}
 \put(20.75,55.5){\line(1,0){2.5}}
   \put(27.5,55.5){\circle*{1.5}}
   \put(33,55.5){\circle*{1.5}}
   \put(38.5,55.5){\circle*{1.5}}
   \put(44,55.5){\circle*{1.5}}
 \put(20.75,50){\line(1,0){2.5}}
   \put(33,50){\circle*{1.5}}
   \put(38.5,50){\circle*{1.5}}
   \put(44,50){\circle*{1.5}}
   \put(49.5,50){\circle*{1.5}}
 \put(20.75,44.5){\line(1,0){2.5}}
    \put(38.5,44.5){\circle*{1.5}}
   \put(44,44.5){\circle*{1.5}}
   \put(49.5,44.5){\circle*{1.5}}
   \put(55,44.5){\circle*{1.5}}
 \put(20.75,39){\line(1,0){2.5}}
   \put(38.5,39){\circle*{1.5}}
   \put(44,39){\circle*{1.5}}
   \put(49.5,39){\circle*{1.5}}
   \put(55,39){\circle*{1.5}}
   \put(60.5,39){\circle*{1.5}}
 \put(20.75,33.5){\line(1,0){2.5}}
   \put(44,33.5){\circle*{1.5}}
   \put(49.5,33.5){\circle*{1.5}}
   \put(55,33.5){\circle*{1.5}}
   \put(60.5,33.5){\circle*{1.5}}
   \put(66,33.5){\circle*{1.5}}
 \put(20.75,28){\line(1,0){2.5}}
   \put(44,28){\circle*{1.5}}
   \put(49.5,28){\circle*{1.5}}
   \put(55,28){\circle*{1.5}}
   \put(60.5,28){\circle*{1.5}}
   \put(66,28){\circle*{1.5}}
   \put(71.5,28){\circle*{1.5}}
 \put(20.75,22.5){\line(1,0){2.5}}
   \put(55,22.5){\circle*{1.5}}
   \put(60.5,22.5){\circle*{1.5}}
   \put(66,22.5){\circle*{1.5}}
   \put(71.5,22.5){\circle*{1.5}}
   \put(77,22.5){\circle*{1.5}}
 \put(20.75,17){\line(1,0){2.5}}
    \put(66,17){\circle*{1.5}}
    \put(71.5,17){\circle*{1.5}}
    \put(77,17){\circle*{1.5}}
    \put(82.5,17){\circle*{1.5}}
 \put(20.75,11.5){\line(1,0){2.5}}
    \put(71.5,11.5){\circle*{1.5}}
    \put(77,11.5){\circle*{1.5}}
    \put(82.5,11.5){\circle*{1.5}}
    \put(88,11.5){\circle*{1.5}}
 \put(27.5,4.75){\line(0,1){2.5}}
 \put(33,4.75){\line(0,1){2.5}}
 \put(38.5,4.75){\line(0,1){2.5}}
 \put(44,4.75){\line(0,1){2.5}}
 \put(49.5,4.75){\line(0,1){2.5}}
 \put(55,4.75){\line(0,1){2.5}}
 \put(60.5,4.75){\line(0,1){2.5}}
 \put(66,4.75){\line(0,1){2.5}}
 \put(71.5,4.75){\line(0,1){2.5}}
 \put(77,4.75){\line(0,1){2.5}}
 \put(82.5,4.75){\line(0,1){2.5}}
 \put(88,4.75){\line(0,1){2.5}}
   \put(88,6){\circle*{1.5}}
 \put(93.5,4.75){\line(0,1){2.5}}
   \put(93.5,6){\circle*{1.5}}
\end{picture}
\end{example}
While $\pi_k(T)$ is not a convex function over the range $0 \le k
\le c(T)$ in the last example, the calculated set $\cI(T)$ does at
least contain all of the lattice points in its own convex hull. To
expect this convexity to hold for every tree would be overly
optimistic, however: if we carry out the same calculation for the
larger tree $ \bigvsum\limits_{i=1}^5 S_4$ (on $16$ vertices instead
of $13$) we find that the points $(11, 1)$ and $(5, 5)$ both belong
to the inertia set, but their midpoint $(8, 3)$ does not.
\begin{question}
What is the computational complexity of determining the partial
inertia set of a tree? The examples above pose no difficulty, but
they do show that the greedy algorithm for $\MD{k}{T}$ fails even
for $T$ a tree. Computing all $n$ values of $\MD{k}{G}$ for a
general graph $G$ is NP-hard because it can be used to calculate the
independence number: $k + \MD{k}{G} = n$ if and only if there is an
independent set of size $n-k$.
\end{question}

In the next section we will consider more general graphs,
rather than restricting to trees and forests, and we will
see that even convexity of partial inertias within a
single stripe can fail in the broader setting.


\section{Beyond the forest}
\label{S7} In this section we investigate, over the set of all
graphs, what partial inertia sets---or more specifically, what
complements of partial inertia sets---can occur. Once a graph $G$ is
allowed to have cycles, we can no longer assume that $\hcI(G) =
\cI(G)$ by diagonal congruence. It happens, however, that each of
the results in this section is the same in the complex Hermitian
case as in the real symmetric case. For the two versions of each
question we will therefore demonstrate whichever is the more
difficult of the two, proving theorems over the complex numbers but
providing counterexamples over the reals.

It is convenient at this point to introduce a
way of representing
the complements of partial inertia sets.

\begin{definition}
A \textit{partition} is a finite (weakly) decreasing sequence of
positive integers. The first integer in the sequence is called the
\textit{width} of the partition, and the number of terms in the
sequence is called the \textit{height} of the partition.
\end{definition}

It is traditional to depict partitions with box diagrams.
In order to agree with our diagrams of partial inertia
sets, we choose the convention of putting the
longest row of boxes on the bottom of the stack;
for example, the decreasing sequence $(5,4,1)$ is
shown as the partition
\partboxes{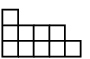}.
Given a box diagram of height $h$ and width $w$,
we index the rows by $0, 1, \ldots, h-1$ from
bottom to top and the columns by
$0, 1, \ldots, w-1$ from left to right.

\begin{definition}
Given a partition $\pi = (\pi_0, \pi_1, \ldots, \pi_{k-1})$,
let $\ell = \pi_0$, and for
$i \in \{0,1,\ldots,\ell-1\}$ let
$\pi^*_i = |\{j :\ \pi_j \ge i+1 \}|$,
i.e. the number of boxes in column $i$ of the
box diagram of $\pi$.
Then $\pi^* = (\pi^*_0, \pi^*_1, \ldots, \pi^*_{\ell-1})$
is called the conjugate partition of $\pi$.
A partition $\pi$ is \textit{symmetric} if $\pi = \pi^*$.
\end{definition}

For example, we have $(5,4,1)^* = (3,2,2,2,1)$
and $(3,3,2)^* = (3,3,2)$, so the partition
with box diagram
\partboxes{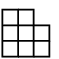} is symmetric.
It is easy to recognize symmetric partitions visually,
since a partition is symmetric if and only if
its box diagram has a diagonal axis of symmetry.

In this section we will describe $\cI(G)$,
for a graph $G$ on $n$ vertices, in terms of
its complement $\NN^2_{\leq n} \setminus \cI(G)$.
Definition~\ref{D19} gives us a natural way to describe
the shape of
this complement
as a partition.
We first extend to the Hermitian case
(distinguishing from the real symmetric case
as usual by prepending an `$\mathrm{h}$').
\begin{definition}
\label{D19b}
For any graph $ G $ on $ n $ vertices and $ k\in\{0,\ldots,n\} $,
let
\[
\hpi_k(G)=\min\left \{i :\ (i,k)\in\hcI(G)\right \}.
\]
\end{definition}

The partition corresponding to
a partial inertia set
is a list of as many of the values of $\pi_i(G)$
as are positive.

\begin{definition}
\label{Dparti}
Given a graph $G$, let $k = \pi_0(G)$
and let $h = \hpi_0(G)$.
Then
the \textit{inertial partition}
of $G$, denoted $\parti(G)$,
is the partition
\[
(\pi_0(G), \pi_1(G), \ldots, \pi_{k-1}(G)).
\]
The \textit{Hermitian inertial partition}
of $G$, denoted $\hpi(G)$,
is the partition
\[
(\hpi_0(G), \hpi_1(G), \ldots, \hpi_{h-1}(G)).
\]
\end{definition}

It would perhaps be more accurate to call these the
\textit{partial inertia complement partition} and
\textit{Hermitian partial inertia complement partition,}
but we opt for the abbreviated names.

The Northeast Lemma ensures that the
inertial partition and Hermitian inertial
partition of a graph are
in fact partitions, and by Observation~\ref{Ob8}
the partitions $\parti(G)$ and $\hpi(G)$
are always symmetric.
This symmetry is the reason why
$k=\pi_0(G)$ is the correct point of truncation:
$\pi_{k-1}(G) > 0$, but $\pi_k(G) = 0$.
\begin{remark}
\label{Rboxes}
If one starts with the entire first quadrant of the plane
$\RR^2$
and then removes everything ``northeast'' of any point
belonging to $\cI(G)$, the remaining
``southwest complement'' has the same shape
as the box diagram of $\pi(G)$.
The same applies of course to $\hcI(G)$
and $\hpi(G)$.
\end{remark}

The partial inertia sets $\cI(G)$ and $\hcI(G)$
can be reconstructed
from the partitions $\parti(G)$ and $\hpi(G)$,
respectively, if the number of
vertices of $G$ is also known.
The addition of an isolated vertex
to a graph $G$ does not change $\parti(G)$.

We begin to investigate the following problem:
\begin{question}
[Inertial Partition Classification Problem]
For which symmetric partitions $\pi$ does there exist
a graph $G$ for which $\parti(G) = \pi$?
\end{question}

Rather than examining all possible partial inertias for a
particular graph, we are now examining what restrictions
on partial inertias (or rather excluded partial inertias)
may hold over the class of all graphs.

The Hermitian Inertial Partition Classification Problem
is the same question with $\hpi(G)$ in the place of $\parti(G)$.
While it is known that there are graphs
$G$ for which $\cI(G)$ is a
strict subset of $\hcI(G)$,
it is not known whether there
are partitions $\pi$ that are inertial partitions but
not Hermitian inertial partitions, or vice versa.

At the moment we are only able to give a complete answer to the
Inertial Partition Classification Problem for symmetric partitions
of height no greater than~3.  We first list examples for a few
symmetric partitions that are easily obtained. Of course, adding an
isolated vertex to any example gives another example for the same
partition. For simplicity we will identify the partition $\pi(G)$
with its box diagram.

\begin{itemize}
\item For height 0, $\parti(G)$ is the empty partition if $G$ has no edges.
\item For height 1, $\parti(K_n) = \partboxes{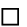}$ for any $n>1$.
\item For height 2, $\parti(P_3) = \partboxes{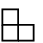}$.
\item For height 3, $\parti(S_4) = \partboxes{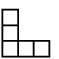}$ and $\parti(P_4) = \partboxes{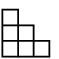}$.
\end{itemize}

The partitions already listed cover every possible case, up to height~$3$,
of an inertia-balanced graph, and leave three non-inertia-balanced
partitions unaccounted for:
\[
\partboxes{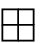},
\ \ \partboxes{332},
\mbox{ and }
\ \partboxes{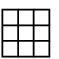}.
\]

The following theorem eliminates cases \partboxes{22} and \partboxes{333},
as well as every larger square partition.

\begin{theorem}
\label{Th_nosquares}
Let $G$ be a graph and let $M \in \Herm(G)$ be a
Hermitian matrix with partial inertia $(k,0)$, $k > 1$.
Then there exists a matrix $M^\prime \in \Herm(G)$
with partial inertia $(r,s)$ satisfying
$r < k$ and $s < k$.  Furthermore, if
$M$ is real then $M^\prime$ can be taken
to be real.
\end{theorem}

\begin{corollary}[No Square Partitions]
\label{Cnosquares}
For any $k > 1$, the
square partition $\pi = (k, k, \ldots, k)$
of height $k$ and width $k$ is not
the inertial partition of any graph $G$,
and is not the Hermitian inertial partition
of any graph $G$.
\end{corollary}

\begin{proof}[Proof of Theorem~\ref{Th_nosquares}]
Let $G$ be a graph on $n$ vertices
and suppose that $M \in \Herm(G)$
is a Hermitian matrix with partial inertia
$(k,0)$.  The matrix $M = \bigl[m_{ij}\bigr]$ is
thus positive semidefinite of rank $k$,
and can be factored as
$A^*A$ for some $k\times n$ complex matrix $A = \bigl[a_{ij}\bigr]$.
If $M$ is real symmetric, then $A$ can
be taken to be real.

We wish to construct a matrix $M^\prime \in \Herm(G)$
with strictly fewer than $k$ positive eigenvalues
and also strictly fewer than $k$ negative
eigenvalues.
By Proposition~\ref{P2}, we will have accomplished
our purpose if we can find $(k-1)\times n$ matrices
$B = \bigl[b_{ij}\bigr]$
and
$C = \bigl[c_{ij}\bigr]$
such that $B^*B-C^*C \in \Herm(G)$,
with the requirement that $B$ and $C$ be
real if $M$ is real.

We need to impose some mild general-position
requirements
on the first two rows of the matrix $A$,
which we accomplish by replacing $A$ by
$UA$, where $U$ is a unitary matrix
and where $U$ is real
(and hence orthogonal) in the case that
$A$ is real.
This is a permissible substitution because
$(UA)^*UA = A^*U^*UA = A^*A = M$.

The first general-position requirement is that, for integers $1 \le
j \le n$, $a_{1j} \neq 0$ and $a_{2j} \neq 0$ unless column $j$ of
$A$ is the zero column. The second requirement, which we will
justify more carefully, is that the set of ratios
$\{a_{1j}/a_{2j}\}$ be disjoint from the set of conjugate
reciprocals $\{\overline{a}_{2i}/\overline{a}_{1i}\}$, or
equivalently
\[
   \overline{a}_{1i}a_{1j} \neq \overline{a}_{2i}a_{2j}
\]
for any
$1 \leq i,j \leq n$ where neither $i$ nor $j$
corresponds to a zero column.

Now we prove the existence of a
unitary matrix $U$ with the desired
properties.  To do so, we
temporarily reserve
the symbol $i \in \CC$ to represent
a solution to $i^2 + 1 = 0$.
For the duration of this argument,
$j$ will represent any integer
$1 \le j \le n$ such that column $j$
of $A$ is not the zero column.

Let $x_j$ represent the vector
$\left [ \begin{array}{c} a_{1j} \\
a_{2j} \end{array} \right ]$.
If $\CC^*$ represents the set of
nonzero complex numbers,
then our first general position assumption
already guarantees $x_j \in (\CC^*)^2$.
We now define three functions
$z,\overline{z},w : (\CC^*)^2 \rightarrow \CC^*$
by
\[
z \left(
\left [ \begin{array}{c} p \\
q \end{array} \right ]
\right ) = p/q,
\ \
\overline{z} \left(
\left [ \begin{array}{c} p \\
q \end{array} \right ]
\right ) = \overline{p}/\overline{q},
\ \mbox{ and }\
w \left(
\left [ \begin{array}{c} p \\
q \end{array} \right ]
\right ) = \overline{q}/\overline{p}.
\]
Our task is to find a unitary matrix $U_1$,
orthogonal in the case that $A$ is real,
such that the sets $\{z(U_1x_j)\}$
and $\{w(U_1x_j)\}$ are disjoint.
In this case we can achieve the desired
general position of $A$ by replacing
it with $UA$, where
$U = U_1 \oplus I_{k-2}$.

Now consider the unitary matrices
\[
R_\theta =
  \left [
    \begin{array}{cc}
      e^{i \theta/2} & 0 \\
      0 & e^{-i \theta/2}
    \end{array}
  \right ]
\ \mbox{ and } \  Q = \frac{1}{\sqrt{2}}
  \left [
    \begin{array}{cc}
      1 & i \\
      i & 1
    \end{array}
  \right ].
\]
These matrices transform complex ratios as follows:
\[ z(R_\theta x) = e^{i \theta}z(x), \ \
   z(Q x) = \frac{z(x) + i}{i z(x) + 1} . \]
We have $Qx_j \in (\CC^*)^2$ as long as
$z(x_j) \not \in \{ i, -i \}$
and in particular as long as $z(x_j)$
is not pure imaginary,
which is automatically true in the case
$A$ is real.
In the case where $A$ is not real, we
can assume without loss of generality
that no $z(x_j)$ is pure imaginary
after uniformly multiplying on the left
by an appropriate choice of $R_\theta$.

Given $x$ and $y$ in $(\CC^*)^2$ such that
neither $z(x)$ nor $z(y)$ is pure imaginary,
$z(x) = w(y)$ if and only if
$z(Qx) = \overline{z}(Qy)$.
We have reduced the problem to that of
finding
a unitary matrix $U_1$, orthogonal
in the case $A$ is real, such that
the sets $\{ z(QU_1x_j) \}$
and $\{ \overline{z}(QU_1x_j) \}$
are disjoint.
In fact we will establish the stronger
condition that the two finite subsets
of the unit circle
\[
\left \{ \frac{z(QU_1x_j)}{|z(QU_1x_j)|} \right \} \ \mbox{ and }\
\left \{ \frac{\overline{z}(QU_1x_j)}{|\overline{z}(QU_1x_j)|}
\right \}
\]
are disjoint.
Let
\[
U_1 = Q^*R_\theta Q
=
\left [
  \begin{array}{cc}
    \cos \theta/2 & -\sin \theta/2 \\
    \sin \theta/2 &  \cos \theta/2
  \end{array}
\right ].
\]
Then $U_1$ is orthogonal,
and
\[
\frac{z(QU_1x_j)}{|z(QU_1x_j)|}
= e^{i \theta} \frac{z(Qx_j)}{|z(Qx_j)|}.
\]
Our general-position requirement for $A$ thus
reduces to the following fact: Given a finite
subset $P$ of the unit circle, there is some
$\theta$ such that $e^{i \theta}P$ is disjoint
from its set of conjugates $e^{-i \theta}\overline{P}$
and from the set $\{i,-i\}$.
To be concrete, if $\epsilon$ is the minimum
nonzero angle between any element of $P$
and any element of $\overline{P}$ or
$\{i,-i\}$, $\theta = \epsilon/3$
will suffice.  This concludes the argument
justifying our assumption of general position for $A$.

We now construct the matrices $B$ and $C$. Each column $j$ of the
matrices $B$ and $C$ for $1 \leq j \leq n$ is as follows:
\begin{itemize}
\item
$b_{1j} = a_{1j}^2$.
\item
$b_{ij} = a_{1j}a_{(i+1) j}$
for $i \in \{2, \ldots, k-1\}$.
\item
$c_{1j} = a_{2j}^2$.
\item
$c_{ij} = a_{2j}a_{(i+1) j}$
for $i \in \{2, \ldots, k-1\}$.
\end{itemize}

Now consider an arbitrary entry $m^\prime_{ij}$ of the matrix $M^\prime = B^*B - C^*C$;
this takes the form
\begin{eqnarray*}
m^\prime_{ij} & = & \overline{a}_{1i}^2 a_{1j}^2 +
\overline{a}_{1i}a_{1j}   \overline{a}_{3i} a_{3j}+ \cdots +
\overline{a}_{1i}a_{1j}   \overline{a}_{ki} a_{kj} \\
&& -\
 \overline{a}_{2i}^2 a_{2j}^2 -
 \overline{a}_{2i}a_{2j}   \overline{a}_{3i} a_{3j}- \cdots -
 \overline{a}_{2i}a_{2j}   \overline{a}_{ki} a_{kj},
\end{eqnarray*}
which factors as
\begin{eqnarray*}
m^\prime_{ij} & = &
(\overline{a}_{1i} a_{1j} - \overline{a}_{2i} a_{2j})
( \overline{a}_{1i} a_{1j} +  \overline{a}_{2i} a_{2j}
   +  \overline{a}_{3i} a_{3j} + \cdots +  \overline{a}_{ki} a_{kj}) \\
& = & (\overline{a}_{1i} a_{1j} - \overline{a}_{2i} a_{2j}) m_{ij}.
\end{eqnarray*}
In case either column $i$ or column $j$ of $A$ is the zero column,
we have $m^\prime_{ij} = 0 = m_{ij}$, and in all other cases we have,
by the generic requirement
\[
   \overline{a}_{1i}a_{1j} \neq \overline{a}_{2i}a_{2j},
\]
that $m^\prime_{ij} = 0$ if and only if $ m_{ij} = 0$.
It follows that $M^\prime$ is a matrix in $\Herm(G)$, and
by construction $M^\prime$ has at most $k-1$ positive
eigenvalues and at most $k-1$ negative eigenvalues.
Furthermore, if $M$ is real symmetric then
so is $M^\prime$.
\end{proof}

We have determined which inertial partitions occur for all
partitions up to height $3$ except for one: the partition
\partboxes{332}. Perhaps surprisingly, there is indeed a graph, on
$12$ vertices, that achieves this non-inertia-balanced partition in
both the real symmetric and Hermitian cases.

\begin{theorem}
\label{Th_exists332}
There exists a graph $G_{12}$ on $12$ vertices such that
$\parti(G_{12}) = \hpi(G_{12}) = \partboxes{332}$, the partition $(3,3,2)$.
\end{theorem}

The counterexample graph $G_{12}$ will be defined directly in terms
of a real symmetric matrix with partial inertia $(3,0)$; we will
then show that $(2,1)$ is not in $\hcI(G_{12})$.

Consider a cube centered at the origin of $\RR^3$,
and choose a representative vector for each line
that passes through an opposite pair of faces, edges,
or corners of the cube.
\[
\setlength{\unitlength}{1pt}
\begin{picture}(240,200)
 \put(0,0){\includegraphics{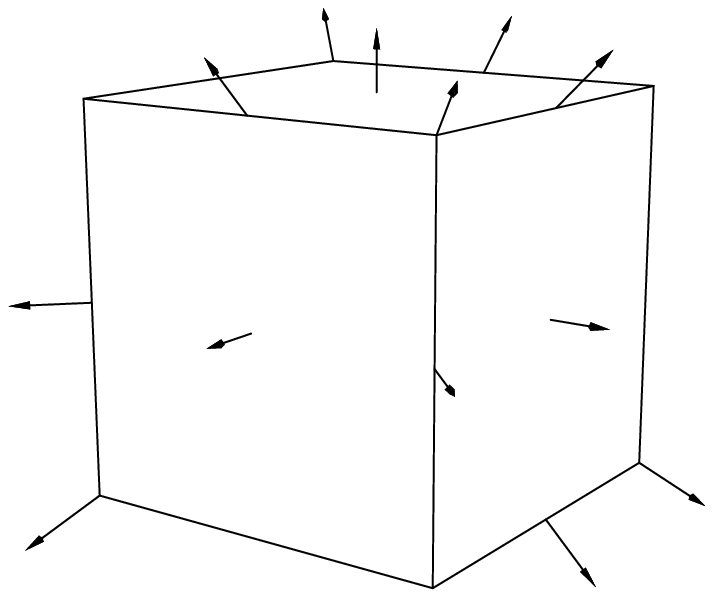}}

\put( 64 , 81){$ x $}
\put( 188.5 , 86.5){$ y $}
\put( 117.5 , 178){$ z $}
\put( 189 , 170){$ 1 $}
\put( 64 , 168){$ 2 $}
\put( 144 , 65.5){$ 3 $}
\put( 184 , 8){$ 4 $}
\put( 159 , 181){$ 5 $}
\put( 7 , 92){$ 6 $}
\put( 12 , 20){$ 7 $}
\put( 215.5 , 34){$ 8 $}
\put( 102 , 184){$ 9 $}
\put( 144 , 152){$ 10 $}

\end{picture}
\]

These 13 vectors give us
the columns of a matrix
\[
\raisebox{0pt}[1.2\totalheight]{
$
M_{13} =
\left[
\raisebox{.15\totalheight}[-.3\totalheight]{
$
\begin{array}{rrrrrrrrrrrrc}
  x & y & z & 1 & 2 & 3 & 4    & 5    &    6 & 7 &    8 &    9 & 10 \\
  1 & 0 & 0 & 0 & 1 & 1 & 0    & \!-1 &    1 & 1 & \!-1 & \!-1 & 1 \rule{0pt}{3ex} \\
  0 & 1 & 0 & 1 & 0 & 1 & 1    & 0    & \!-1 & \!-1 & 1 & \!-1 & 1 \\
  0 & 0 & 1 & 1 & 1 & 0 & \!-1 & 1    &    0 & \!-1 & \!-1 & 1 & 1
\end{array}
$
}
\right]
$
},
\]
which columns we index by the set of symbols
$\{x,y,z,1,\ldots,10\}$. The matrix $M_{13}^TM_{13}$ is real
symmetric and positive semidefinite of rank $3$, and thus has
partial inertia $(3,0)$.  We define $G_{13}$ as the graph on 13
vertices (labeled by the same 13 symbols) for which $M_{13}^TM_{13}
\in \Sym(G_{13})$; distinct vertices $i$ and $j$ of $G_{13}$ are
adjacent if and only if columns $i$ and $j$ of $M_{13}$ are not
orthogonal.  We note in passing that the subgraph of $G_{13}$
induced by vertices labeled 1--10 is the line graph of $K_5$, or the
complement of the Petersen Graph. We now define the graph $G_{12}$
(as promised in Theorem~\ref{Th_exists332}) as the induced subgraph
of $G_{13}$ obtained by deleting the vertex labeled $10$.

Before proving Theorem \ref{Th_exists332}, we prove a lemma
about a smaller graph $G_{10}$ that is obtained from $G_{12}$
by deleting the vertices labeled $6$ and $9$ (while retaining
the labels of the other vertices).
The vertices of $G_{10}$ are
thus labeled
$\{x,y,z,1,2,3,4,5,7,8\}$ (notice that this set skips index $6$).

\begin{lemma}
\label{LG10}
Let $A = \bigl [ a_{ij} \bigr ]$ be a Hermitian matrix in $\Herm(G_{10})$
of rank no more than $3$.  Then the first two
diagonal entries $a_{xx}$ and $a_{yy}$
are both nonzero and have the same sign.
\end{lemma}

\begin{proof}
Let $d_x$, $d_y$, and $d_z$ be the first three diagonal entries of $A$:
\[d_x=a_{xx},\ \  d_y=a_{yy},\ \  d_z=a_{zz}.\]  We show first that all
three of these entries are nonzero.  For this purpose
it suffices to consider only the first six rows and columns of $A$,
corresponding to the graph $G_6 = $
\[
\setlength{\unitlength}{1pt}
\begin{picture}(68,59)
 \put(0,0){\includegraphics{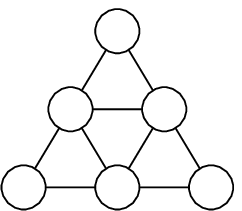}}
 \put(31,49){$z$}
 \put(18,26){2}
 \put(44.5,26){1}
 \put(4,5){$x$}
 \put(30.5,3.5){3}
 \put(58,5){$y$}
\end{picture}
\]
sometimes called the supertriangle graph.  The
automorphism group of $G_6$
realizes any permutation
of the vertices $\{x,y,z\}$ (as do the automorphism
groups of $G_{13}$ and of $G_{12}$, but
not that of $G_{10}$).

The principal submatrix of $A$
on rows and columns $\{x,y,z,1,2,3\}$, like $A$ itself, has rank
at most $3$.  Suppose that we had $d_y = 0$ while $d_z \neq 0$.
Then the $4\times 4$ submatrix on rows $\{x,y,z,1\}$ and
columns $\{y, z,1,2\}$ would be combinatorially nonsingular
(that is, permutation equivalent to an upper-triangular matrix
with nonzero entries on the diagonal), contradicting
that $\rank(A) \leq 3$.
By the symmetries of $G_6$, we could have chosen
any pair instead of $d_y=0,d_z \neq 0$, and thus if any one of the
three quantities $d_x$, $d_y$, or $d_z$ is equal to zero,
then all three must be.
But if all three of the first diagonal entries were zero,
then the $4\times 4$ principal submatrix
on rows and columns $\{x,y,1,2\}$
would be combinatorially nonsingular.
It follows that in the
$10\times 10$ rank-3 matrix $A$, the first three diagonal
entries $d_x$, $d_y$, and $d_z$ are all nonzero.

Considering once more the full matrix $A$,
let
$\beta = a_{x2}/a_{z2}$
and
$\gamma = a_{y1}/a_{z1}$,
so the first three rows of $A$ can be written:
\[
\left[
\begin{array}{ccccccccccc}
  d_x &  0  &  0  &    0    &\beta a_{z2} & a_{x3} &    0   & a_{x5} & a_{x7} & a_{x8} \\
   0  & d_y &  0  &\gamma a_{z1} &    0   & a_{y3} & a_{y4} &    0   & a_{y7} & a_{y8} \\
   0  &  0  & d_z & a_{z1} & a_{z2}       &    0   & a_{z4} & a_{z5} & a_{z7} & a_{z8}
\end{array}
\right].
\]
Since $d_x$, $d_y$, and $d_z$ are nonzero and $A$
has rank at most $3$,
every other row of $A$
can be obtained from these first three rows
by taking a linear combination,
and the coefficients
of the linear combination are determined
by entries in
the first three columns.
Every entry of $A$ is thus
determined by
the variables appearing
in the $3 \times 10$ matrix above.
For any $i$ and $j$ in the set $\{1,2,3,4,5,7,8\}$,
we have
\[
a_{ij}
  = \frac{\overline{a}_{xi}}{d_x} a_{xj}
  + \frac{\overline{a}_{yi}}{d_y} a_{yj}
  + \frac{\overline{a}_{zi}}{d_z} a_{zj}.
\]

In those cases where $i \neq j$ and $ij$
is not an edge of $G_{10}$, the entry
$a_{ij} = 0$
gives an equation on the entries of the first three rows.
Using several such equations, we deduce
that $d_xd_y > 0$, as follows:
\begin{enumerate}
\item
  The entries $a_{27} = 0$ and $a_{18}=0$ give us the pair of equations
  \[
  d_x a_{z7} = -d_z \overline{\beta} a_{x7}
  \mbox{\ and\ }
  d_y a_{z8} = -d_z \overline{\gamma} a_{y8}.
  \]
\item
  Combining the equations from $a_{37} = 0$ and $a_{38} = 0$, we have
  \[
  a_{x7} a_{y8} = a_{y7} a_{x8}.
  \]
\item
  Combining the equations from $a_{52} = 0$ and $a_{58} = 0$, we have
  \[
  a_{x8} = \beta a_{z8}.
  \]
\item
  Combining the equations from $a_{41} = 0$ and $a_{47}=0$, we have
  \[
  a_{y7} = \gamma a_{z7}.
  \]
\end{enumerate}
Multiplying the first pair of equations and then substituting in
each of the remaining equations in order, then canceling the nonzero
term $a_{z7} a_{z8}$, we arrive finally at
  \[
  d_x d_y = d_z^2\, \beta \overline{\beta}\, \gamma \overline{\gamma},
  \]
a positive quantity.  This proves that, in any Hermitian matrix
$A \in \Herm(G_{10})$ of rank no more than $3$,
the first two diagonal entries are nonzero and have the same sign.
\end{proof}

\begin{proof}[Proof of Theorem~\ref{Th_exists332}]
We review the definition of the graph
$G_{12}$ that will provide the claimed example:
starting from the diagram of the cube
with a labeled vector for every pair of faces, edges, or corners, we omit
the vector~$10$ and connect pairs of vertices from the set
$\{x,y,z,1,2,3,4,5,6,7,8,9\}$ whenever their corresponding
vectors are not orthogonal.

Letting $M_{12}$ be the submatrix of $M_{13}$ obtained by deleting
the last column (labeled $10$), we have $M_{12}^TM_{12} \in
\Sym(G_{12})$ and thus $(3,0) \in \cI(G_{12})$ and $(3,0) \in
\hcI(G_{12})$. It follows from Theorem~\ref{Th_nosquares} that the
point $(2,2)$ also belongs to $\cI(G_{12})$ and $\hcI(G_{12})$. To
show that $\parti(G_{12}) = \hpi(G_{12}) = \partboxes{332}$, it
suffices by the Northeast Lemma and symmetry to show that $(2,1)
\not \in \hcI(G_{12})$.

Let $A$ be any matrix of rank $3$ in $\Herm(G_{12})$,
and let $d_x$, $d_y$, and $d_z$ be the first three diagonal
entries of $A$.
Omitting rows and columns $6$ and $9$ gives us a matrix of
rank no more than $3$ in $\Herm(G_{10})$, and so Lemma~\ref{LG10}
tells us that
$d_x$ and $d_y$ are nonzero and have the same sign.
However, the automorphism
group of $G_{12}$ inherits all the symmetries of a cube with
one marked corner, and thus anything true of the pair
of vertices $\{x,y\}$
is also true of the pair $\{y,z\}$,
so
$d_y$ and $d_z$ are also nonzero and have the same sign.
More explicitly, using the symmetry of a counterclockwise rotation of
the cube around the corner marked $10$,
we delete rows and columns $4$ and $7$
(instead of $6$ and $9$)
and reorder the remaining rows and columns as
$(y,z,x,2,3,1,5,6,8,9)$
to yield a different matrix belonging to
$\Herm(G_{10})$, and invoke Lemma~\ref{LG10}
again to obtain $d_yd_z > 0$, showing that
the three diagonal entries $d_x$, $d_y$, and $d_z$
all have the same sign.

The principal submatrix of $A$ on rows and columns $\{x,y,z\}$
has either three positive eigenvalues or three negative eigenvalues,
and so by interlacing the partial inertia of $A$
must be either $(3,0)$ or $(0,3)$.  This completes the
proof that the graph $G_{12}$ achieves the inertial partition
and Hermitian inertial partition $(3,3,2)$, and is not
inertia-balanced.
\end{proof}

Of course the same argument also shows that
$\parti(G_{13}) = \hpi(G_{13}) = \partboxes{332}$,
but we are interested in the smallest possible
graph that is not inertia-balanced.
The following proposition justifies our claim
that $G_{12}$ is at least locally optimal.

\begin{theorem}
\label{PG12min}
Every proper induced subgraph of $G_{13}$
is either isomorphic to $G_{12}$ or
is inertia-balanced and Hermitian inertia-balanced.
\end{theorem}

\begin{proof}
By Theorem~\ref{Th_nosquares} every graph $G$ with
$\mr(G) < 3$ is inertia-balanced and every
graph $G$ with $\hmr(G) < 3$ is
Hermitian inertia-balanced.
It thus suffices to show that
for every proper induced subgraph
$F$ of $G_{13}$ other than $G_{12}$,
$(2,1) \in \cI(F)$
unless $|F| < 3$.

Recall that $G_{13}$ is defined by orthogonality relations between
the columns of the matrix
\[
M_{13}\! = \!
\left[
\begin{array}{rrrrrrrrrrrrc}
  1 & 0 & 0 & 0 & 1 & 1 & 0    & \!\!-1 &    1 & 1 & \!\!-1 & \!\!-1 & 1 \\
  0 & 1 & 0 & 1 & 0 & 1 & 1    & 0    & \!\!-1 & \!\!-1 & 1 & \!\!-1 & 1 \\
  0 & 0 & 1 & 1 & 1 & 0 & \!\!-1 & 1    &    0 & \!\!-1 & \!\!-1 & 1 & 1
\end{array}
\right],
\]
corresponding to various axes of symmetry of a cube.
In other words,
$M_{13}^T I_3 M_{13} \in \Sym(G_{13})$ (where the identity matrix $I_3$
imposes the standard positive definite inner product on $\RR^3$)
and so $(3,0) \in \cI(G_{13})$.

The automorphism group of $G_{13}$ has three orbits, corresponding
to the faces ($x$, $y$, and $z$), edges (1, 2, 3, 4, 5, and 6), and
corners (7, 8, 9, and 10) of the cube. The deletion of any corner
yields $G_{12}$ (perhaps with a different labeling) and there is, up
to isomorphism, only one way to delete two corners.  Every proper
induced subgraph of $G_{13}$ other than $G_{12}$ is thus isomorphic
to an induced subgraph of $G_{13} - x$, of $G_{13} - 3$, or of
$G_{13} - \{7,8\}$. Letting $G$ be each of these three graphs in
turn, we exhibit for each a diagonal matrix $D$ with $\pin(D) =
(2,1)$ and a real matrix $M$ such that $M^{\,T}\!DM \in \Sym(G)$.

\bigskip

\noindent $G = G_{13} - x,\  (2,1) \in \cI(G)$:
\[
D =
\left[
\begin{array}{rrr}
  3 & 0 & 0 \\
  0 & 1 & 0 \\
  0 & 0 &\!\!-2
\end{array}
\right],\ \
M =
\left[
\begin{array}{rrrrrrrrrrrrr}
 . & 0 & 2 & \!\!-1 & 1 &  1 & \!\!-1 & 0 & 1 &  0 & 1 & 0 &  1 \\
 . & 1 & 0 &  1 & 0 & \!\!-1 & \!\!-1 & 0 & 1 & \!\!-2 & 3 & 2 & \!\!-3 \\
 . & 0 & 3 &  1 & 0 &  1 &  1 & 1 & 1 &  1 & 0 & 1 &  0
\end{array}
\right]
\]

\bigskip

\noindent $G = G_{13} - 3,\  (2,1) \in \cI(G)$:
\[
D =
\left[
\begin{array}{rrr}
  1 & 0 & 0 \\
  0 & 1 & 0 \\
  0 & 0 &\!\!-1
\end{array}
\right],\ \
M =
\left[
\begin{array}{rrrrrrrrrrrrr}
 1 & 0 & 0 & 0 & 2 & . & 0 & 1 & \!\!-1 & 1 & 4 & 1 & 2 \\
 0 & 1 & 0 & 2 & 0 & . & 1 & 0 &  1 & 4 & 1 & 1 & 2 \\
 0 & 0 & 1 & 1 & 1 & . & 2 & 2 &  0 & 2 & 2 & 2 & 1
\end{array}
\right]
\]

\bigskip

\noindent $G = G_{13} - \{7,8\},\  (2,1) \in \cI(G)$:
\[
D =
\left[
\begin{array}{rrr}
  1 & 0 & 0 \\
  0 & 1 & 0 \\
  0 & 0 &\!\!-1
\end{array}
\right],\ \
M =
\left[
\begin{array}{rrrrrrrrrrrrr}
 1 & 0 & 0 & 0 & 2 & 1 & 0 & 1 & \!\!-1 & . & . & 1 & 2 \\
 0 & 1 & 0 & 2 & 0 & 1 & 1 & 0 &  1 & . & . & 1 & 2 \\
 0 & 0 & 1 & 1 & 1 & 0 & 2 & 2 &  0 & . & . & 2 & 1
\end{array}
\right]
\]

\bigskip

For each value of $M$ and $D$, the matrix $M^{\,T}\!DM$ has
partial inertia $(2,1)$ and belongs to $\Sym(G)$
for the desired subgraph $G$.
If $F$ is a proper induced subgraph of $G_{13}$
other than $G_{12}$, then $F$ is an induced subgraph
of one of these three graphs.
Part (a) of Proposition~\ref{P22}
allows us to delete vertices
from any one of the three graphs and keep the partial
inertia $(2,1)$ as long as at least $3$ vertices
remain, which gives us
$(2,1) \in \cI(F)$ unless $|F| < 3$.
\end{proof}

\begin{question}
Is $G_{12}$ the unique graph on fewer than $13$ vertices
that is not inertia-balanced?
\end{question}

Theorems~\ref{Th_nosquares} and~\ref{Th_exists332} only permit
us to answer the Inertial Partition Classification Problem
and Hermitian Inertial Partition Classification Problem
up to height $3$.
We have shown examples of constructing a graph
whose minimum rank realization with a particular
partial inertia is sufficiently ``rigid''
to prevent intermediate partial inertias
of the same rank
between the matrix and its negative.
If the rank is allowed to increase, though,
it is much less clear what
restrictions can be made.
The next difficult
question appears to be whether
$(4,4,4,3) = \partboxes{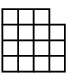}$ is an inertial partition.

\begin{question}
Let $G$ be a graph and let $M$ be a matrix
in $\Sym(G)$ with $\pin(M) = (4,0)$.  Must there
exist a matrix $M^\prime \in \Sym(G)$ with
$\pin(M^\prime) = (3,2)$?
\end{question}

On the one hand, partial inertia $(3,2)$ is of higher rank
than partial inertia $(4,0)$, which means that any proof
along the lines of Theorem~\ref{Th_exists332}---a proof
that a particular arrangement of orthogonality relations
of vectors in $\RR^4$ could not be duplicated
in $\RR^5$ with an indefinite inner product of
signature $(3,2)$---would
have an extra degree of freedom to contend with.
On the other hand, there
seems to be little hope of constructing the matrix
$M^\prime$ directly from $M$
using any sort of continuous map such as that
employed in the proof of Theorem~\ref{Th_nosquares}.

\end{document}